\pdfoutput=1
\documentclass[11pt]{amsart}

\usepackage[letterpaper,margin=1in]{geometry}
\usepackage{amsaddr}
\usepackage[T1]{fontenc}
\usepackage{lmodern}
\usepackage{newpxtext}
\usepackage{newpxmath}
\DeclareEncodingSubset{TS1}{zpltlf}{0}
\DeclareEncodingSubset{TS1}{zpllf}{0}
\usepackage{microtype}
\usepackage{graphicx}
\usepackage[inline]{enumitem}
\usepackage{booktabs}
\usepackage{cite}
\usepackage{algorithm}
\usepackage{algpseudocode}
\usepackage{xcolor}
\usepackage{siunitx}
\definecolor{tuftsblue}{HTML}{3E8EDE}
\usepackage[colorlinks,allcolors=tuftsblue]{hyperref}
\usepackage{cleveref}

\crefformat{equation}{(#2#1#3)}
\Crefformat{equation}{(#2#1#3)}
\crefrangeformat{equation}{(#3#1#4)--(#5#2#6)}
\Crefrangeformat{equation}{(#3#1#4)--(#5#2#6)}
\crefmultiformat{equation}{(#2#1#3)}{ and~(#2#1#3)}{, (#2#1#3)}{, and~(#2#1#3)}
\Crefmultiformat{equation}{(#2#1#3)}{ and~(#2#1#3)}{, (#2#1#3)}{, and~(#2#1#3)}

\makeatletter
\renewcommand{\@captionfont}{\normalfont\small}
\renewcommand\floatc@ruled[2]{{\small{\@fs@cfont #1} #2\par}}
\def\theHALG@line{\thealgorithm.\arabic{ALG@line}}
\makeatother

\theoremstyle{plain}
\newtheorem{lemma}{Lemma}
\newtheorem{proposition}[lemma]{Proposition}
\theoremstyle{remark}
\newtheorem{remark}{Remark}
\newtheorem{assumption}{Assumption}
\crefname{assumption}{Assumption}{Assumptions}
\crefname{appendix}{Appendix}{Appendices}
\crefname{subappendix}{Appendix}{Appendices}

\DeclareMathOperator*{\minopt}{minimize}
\newcommand{\scr}[1]{\ensuremath{\mathsf{#1}}}
\newcommand{\defeq}{\ensuremath{\mathrel{\mathop:}=}}
\newcommand{\eqdef}{\ensuremath{=\mathrel{\mathop:}}}

\newcommand{\tabadjust}{\centering\footnotesize\renewcommand{\arraystretch}{0.8}}
\newcommand{\algadjust}{\small}
\newcommand\fnum[1]{\num[scientific-notation=fixed,round-precision=3,round-mode=places]{#1}}
\newcommand\snum[1]{\num[round-precision=3,round-mode=places,tight-spacing=true,retain-explicit-plus=true,output-exponent-marker=\text{e}]{#1}}
\newcommand\inum[1]{\num[scientific-notation=fixed,round-precision=0,round-mode=places,group-separator={,}]{#1}}

\begin{document}

\title[Diffeomorphic Flows in Classification and Generative Modeling]{Optimal Control of RKHS Diffeomorphic Flows for Classification and Generative Modeling}

\author{Christopher Chukwuemeka}
\address{Department of Mathematics, Tufts University, Medford, MA, USA}
\email{christopher.chukwuemeka@tufts.edu}

\author{Andreas Mang}
\address{Department of Mathematics, Tufts University, Medford, MA, USA}
\email{andreas.mang@tufts.edu}



\begin{abstract}
We develop a unified optimal-control framework that solves supervised classification and likelihood-based density estimation with a single diffeomorphic flow whose velocity lives in a reproducing kernel Hilbert space (RKHS); the two tasks share all structure except the objective. Variable projection eliminates the terminal model (a classifier, or a Gaussian-mixture (GMM) reference), leaving a reduced control problem whose gradient comes from an analytic discrete adjoint. The RKHS parameterization makes the flow's divergence exact and deterministic at the cost of one kernel matrix--vector product, avoiding the stochastic trace estimators and potential-flow restrictions of existing continuous normalizing flows. Because the objective and gradient are deterministic, matrix-free reduced Gauss--Newton (GN) and related second-order solvers are practical for both tasks, which we benchmark against a range of first- and second-order optimizers; the learned map is a diffeomorphism with explicitly controlled regularity.

In experiments from dimension 2 to 256, the reduced GN solvers use the fewest iterations of any optimizer tested and, together with joint L-BFGS, the least work; a multiscale kernel improves the density fit, and a trainable GMM reference reduces the deformation needed for multimodal transport at comparable likelihood while lowering the likelihood on real data. The full-model solve is quadratic in the particle count and sub-linear in dimension; a control-point reduction makes the particle-count scaling linear at near-lossless accuracy.
\end{abstract}

\maketitle

\section{Introduction}
\label{sec:intro}

Diffeomorphic learning~\cite{younes2020:diffeomorphic} transforms data by a smooth invertible map $\phi: \mathcal{X} \to \mathcal{Z}$ (a diffeomorphism, i.e., a bijection with a smooth inverse) before a terminal model is applied. Untangling the data distribution into a more tractable representation that preserves a one-to-one correspondence with, and often the geometry of, the feature space $\mathcal{X}$, our construction serves two tasks:
\begin{enumerate*}[label={\it(\roman*)}]
\item for \emph{supervised classification} the transformed features become linearly separable and
\item for \emph{generative modeling} the transformed data are made to match a simple reference distribution from which one can both sample and evaluate likelihoods\end{enumerate*}.
The proposed framework bridges ideas from flows of diffeomorphisms~\cite{trouve1998:diffeomorphisms, dupuis1998:variational, beg2005:computing, glaunes2008:large, mang2023:operator, zhang2021:diffeomorphic, mang2019:claire, mang2017:lagrangian}, non-parametric prediction~\cite{ruiz2023:neural, haber2017:stable, chen2018:neural, rosso2025:weight}, and density estimation and variational inference based on normalizing flows (NFs)~\cite{papamakarios2021:normalizing, rezende2016:variational, marzouk2024:distribution, parno2018:transport,liu2016:stein}.

\begin{figure}
\centering
\includegraphics[height=3.2cm]{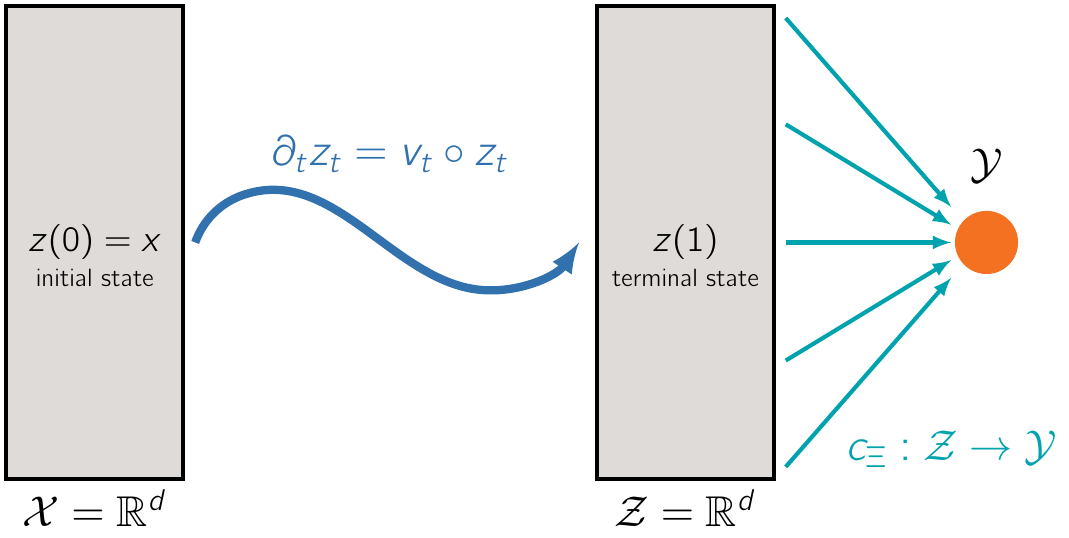}\qquad
\includegraphics[height=3.2cm]{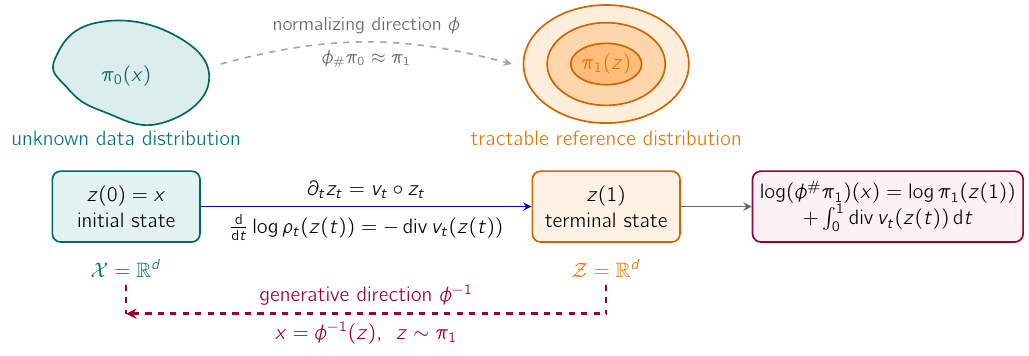}
\caption{Optimal control of diffeomorphic flows for classification and density estimation. In both tasks the features $x = z(0)$ are transported by the flow~\cref{eq:flow} to the terminal state $z(1)$. Left, classification: $z(1)$ is passed through a trainable classifier $c_{\Xi} : \mathcal{Z} \to \mathcal{Y}$. Right, density estimation: the flow pushes $\pi_0$ toward a reference $\pi_1$; the log-likelihood accumulates the divergence of $v_t$ along the trajectory (with $\rho_t \defeq (z_t)_\#\pi_0$ the transported density), and samples are generated by integrating backward, $x = \phi^{-1}(z)$ with $z \sim \pi_1$.}
\label{fig:framework}
\end{figure}

We frame both tasks as an optimal control problem~\cite{ruiz2023:neural, haber2017:stable, chen2018:neural}; \Cref{fig:framework} illustrates the two pipelines. The common input is a set of $n_s \in \mathbb{N}$ feature samples $x = (x_i)_{i = 1}^{n_s} \subset \mathcal{X}$, $\mathcal{X} \subseteq \mathbb{R}^d$, drawn from a data distribution $\pi_0$; the flow maps them into $\mathcal{Z} = \mathbb{R}^d$. In \emph{supervised classification}, each sample carries a label $y = (y_i)_{i = 1}^{n_s} \subset \mathcal{Y}$, where $\mathcal{Y}$ is an $n_c$-dimensional real-valued vector space or a probability simplex over $n_c$ classes, and the goal is to predict $y$ from $x$. In \emph{generative modeling} no labels are available; the goal is to estimate $\pi_0$ by transporting it onto a tractable reference $\pi_1$, so that one can evaluate likelihoods and draw samples. Our formulation seeks a diffeomorphic map $\phi$ that minimizes a task-dependent objective: a classification loss against $y$, or the likelihood under $\pi_1$.

Our architecture extends residual networks~\cite{he2016:deep, zagoruyko2016:wide}, interpreted as dynamical systems (\emph{neural} ordinary differential equations, ODEs~\cite{e2017:proposal, haber2017:stable}), to continuous time, where after discretization the layer count equals the number of ODE right-hand-side evaluations~\cite{haber2017:stable, chang2017:multi}. Following~\cite{younes2020:diffeomorphic}, we realize $\phi$ as the time-one flow of an ODE, modeling the transported features as $\mathbb{R}^d$-diffeomorphisms $z_t \defeq z(t,\cdot\,)$, $t \mapsto z_t$, given by the solution of
\begin{equation}\label{eq:flow}
\partial_t z_t = v_t \circ z_t \;\; \text{in}\; (0,1],
\quad z_t = \operatorname{id}_{\mathcal{X}} \;\text{for}\;t=0,
\end{equation}

\noindent where $v \in \mathcal{V}$, $\mathcal{V} \defeq L^2([0,1],\mathcal{H})$, $t \mapsto v_t \defeq v(t,\cdot\,)$, $v_t \in \mathcal{H}$, is a smooth, time-dependent vector field that will enter our problem formulation as a control variable, $\mathcal{H}$ is an adequately chosen Hilbert space, and $\operatorname{id}_{\mathcal{X}} : \mathcal{X} \to \mathcal{X}$, $\operatorname{id}_{\mathcal{X}}(x) = x$, is the identity transformation in $\mathcal{X}$.

We impose smoothness on $v$ so the flow is well posed. Under Carath\'eodory conditions ($t\mapsto v_t$ measurable with $\int_0^1\operatorname{lip}(v_t)\,\mathrm{d}t<\infty$ and $\int_0^1\|v_t(x_0)\|\,\mathrm{d}t<\infty$ for some $x_0$, both implied below by $\mathcal{H}\hookrightarrow C^1_0$), the ODE~\cref{eq:flow} has a unique solution~\cite{coddington1956:theory, butcher2016:numerical}, and the induced map $z_t$ is a bi-Lipschitz \emph{homeomorphism} with $\operatorname{lip}(z_t^{\pm1})\le\exp\int_0^1\operatorname{lip}(v_s)\,\mathrm{d}s$~\cite{calcaterra2008:lipschitz}; differentiability of the flow requires $v_t\in C^1(\mathcal{X},\mathbb{R}^d)$. The Lipschitz class alone bounds neither this exponent nor the inverse's differentiability; confining $v$ to $\mathcal{H}$ does both: each $z_t$ is a $C^q$ diffeomorphism, and the embedding $\mathcal{H}\hookrightarrow C^1_0$ together with Cauchy--Schwarz turns the geodesic-energy penalty into the bound $\int_0^1\operatorname{lip}(v_s)\,\mathrm{d}s \le c\,(\int_0^1\|v_s\|_{\mathcal{H}}^2\,\mathrm{d}s)^{1/2}$. Thus, the bi-Lipschitz constant, and by that the conditioning of the inverse, stays controlled over the fixed unit horizon. The learning task is then fixed by the objective on the terminal state $z(1)$. For classification, we feed $z(1)$ through a trainable classifier $c_{\Xi} : \mathcal{Z} \to \mathcal{Y}$ parameterized by $\Xi$. Training then seeks controls $v \in \mathcal{V}$, states $z \in \mathcal{T}$ (trajectories $t \mapsto z_t$ of maps $z_t : \mathcal{X} \to \mathcal{Z}$), and classifier parameters $\Xi \in \mathcal{W}$ solving
\begin{equation}\label{eq:varopt-training}
\begin{aligned}
\minopt_{v \in \mathcal{V},\; z \in \mathcal{T}, \;\Xi \in \mathcal{W}} \quad &
\operatorname{loss}(c_{\Xi}(z(1)), y)
+ \operatorname{dist}^2(\operatorname{id}_{\mathcal{X}}, z(1))
+ \operatorname{reg}(\Xi) \\
\text{subject to} \qquad & z\,\,\text{satisfies \cref{eq:flow}}.
\end{aligned}
\end{equation}

\noindent The three terms are a measure for the discrepancy between predictions $c_{\Xi}(z(1))$ and labels $y \in \mathcal{Y}$, a squared Riemannian distance that regularizes the transport map, and a convex penalty on $\Xi$. We will specify these building blocks precisely in \Cref{sec:methods}.

The same mathematical formulation describes generative modeling: replacing the supervised loss by a likelihood matching the terminal state $z(1)$ to a tractable reference density turns \cref{eq:varopt-training} into continuous-time density estimation.

Confining $v$ to a reproducing kernel Hilbert space (RKHS) and penalizing the transport by its geodesic energy makes the learned map a diffeomorphism with explicitly controlled regularity.

\subsection{Related Work}\label{sec:related-work}

Our work builds on the diffeomorphic learning framework of~\cite{younes2020:diffeomorphic}, which formulates supervised learning as a control problem over RKHS-parameterized flows of diffeomorphisms. This work has recently been extended to multivariate regression via affine--diffeomorphic sequences~\cite{lohr2024:finemorphs}. We extend it to likelihood-based density estimation and change how it is solved: where~\cite{younes2020:diffeomorphic} optimizes the flow controls and terminal classifier jointly with first-order methods, we eliminate the terminal block by variable projection (VarPro) and solve the reduced control problem with second-order Gauss--Newton (GN) solvers. We also introduce a sparse, reduced parameterization and improve expressivity.

Neural ODEs~\cite{e2017:proposal,haber2017:stable,chen2018:neural} provide a continuous-depth view of deep neural networks (NNs), with the forward pass defined by an ODE and trained through adjoint sensitivities. Continuous normalizing flows (CNFs)~\cite{chen2018:neural} extend this to density estimation, expressing the log-density evolution through the divergence of a vector field $v$. The original formulation computes the Jacobian trace explicitly at $\mathcal{O}(d^2)$ cost, which~\cite{grathwohl2018:ffjord} reduces to $\mathcal{O}(d)$ with Hutchinson's stochastic trace estimator, enabling unconstrained NNs but introducing variance and often stiff dynamics. Subsequent work counteracts this stiffness with optimal-transport-inspired penalties~\cite{finlay2020:how} or a potential-flow structure allowing exact, linear-cost divergence evaluation~\cite{onken2021:ot}. Other works that formulate training of deep NNs as optimal control problems are~\cite{benning2019:control, haber2017:stable, ruiz2023:neural}; \cite{benning2019:control} derives optimality conditions with structure-preserving discretizations. We adopt this discretize-then-optimize viewpoint. Our work also shares features with structure-preserving invertible or equivariant networks~\cite{celledoni2021:structure}.

Related generative approaches (flow matching~\cite{lipman2023:flow}, rectified flow~\cite{liu2023:rectified}, and stochastic interpolants~\cite{albergo2023:interpolants}) learn a velocity field $v$ by regressing against a prescribed probability path between data and reference, thereby avoiding likelihood-based training, with optimal-transport couplings straightening the trajectories~\cite{tong2024:improving}. Score-based diffusion models~\cite{ho2020:ddpm, song2021:score} learn the score of a noising process, whose probability-flow ODE gives a deterministic, likelihood-bearing counterpart. Recent analyses link bi-Lipschitz NF expressivity to diffusion through the induced diffeomorphic maps~\cite{iske2026:bilipschitz}. Flow matching has been adapted to identifying dynamical systems~\cite{rout2025:flow}. These methods emphasize scalable training and flexible neural parameterizations but typically sacrifice exact likelihood evaluation during training.

Close to our setting are kernel-parameterized diffeomorphic ODE flows. The work in~\cite{raviola2022:kernelode} casts kernel neural ODEs as optimal control over an RKHS parameterization of a velocity for supervised regression; \cite{pandey2025:diffeomorphic} learns such a flow for generative modeling, transporting a Gaussian reference to the data with a-priori error bounds. The latter work minimizes a sample discrepancy (maximum mean discrepancy, or kernel Stein discrepancy~\cite{fisher2021:ksd}) between pushforward and target, avoiding density evaluation. We instead retain a likelihood objective, evaluated through the closed-form RKHS divergence and minimized with an analytic discrete adjoint and reduced second-order GN solvers.

Our framework departs from existing CNF approaches in the divergence it uses, the objective it retains, and how it is trained. The RKHS parameterization makes the divergence exact and deterministic at the cost of one kernel matvec, avoiding the variance of stochastic trace estimators~\cite{grathwohl2018:ffjord} and potential-flow restrictions on $v$~\cite{onken2021:ot, finlay2020:how}. Unlike the simulation-free training of flow matching and score-based diffusion~\cite{lipman2023:flow, song2021:score}, we retain a likelihood objective with a reduced gradient from the flow's discrete adjoint rather than differentiation through a black-box network~\cite{ruiz2023:neural}. Being deterministic, our objective and gradient are amenable to reduced GN solvers, whereas most flow-based approaches optimize a stochastic likelihood estimate, and even second-order neural-ODE optimizers apply their curvature on top of a stochastic objective~\cite{liu2021:snopt}. Confining $v$ to an RKHS and penalizing its geodesic energy makes the learned transport a diffeomorphism with controlled smoothness, which discrete NFs obtain only through constrained invertible architectures~\cite{behrmann2019:iresnet, chen2019:residual} and unconstrained neural flows and diffusion models do not provide.

\subsection{Contributions}\label{sec:contributions}

We propose a unified optimal-control framework for learning and generative modeling based on diffeomorphic flows. Our contributions are:
\begin{itemize}[leftmargin=1em]
\item
\emph{Unified formulation.} We cast supervised classification and likelihood-based density estimation as optimal control problems over RKHS-parameterized diffeomorphic flows. The two tasks differ only in their objective (a terminal classification loss, or a likelihood that also accumulates the flow's divergence along the trajectory), revealing a structural unity between discriminative and generative modeling.
\item
\emph{VarPro solvers.} We eliminate the classifier by variable projection and develop, implement, and benchmark first- and second-order solvers against numerous baselines. Prior VarPro methods eliminate the final affine layer of a NN~\cite{newman2021:varpro, newman2022:slimtrain, dus2026:grassmannian}. Here the eliminated block instead acts on the output of an RKHS-parameterized flow; we compute its reduced gradient with the flow's discrete adjoint rather than by backpropagation through a network. For density estimation the eliminated block is a Gaussian-mixture reference fit by EM, outside the VarPro theory but stable in practice.
\item
\emph{Deterministic exact-divergence CNF.} The RKHS parameterization yields a closed-form expression for the divergence at the cost of a single kernel matvec. Carrying it into the likelihood-based change of variables makes the CNF objective deterministic and its discrete-adjoint reduced gradient exact, avoiding both the variance of stochastic trace estimators~\cite{grathwohl2018:ffjord} and the potential-flow restriction of exact-trace methods~\cite{onken2021:ot}. The resulting transport is also a diffeomorphism whose regularity is explicitly controlled, rather than a map that is only invertible by construction or diffeomorphic in principle.
\item
\emph{Expressivity extensions.} We enhance expressivity on both sides of the transport: a multiscale mixture of isotropic kernels captures deformations across spatial scales, and a trainable Gaussian-mixture (GMM) reference matches multimodal targets so the flow need barely deform them, cutting the round-trip reconstruction error at comparable likelihood on the synthetic targets and lowering the likelihood on real data.
\item
\emph{Control-point reduction.} We cut the per-solve cost from $O(n_s^2)$ to $O(n_s m)$ by expanding the velocity over $m \ll n_s$ control points instead of all $n_s$ particles. This Nystr\"om-style reduction preserves all properties stated above; what it gives up is the exactness of the representer parameterization for the supervised loss, since the kernel centers are decoupled from the data particles.
\end{itemize}

\noindent The code and scripts reproducing the results will be released upon acceptance; \Cref{app:hardware} details the computing environment and reproducibility settings.

\subsection{Limitations}\label{sec:limitations}

The kernel coupling makes each forward and adjoint solve of the full model scale as $O(n_s^2 dn_t)$, which the control-point reduction lowers to $O(n_s m)$. How best to place and adapt the control points for a target accuracy, and how far the reduction carries on very large, high-dimensional data, remain open. Finally, kernel methods degrade in high dimensions: we handle real-data problems through a latent-space embedding rather than operating in the ambient feature space; large-scale, high-dimensional generative modeling remains open.

\subsection{Outline}\label{sec:outline}

\Cref{sec:methods} develops the unified framework. Adopting a discretize-then-optimize perspective, we formulate training as a finite-dimensional optimal control problem over RKHS-parameterized particle flows, derive the Lagrangian and first-order optimality conditions, instantiate the two objectives, and describe the numerical approach. \Cref{sec:results} reports experiments on both classification and generative modeling. Additional derivations, implementation details, and complementary results are in the appendix.

\section{Problem Formulation and Numerics}\label{sec:methods}

We propose a mathematical framework for supervised classification and likelihood-based density estimation based on an optimal control formulation. Our formulation has four task-agnostic components: \begin{enumerate*}[label={\it(\roman*)}] \item a flow of diffeomorphisms generated by a time-dependent vector field $v$; \item an RKHS parameterization that reduces $v$ to a finite set of control coefficients $a$ and turns the control cost into a quadratic kernel form; \item a discretize-then-optimize treatment that yields a finite-dimensional problem; and \item a Lagrangian problem formulation with adjoint-based reduced gradients\end{enumerate*}. The two learning tasks enter through the objective imposed on the transported particles, while sharing the same dynamics, control cost, and numerical treatment.

We develop the approach and instantiate it for classification, following the diffeo\-morphic-learning formulation of~\cite{younes2020:diffeomorphic}, then develop the density-estimation instance.

\subsection{Flows of Diffeomorphisms}

Let $k \defeq (k_1,\ldots,k_d) \in \mathbb{N}^d$ be a multi-index, $\partial^k$ the differential operator of order $|k| = \sum_{i=1}^d k_i$, and $q \in \mathbb{N}$. We write $C^q(\mathcal{X},\mathbb{R})$ for the $q$-times continuously differentiable functions on $\mathcal{X}$, equipped with $\|u\|_{q,\infty} = \max_{0 \leq |k| \leq q}\|\partial^k u\|_{\infty}$, and denote by $C^q_0(\mathcal{X},\mathbb{R}^d)$ the class-$C^q$ vector fields that vanish at infinity together with their derivatives up to order $q$. Throughout the well-posedness discussion we take $\mathcal{X} = \mathbb{R}^d$, so that the flow, the diffeomorphism group, and the geodesic distance are all defined on the same space; the data themselves occupy a bounded subset of $\mathbb{R}^d$.

Flows of $\mathbb{R}^d$-diffeomorphisms $z_t \defeq z(t,\cdot\,)$ are given by solutions of~\cref{eq:flow}, where $v \in \mathcal{V}$, $\mathcal{V} \defeq L^2([0,1],\mathcal{H})$, $t \mapsto v_t \defeq v(t,\cdot\,)$, $v_t \in \mathcal{H}$; $\mathcal{H}$ is a Hilbert space continuously embedded into $C^q_0(\mathcal{X}, \mathbb{R}^d)$. The regularity required for the flow to be well-posed is determined by the Sobolev embedding theorem and depends on $d$.

\begin{assumption}[Regularity]\label[assumption]{asm:regularity}
$\mathcal{H}$ is continuously embedded into $C^q_0(\mathcal{X}, \mathbb{R}^d)$ for some $q \geq 1$, and the velocity field satisfies $v \in L^2([0,1],\mathcal{H})$.
\end{assumption}

\noindent Under \Cref{asm:regularity}, \cref{eq:flow} admits a unique solution whose displacement satisfies $z - \operatorname{id}_{\mathbb{R}^d} \in C([0,1], C^q_0(\mathbb{R}^d, \mathbb{R}^d))$, each map $z_t : \mathbb{R}^d \to \mathbb{R}^d$ is a $C^q$-diffeomorphism with a $C^q$ inverse (the flow of $-v$ run backward in time), and the displacement map $v \mapsto z^v - \operatorname{id}_{\mathbb{R}^d}$ is locally Lipschitz from $L^2([0,1],\mathcal{H})$ into $C([0,1], C^0_0(\mathbb{R}^d,\mathbb{R}^d))$ with the sup norm. This is an established result~\cite{dupuis1998:variational, trouve1998:diffeomorphisms}; see also~\cite[Chap.~7]{younes2019:shapes} and~\cite{glaunes2008:large} for the matrix-kernel setting.

The admissible maps $\phi \defeq z_1$ at $t=1$ form a subgroup $\mathcal{G}$ of diffeomorphisms of $\mathbb{R}^d$, equipped with the right-invariant metric given by the minimal geodesic path length between elements. The squared geodesic distance between $\operatorname{id}_{\mathcal{X}}$ and $\phi \in \mathcal{G}$ is
\begin{equation}\label{eq:riemannian-dist}
\operatorname{dist}^2_{\mathcal{G}}(\operatorname{id}_{\mathcal{X}},\phi) \defeq
\inf_{v\in \mathcal{V}}
\left\{
\int_0^1 \| v_t \|^2_{\mathcal{H}} \,\mathrm{d} t
: \phi = z_1,\,\,
z \;\text{satisfies}\; \cref{eq:flow}
\right\}.
\end{equation}

\noindent By Cauchy--Schwarz the energy infimum equals the squared length infimum (equality at constant speed, which leaves $\phi=z_1$ unchanged), so $\operatorname{dist}_{\mathcal{G}}$ is a right-invariant distance. This specifies the first building block of~\cref{eq:varopt-training}.

\subsection{Terminal Objective: Supervised Classification}\label{sec:classifier}

We specify the classifier $c_{\Xi}$ and its loss in~\cref{eq:varopt-training}. We model class membership by a probability simplex; a natural loss is the cross-entropy given by the negative log-likelihood (NLL) of a multi-class logistic regression model,
\begin{equation}\label{eq:loss}
\operatorname{loss}(c_{\Xi}(z(1)), y)
= -\sum_{i = 1}^{n_s}
\Big\langle
y_i, \log\frac{\exp(\Xi z_i(1))}{\langle e_{n_c},\exp(\Xi z_i(1))\rangle}
\Big\rangle,
\end{equation}

\noindent $e_{n_c} = (1,\ldots,1) \in \mathbb{R}^{n_c}$. This loss is smooth and convex for fixed $z$. The coupling between $\Xi$ and $z$ and~\cref{eq:flow} render the overall problem non-linear and non-convex.

We regularize the parameter matrix $\Xi \in \mathbb{R}^{n_c \times d}$ with a Tikhonov-type model $\operatorname{reg}(\Xi) = \lambda_{\scr{cls}}\|B(\Xi - \bar{\Xi})\|_F^2$ with regularization parameter $\lambda_{\scr{cls}} > 0$, regularization operator $B \succ 0$, and reference model $\bar{\Xi} \in \mathbb{R}^{n_c \times d}$. We select $\bar{\Xi} = 0$ and $B = \operatorname{diag}(1,\ldots,1) \in \mathbb{R}^{n_c \times n_c}$. This defines the remaining two components of~\cref{eq:varopt-training}.

\subsection{RKHS Representation}

The vector field $v$ lives in a vector-valued RKHS $\mathcal{H}$~\cite{aronszajn1950:theory, scholkopf2002:learning}. We parameterize it by a finite kernel expansion over the training particles,
\[
\textstyle v_t(x) = \sum_{i=1}^{n_s} \kappa(x,z_i(t)) a_i(t) \quad \text{for
all}\; x \in \mathcal{X},
\]

\noindent with a matrix-valued kernel $\kappa : \mathcal{X} \times \mathcal{X} \to \mathbb{R}^{d\times d}$, data points $\{z_i\}_{i=1}^{n_s} \subset \mathcal{X}$, and time-dependent coefficients $a(t) = (a_i(t))_{i=1}^{n_s}$, $a_i(t) \in \mathbb{R}^d$. This reduces the search for $v$ to the coefficients $a_i : [0,1]\to\mathbb{R}^d$. The RKHS construction and the representer argument justifying this parameterization are given in \Cref{app:kernel-def,app:representer}.

We simplify the construction by modeling the kernel as $\kappa(\,\cdot\,, \,\cdot\,) = \operatorname{ker}_{\sigma}(\,\cdot\,, \,\cdot\,) I_d$, $I_d \defeq \operatorname{diag}(1,\ldots,1) \in \mathbb{R}^{d \times d}$, $\operatorname{ker}_{\sigma} : \mathcal{X} \times \mathcal{X} \to \mathbb{R}$, parameterized by the bandwidth $\sigma > 0$. The choice of the kernel determines the associated RKHS; we select
\begin{equation}\label{eq:gaussiankernel}
\operatorname{ker}_{\sigma}(u,w)
= \exp\left(
-\tfrac{1}{2} \left\langle u - w, \Sigma^{-1}(u - w)\right\rangle
\right),
\end{equation}

\noindent $\Sigma = \sigma^2 I_d$, for any $u,w\in\mathcal{X}$. This kernel is symmetric, positive definite, and yields a metric invariant under rotation and translation; its RKHS satisfies \Cref{asm:regularity} for every $q\ge1$ (see \Cref{app:kernel-def}). The Hilbert space $\mathcal{H}_{\sigma}^d$ of smooth vector fields in $\mathcal{X}$ (the space $\mathcal{H}$ of \Cref{asm:regularity}) is defined as
\begin{equation}\label{eq:rkhs}
\mathcal{H}_{\sigma}^d
\textstyle
\defeq
\overline{
\left\{
w : w(\,\cdot\,) = \sum_{i=1}^n \ker_{\sigma}(x_i,\,\cdot\,) a_i,
\,\, a_i \in \mathbb{R}^d,
\,\, x_i \in \mathcal{X},
\,\,n \in \mathbb{N}
\right\}}^{\,\|\cdot\|_{\mathcal{H}_\sigma^d}}.
\end{equation}

\noindent The norm on $\mathcal{H}_{\sigma}^d$ in~\cref{eq:rkhs} is defined by the square root of
\begin{equation}\label{eq:rkhs-norm}
\langle v_t(\,\cdot\,), v_t(\,\cdot\,)\rangle_{\mathcal{H}_{\sigma}^d}
= \textstyle\sum_{k=1}^{n_s}\sum_{l=1}^{n_s}
\operatorname{ker}_{\sigma}(z_k(t),z_l(t))\,\langle a_k(t),a_l(t)\rangle.
\end{equation}

Associating with each feature $x_k \in \mathbb{R}^d$ a time-dependent state $z_k(t) \in \mathbb{R}^d$ and control $a_k(t) \in \mathbb{R}^d$, the RKHS representation turns the flow~\cref{eq:flow} of $(x_k)_{k=1}^{n_s}$ into
\begin{equation}\label{eq:rkhs-flow}
\textstyle\partial_t z_k(t) = \sum _{l=1}^{n_s} \operatorname{ker}_{\sigma}(z_k(t), z_l(t)) a_l(t)
\end{equation}

\noindent with initial condition $z_k(0) = x_k$, $k=1,\ldots,n_s$. This representation is exact for the supervised loss, whose dependence on $v$ is only through the particle trajectories and the control cost; for density estimation the same expansion is adopted as a modeling choice rather than a consequence of optimality, with the universality of the RKHS ensuring an expressive model class (see \Cref{app:representer}).

\subsection{Numerical Discretization}

The input features $\{x_k\}_{k=1}^{n_s}$ form a long vector $x^0 = (x_i)_{i=1}^{n_s} \in \mathbb{R}^{dn_s}$ at time $t=0$. We subdivide the time interval $[0,1]$ into $n_t$ steps of size $h_t > 0$ and define the discrete state and control vectors by
$a = (a^{j})_{j = 1}^{n_t}
= (a^{j}_1, \ldots, a^{j}_{n_s})_{j = 1}^{n_t}
\in \mathbb{R}^{dn_sn_t}$,
$z = (z^{j})_{j = 1}^{n_t+1}
= (z^{j}_1, \ldots, z^{j}_{n_s})_{j = 1}^{n_t+1}
\in \mathbb{R}^{dn_s(n_t+1)}$. Inserting the RKHS representation of $v_t$ into~\cref{eq:flow} and using an explicit Euler time integration scheme\footnote{We also tested 2nd- and 4th-order Runge--Kutta schemes but observed no significant gains from higher-order time integration for the classification problem.}, the evolution of the features $x^0$ is modeled as
\begin{equation}\label{eq:euler-step-state}
\textstyle z^{j+1}_k = z^j_k + h_t
\sum_{l=1}^{n_s} \operatorname{ker}_{\sigma}(z^j_k, z^j_l)a^j_l,
\quad j = 1,\ldots, n_t,
\end{equation}

\noindent with $z^1_k = x^0_k$, $k = 1,\ldots,n_s$. The terminal state $z^{n_t+1} \eqdef x^1$ at $t=1$ is fed through the (linear) classifier $c_{\Xi}$ in~\cref{eq:loss} parameterized by $\Xi = \left( \xi_1, \ldots, \xi_{n_c} \right) \in \mathbb{R}^{n_c \times d}$ with $\xi_1 = 0$ fixing class~1 as the reference, to map the transported features to their labels $y = (y_k)_{k=1}^{n_s} \in \mathbb{R}^{n_cn_s}$.

We discretize the energy in~\cref{eq:riemannian-dist} by approximating the time integral by a left Riemann summation. With~\cref{eq:rkhs-norm} we obtain the control cost
\begin{equation}
\textstyle
\operatorname{cost}(z,a)
= h_t \sum_{j=1}^{n_t} \sum_{k=1}^{n_s}\sum_{l=1}^{n_s}
\operatorname{ker}_{\sigma}(z_k^j,z_l^j)
\langle a_k^j, a_l^j \rangle
= \langle a, K(z) a \rangle,
\label{eq:discrete_cost}
\end{equation}

\noindent where $K(z) = h_t \operatorname{diag}(K^1, \ldots, K^{n_t})$, with $K^j = (K_{kl}^j I_d)_{k, l = 1}^{n_s,n_s}$, $K_{kl}^j = \operatorname{ker}_{\sigma}(z_k^j,z_l^j)$, so that $K(z)$ acts on $a \in \mathbb{R}^{dn_sn_t}$. The control cost thus depends on the state as well as on the control; in the Lagrangian of \Cref{app:kkt_diffeo} the two are independent variables, and once the state is eliminated through the forward solve $z = z(a)$ we abbreviate $K(a) \defeq K(z(a))$.

\subsection{Optimization Problem}

We use a discretize-then-optimize approach; it yields exact discrete gradients and avoids adjoint-consistency issues of the optimize-then-discretize alternative~\cite{onken2020:dovsod, gholami2019:anode}. The unknowns are the controls $a$ prescribing the dynamics and the classifier parameters $\Xi$; our architecture is parameterized via $\Theta = (a,\Xi) \in \mathbb{R}^{dn_sn_t} \times \mathcal{W}$, where $\mathcal{W} \defeq \{\Xi \in \mathbb{R}^{n_c\times d} : \xi_1 = 0\}$ is the classifier feasible set fixed by the reference-class normalization introduced above. The inputs are the training data $\Psi \defeq \{(x_i,y_i)\}_{i = 1}^{n_s} \subset \mathcal{X} \times \mathcal{Y}$ consisting of $n_s$ feature-label pairs $(x_i,y_i)$. We define the objective function $f: \mathbb{R}^{dn_s(n_t+1)} \times \mathbb{R}^{dn_sn_t} \times \mathbb{R}^{n_c\times d}\to \mathbb{R}$,
\[
f(z,\Theta) \defeq
-\lambda_{\scr{loss}}\sum_{i = 1}^{n_s}
\Big\langle
y_i, \log \frac{\exp(\Xi z_i^{n_t+1})}{\langle e_{n_c}, \exp(\Xi z_i^{n_t+1})\rangle}
\Big\rangle
+ \langle a, K a \rangle + \lambda_{\scr{cls}}\|\Xi\|_F^2.
\]

\noindent Here, $\lambda_{\scr{loss}} > 0$ is the data-misfit weight balancing the terminal loss against the RKHS kinetic energy $\langle a, K a \rangle$ and the classifier regularization with weight $\lambda_{\scr{cls}} > 0$. We determine the trainable parameters $\Theta$ by solving the constrained problem
\begin{equation}\label{eq:varopt-learning-discrete}
\minopt_{z,\,\Theta}\;\;f(z,\Theta)
\quad \text{subject to}\quad z\;\text{satisfies}\;\cref{eq:euler-step-state}.
\end{equation}

\subsection{Optimality Conditions}\label{sec:optimality-conditions}

We solve~\cref{eq:varopt-learning-discrete} by the method of Lagrange multipliers; we introduce the unknown
$p = (p^{j})_{j = 1}^{n_t+1}
= (p^{j}_1, \ldots, p^{j}_{n_s})_{j = 1}^{n_t+1}
\in \mathbb{R}^{(n_t+1)dn_s}$
and form the Lagrangian $\ell : \mathbb{R}^{(n_t+1)dn_s} \times \mathbb{R}^{(n_t+1)dn_s} \times \mathbb{R}^{n_tdn_s} \times \mathbb{R}^{n_c \times d}\to \mathbb{R}$,
\begin{equation}\label{eq:lagrangian}
\ell(z,p,\Theta) =
f(z,\Theta)
+
\sum_{k=1}^{n_s}
\Big\langle
p_k^1, z_k^1 - x^0_k
\Big\rangle
+ \sum_{j=1}^{n_t}
\sum_{k=1}^{n_s}
\Big\langle
p_k^{j+1}, z^{j+1}_k \!- \psi(z^j,a^j)
\Big\rangle,
\end{equation}

\noindent where $\psi(z^j,a^j)$ is the right hand side of~\cref{eq:euler-step-state}. To derive the optimality conditions for any optimal $\Theta^\star$ and primal-dual optimal pair $(z^\star,p^\star)$, we compute derivatives of $\ell$ in \cref{eq:lagrangian} with respect to the control $a$, the parameters $\Xi$, the state $z$, and the adjoint $p$. For simplicity of presentation we drop the optimality superscript $\star$.

Differentiating $\ell$ with respect to $a$ yields the reduced gradient; collecting the per-particle terms into a vector $g_a \defeq (g^j_k(a^j))_{k,j=1}^{n_s,n_t} \in \mathbb{R}^{dn_sn_t}$, it is given by
\begin{equation}\label{eq:reduced-gradient}
\textstyle g^j_k(a^j)
\defeq
h_t\sum_{l=1}^{n_s}\operatorname{ker}_{\sigma}(z_k^j,z_l^j) (2 a_l^j - p_l^{j+1}).
\end{equation}

A common alternative is the \emph{natural gradient}~\cite{younes2019:shapes}, obtained by differentiating with respect to the Riemannian metric induced by~\cref{eq:rkhs-norm}, i.e., by multiplying the Euclidean gradient by the inverse of the metric tensor (here the matrix $K$ of~\cref{eq:discrete_cost}, which carries the factor $h_t$); this gives
\begin{equation}\label{eq:reduced-gradient-natural}
g_k^{\scr{nat},j}(a^j) \defeq 2 a_k^j - p_k^{j+1}.
\end{equation}

At optimality, $g_a$ needs to vanish. We evaluate~\cref{eq:reduced-gradient} for candidate $a$ by solving the associated state and adjoint equations to obtain $z$ and $p$. Differentiating~\cref{eq:lagrangian} with respect to $p$ recovers the strong form of the state equation~\cref{eq:euler-step-state}; differentiating with respect to $z$ gives the adjoint equations,
\begin{equation}\label{eq:adjoint}
\begin{aligned}
p_k^{j-1} &= p_k^j + h_t \sum_{l=1}^{n_s}
\big(
\big\langle p_k^j - 2 a_k^{j-1}, a_l^{j-1} \big\rangle + \big\langle p_l^j, a_k^{j-1}\big\rangle
\big)
\partial_1 \operatorname{ker}_{\sigma}(z_k^{j-1}, z_l^{j-1}),
\\
p_k^{n_t+1} &= -\lambda_{\scr{loss}}\,\Xi^\mathsf{T}\!
\left(
\frac{\exp(\Xi z_k^{n_t+1})}{\langle e_{n_c}, \exp(\Xi z_k^{n_t+1})\rangle} - y_k
\right).
\end{aligned}
\end{equation}

\noindent Here, $\partial_1$ denotes differentiation in the first argument; for the Gaussian kernel in~\cref{eq:gaussiankernel}, we obtain $\partial_1\operatorname{ker}_{\sigma}(z_k^{j-1},z_l^{j-1}) = -\operatorname{ker}_{\sigma}(z_k^{j-1},z_l^{j-1})\Sigma^{-1}(z_k^{j-1}-z_l^{j-1})$. The adjoint equation~\cref{eq:adjoint} is solved backward in time; its recursion (first line) is shared by both learning tasks (density estimation adds an interior divergence source; see \Cref{app:kkt_cnf}), while the terminal condition (second line) is task-specific.

Our approach for the dynamics corresponds to a reduced-space formulation~\cite{biros2005:parallel, mang2015:inexact}; we iterate on the control $a$ alone and enforce the state and adjoint equations exactly. Differentiating~\cref{eq:lagrangian} with respect to the classifier along its feasible set $\mathcal{W}$, i.e., in its free rows $\bar\Xi = (\xi_2,\ldots,\xi_{n_c})$, yields a second reduced gradient $\partial_{\bar\Xi} f$ that also vanishes at optimality; the derivation of the adjoint, the reduced gradient, and the classifier gradient is given in \Cref{app:kkt_diffeo}. Evaluated at the state $z^\star(a)$ and adjoint $p^\star(a,\Xi)$ that solve~\cref{eq:euler-step-state,eq:adjoint}, the right-hand side of~\cref{eq:reduced-gradient} is the exact gradient of the state-eliminated objective $F(a,\Xi) \defeq f(z^\star(a),a,\Xi)$ with respect to $a$ (and the classifier expression equals $\partial_{\bar\Xi} F$). We reserve $\widehat{f}$ for the fully reduced objective of~\cref{eq:reduced-obj}, in which $\Xi$ is eliminated as well.

\subsection{Control-Point Reduction}\label{sec:control-points}

The expansion~\cref{eq:rkhs-flow} places a control coefficient on every training particle $x_k$, so the field couples all $n_s$ particles and each forward and adjoint solve costs $O(n_s^2 d n_t)$ through the dense $n_s\times n_s$ kernel evaluations; the matrix-free backend removes the memory but not the computational complexity. To reduce these costs we expand $v$ over $m \ll n_s$ \emph{control points} $c = (c_i)_{i=1}^m$ rather than the $n_s$ particles~\cite{durrleman2014:morphometry} to obtain
\begin{equation}\label{eq:cp-velocity}
\textstyle v_t(x) = \sum_{i=1}^m \operatorname{ker}_\sigma(x, c_i(t))\, a_i(t), \qquad m \ll n_s.
\end{equation}

\noindent This corresponds to a subset-of-regressors (Nystr\"om) reduction of the RKHS~\cite{williams2001:nystrom, snelson2006:sparse, quinonero2005:unifying}. The $m$ control points carry the momentum $a$ and flow under their own $m\times m$ dynamics, while the $n_s$ data particles are advected by this shared field through
\begin{equation}\label{eq:cp-flow}
\textstyle c_i^{j+1} = c_i^j + h_t \sum_{l=1}^m \operatorname{ker}_\sigma(c_i^j, c_l^j)a_l^j,
\quad
z_k^{j+1} = z_k^j + h_t \sum_{i=1}^m \operatorname{ker}_\sigma(z_k^j, c_i^j)a_i^j.
\end{equation}

\noindent This expression replaces the full-particle discretization~\cref{eq:euler-step-state}. Because the field depends only on $(c,a)$ and not on the data particles $x$, we gain in several aspects.
\begin{enumerate*}[label={\it(\roman*)}]
\item Each forward and adjoint solve is linear in $n_s$ for fixed $m$ (the control flow is $O(m^2)$) and the reduced control lives in $\mathbb{R}^{dmn_t}$, independent of $n_s$.
\item The transport is still a diffeomorphism, since a finite kernel expansion is an RKHS element; the closed-form divergence used for density estimation carries over verbatim, now summed over the $m$ control points at $O(n_s m d)$ cost.
\item Out-of-sample evaluation advects a new point through the $m$-point field at $O(m)$ cost. The learned model is a compact parametric flow rather than the full $n_s$-particle expansion. The same Lagrangian machinery applies and the reduced GN and VarPro solvers carry over, but the reduced optimality system is \emph{not} identical to the full-model one: it acquires a second adjoint on the $m$ control points and a cross-kernel term through which the $n_s$ data particles enter only in sums (see \Cref{app:cp-kkt}).
\end{enumerate*}

\subsection{Density Estimation via CNFs}
\label{sec:cnf}

We now turn to the density-estimation objective. NFs approximate complex densities by transporting a simple reference distribution through an invertible map. This approach permits exact likelihood evaluation via the change-of-variables formula, unlike generative adversarial networks (no explicit density) or variational autoencoders (VAEs), which give only a likelihood lower bound~\cite{rezende2016:variational, papamakarios2021:normalizing, kobyzev2021:normalizing}. Where discrete flows use architectures with a cheap Jacobian determinant~\cite{dinh2015:nice, dinh2017:density, kingma2018:glow, papamakarios2018:masked}, CNFs replace the discrete composition by an ODE~\cite{chen2018:neural, grathwohl2018:ffjord}, so that the log-density evolves through an instantaneous change-of-variables formula rather than a determinant. This is the flow~\cref{eq:flow}. We reuse the RKHS parameterization, the particle dynamics~\cref{eq:rkhs-flow}, the Euler discretization~\cref{eq:euler-step-state}, the quadratic control cost~\cref{eq:discrete_cost}, and the control-point reduction; only the objective changes, from the supervised loss of~\Cref{sec:classifier} to the likelihood below.

\subsubsection{Likelihood Objective and Exact Divergence}
Let $\pi_0$ denote the (unknown) data distribution on $\mathcal{X}=\mathbb{R}^d$ and $\pi_1$ a tractable reference density. Transporting a data sample $x\sim\pi_0$ along the flow~\cref{eq:flow} and writing $\rho_t \defeq (z_t)_\#\pi_0$ for the transported density (so $\rho_0=\pi_0$ and $\rho_1=\phi_\#\pi_0$), the log-density evolves according to the instantaneous change-of-variables formula~\cite{chen2018:neural}
\begin{equation}
\label{eq:instant_cov_v}
\frac{\mathrm{d}}{\mathrm{d}t} \log \rho_t(z(t)) = -\operatorname{div} v_t(z(t)),
\end{equation}

\noindent which, integrated over $[0,1]$, gives $\log\rho_0$ in terms of $\log\rho_1$ and the accumulated divergence. Substituting the reference $\pi_1$ for the (unknown) terminal density $\rho_1$ turns this into the model log-density, the pullback $\phi^{\#}\pi_1 \defeq (\phi^{-1})_\#\pi_1$ of the reference, evaluated in closed form and equal to the data log-likelihood $\log\pi_0$ once training drives $\rho_1 = \phi_\#\pi_0\to\pi_1$,
\begin{equation}
\label{eq:cnf_formula}
\log (\phi^{\#}\pi_1)(x)
= \log \pi_1(z(1))
+ \int_0^1 \operatorname{div} v_t(z(t))\,\mathrm{d}t.
\end{equation}

As in dynamical optimal transport, we match the pushforward of $\pi_0$ to the reference $\pi_1$ while penalizing the kinetic energy of $v$. The training problem reads
\begin{equation}\label{eq:varopt-cnf}
\minopt_{v \in \mathcal{V}}\;\;
\lambda_{\scr{loss}}\operatorname{dist}_{\scr{kl}}(\phi_\#\pi_0 \,\|\, \pi_1)
+ \int_0^1 \| v_t \|_{\mathcal{H}}^2\,\mathrm{d}t,
\end{equation}

\noindent where $\phi = z(1)$, $\phi  : \mathcal{X}\to\mathcal{Z}$ is the time-one \emph{forward} (normalizing) flow map with $z$ satisfying~\cref{eq:flow}, and $\phi_\#\pi_0$ is its pushforward. Training drives $\phi_\#\pi_0$ toward $\pi_1$, so the learned model density is the pullback $\phi^{\#}\pi_1$, obtained in closed form via~\cref{eq:cnf_formula}. Generation reverses this: drawing $z\sim\pi_1$ and applying the \emph{inverse} map $\phi^{-1}$ (the flow~\cref{eq:flow} integrated backward from $t=1$ to $t=0$) produces samples, so the same flow yields likelihoods through $\phi$ and samples through $\phi^{-1}$.

Approximating the expectation over $\pi_0$ in~\cref{eq:varopt-cnf} by a Monte Carlo average over the samples $\{x_k\}_{k=1}^{n_s}$ and substituting~\cref{eq:cnf_formula} (equivalently, minimizing the expected NLL up to the $v$-independent data entropy; see \Cref{app:kkt_cnf}), we obtain the CNF objective function
\begin{equation}
\label{eq:cnf-continuous-objective}
f(v,z) \defeq
\frac{\lambda_{\scr{loss}}}{n_s}\sum_{k=1}^{n_s}
\left[
-\log \pi_1(z_k(1))
-\!\int_0^1\!\!\operatorname{div} v_t(z_k(t))\,\mathrm{d}t
\right]
+ \!\int_0^1\!\!\| v_t\|_{\mathcal{H}}^2\,\mathrm{d}t,
\end{equation}

\noindent where $\lambda_{\scr{loss}} > 0$ weights the data misfit averaged over the $n_s$ samples against the kinetic energy. Using~\cref{eq:cnf-continuous-objective} the CNF optimal control problem reads
\begin{equation}
\label{eq:rkhs_cnf_continuous}
\minopt_{z \in \mathcal{T},\,v \in \mathcal{V}}\;\; f(v,z)
\quad \text{subject to} \quad z\,\,\text{satisfies \cref{eq:flow}}.
\end{equation}

The velocity, particle dynamics, and control cost carry over from the classification formulation; the new ingredient is the divergence term in $f$, which the RKHS parameterization renders exact and closed-form.

\begin{lemma}[Exact divergence for RKHS velocity fields]\label[lemma]{lem:exact-div}
Let $\operatorname{ker}_\sigma \in C^1(\mathcal{X} \times \mathcal{X}, \mathbb{R})$ and $\kappa(\,\cdot\,,\,\cdot\,) = \operatorname{ker}_\sigma(\,\cdot\,,\,\cdot\,)I_d$. For $v_t(x) = \sum_{l=1}^{n_s} \operatorname{ker}_\sigma(x, z_l(t))a_l(t)$, $a_l(t) \in \mathbb{R}^d$, the divergence at any particle position $z_k(t)$ is
\begin{equation}\label{eq:exact-div}
\textstyle\operatorname{div} v_t(z_k(t))
= \sum_{l=1}^{n_s} \langle(\partial_1 \operatorname{ker}_\sigma)(z_k(t), z_l(t)), a_l(t)\rangle.
\end{equation}

\noindent For the Gaussian kernel~\cref{eq:gaussiankernel}, we obtain
\begin{equation}\label{eq:exact-div-gauss}
\textstyle
\operatorname{div} v_t(z_k(t))
= -\sum_{l=1}^{n_s} \operatorname{ker}_\sigma(z_k(t), z_l(t))
\langle \Sigma^{-1}(z_k(t) - z_l(t)), a_l(t)\rangle.
\end{equation}

\noindent The evaluation of~\cref{eq:exact-div-gauss} for all $k = 1,\ldots,n_s$ costs $\mathcal{O}(dn_s^2 )$, i.e., the same as a single kernel matvec.
\end{lemma}

Both identities follow by termwise differentiation of the kernel expansion, where $\partial_1$ denotes differentiation in the first argument. This is the closed-form divergence of a kernel field, akin to kernelized Stein methods~\cite{liu2016:stein}. The novelty is not the identity but that it renders the CNF objective deterministic and its reduced gradient exact at the $\mathcal{O}(dn_s^2)$ cost that the GN solvers of~\Cref{sec:outer-step} exploit.

Unlike the Hutchinson trace estimator of~\cite{grathwohl2018:ffjord}, which introduces gradient variance, and~\cite{onken2021:ot}, which recovers linear cost only by restricting $v_t$ to potential flows, the expression~\cref{eq:exact-div} is deterministic and exact and places no structural constraint on $v_t$ beyond the RKHS parameterization (see \Cref{app:res-diveval}). Because the objective stays deterministic and its reduced gradient exact, the reduced GN solvers carry over to density estimation with the curvature modifications of \Cref{app:varpro-variants}. The CNF inherits the state equation~\cref{eq:euler-step-state} and control cost~\cref{eq:discrete_cost} of the classification problem; the adjoint and control equations~\cref{eq:adjoint,eq:reduced-gradient} follow after replacing the terminal condition and adding the divergence source term.

\subsubsection{Discretization and Optimality Conditions}

Discretizing $f$ in~\cref{eq:cnf-continuous-objective} using~\cref{eq:euler-step-state} and the closed-form divergence of~\Cref{lem:exact-div}, the discrete objective is
\begin{equation}\label{eq:cnf-discrete-objective}
f(a,z)
=
\frac{\lambda_{\scr{loss}}}{n_s}\sum_{k=1}^{n_s}
\Big(
-\log\pi_1(z_k^{n_t+1})
-h_t\sum_{j=1}^{n_t}\operatorname{div} v_{t^j}(z_k^j)
\Big)
+ \langle a, Ka\rangle,
\end{equation}

\noindent with $\operatorname{div} v_{t^j}(z_k^j)$ given in closed form by~\cref{eq:exact-div} and $K$ the block-diagonal control-cost matrix introduced in~\cref{eq:discrete_cost}. The discrete problem then reads
\begin{equation}
\label{eq:cnf_discrete}
\minopt_{z,\,a}\;\; f(a,z)
\quad \text{subject to}\quad z\;\text{satisfies}\;\cref{eq:euler-step-state}.
\end{equation}

\noindent
Relative to the supervised classification case, its terminal adjoint is now driven by the score $\partial_z\log\pi_1$ rather than the softmax residual, with an added interior divergence source; the optimality conditions are given in \Cref{app:kkt_cnf}.

\subsubsection{GMM Reference Distribution}
\label{sec:gmm-reference}

To model multi-modal target densities, we replace the unimodal Gaussian (UGM) reference of~\cref{eq:rkhs_cnf_continuous} by a GMM as the terminal density $\pi_1$. The flow, the control cost, and the change-of-variables relation are unchanged; only the terminal density and the score $\partial_z\log\pi_1$ entering the terminal adjoint condition are modified. We define $\pi_1$ as a mixture of $n_{\scr{mod}}$ Gaussians,
\[
\textstyle\pi_1(z) = \sum_{i=1}^{n_{\scr{mod}}}
w_i\mathcal{N}(z; \mu_i, \Sigma_i),
\qquad
\mathcal{N}(z; \mu_i, \Sigma_i)
= c_i
\exp\big(
-\tfrac{1}{2} \left\| z - \mu_i\right\|_{\Sigma_i^{-1}}^2
\big),
\]

\noindent where $c_i = (2\pi)^{-d/2}|\Sigma_i|^{-1/2}$, $w_i > 0$, $\sum_{i=1}^{n_{\scr{mod}}} w_i = 1$, $\mu_i \in \mathbb{R}^d$, $\Sigma_i \succ 0$. In the absence of prior information, we use uniform weights $w_i = 1 / n_{\scr{mod}}$. The unknown covariance matrices $\Sigma_i$ may be chosen as isotropic, diagonal, or unstructured.

The terminal contribution $-\log\pi_1(z_k^{n_t+1})$ is evaluated stably by the log-sum-exp identity applied to the component log-densities $\ell_i = \log w_i + \log\mathcal{N}(z_k^{n_t+1};\mu_i,\Sigma_i)$. The score function $\partial_z \log \pi_1(z)$ admits the closed-form expression
\[
\textstyle\partial_z \log \pi_1(z)
=-\sum_{i=1}^{n_{\scr{mod}}}
\gamma_i(z)\Sigma_i^{-1}(z - \mu_i),
\]

\noindent where $\gamma_i(z) = w_i\mathcal{N}(z; \mu_i, \Sigma_i) \big/ \sum_{l=1}^{n_{\scr{mod}}} w_l \mathcal{N}(z; \mu_l, \Sigma_l)$ is the posterior responsibility.

The mixture parameters may be trained. We use the means $\mu$ directly, reparameterize the weights by logits $\eta$ ($w_i=\operatorname{softmax}(\eta)_i$) and the diagonal covariances by log-variances $\zeta$ ($\sigma_{i,j}^2=\exp(\zeta_{i,j})$), so the updates are unconstrained. We set $\Phi=(\eta,\mu,\zeta)$ and $\Theta=(a,\Phi)$ (and $\Theta=a$ when $\pi_1$ is fixed). The parameters $\Phi$ enter only through the terminal NLL $\operatorname{nllike}(\Phi,z)=-(1/n_s)\sum_{k=1}^{n_s}\log\pi_1(z_k^{n_t+1};\Phi)$; the gradients of $\operatorname{nllike}$ with respect to $\Phi$ are available in closed form in terms of the responsibilities $\gamma_i$.

\subsubsection{Multiscale and Sparsity}\label{sec:multiscale}

Motivated by~\cite{sommer2013:sparse}, we extend the single-scale formulation to a multiscale velocity model. Instead of a single kernel bandwidth, we use a set of $n_{\scr{ker}}$ distinct bandwidths $\mathcal{S} = \{\sigma_1,\ldots,\sigma_{n_{\scr{ker}}}\}$, each defining an isotropic kernel $\operatorname{ker}_\sigma : \mathcal{X} \times \mathcal{X} \to \mathbb{R}$, and represent $v$ as a superposition
\[
\textstyle v_{t^j}(x) =
\sum_{\sigma\in\mathcal{S}}
\sum_{l=1}^{n_s}
\operatorname{ker}_{\sigma}(x,z_l^j)a_{l,\sigma}^j,
\]

\noindent where $a_{l,\sigma}^j \in \mathbb{R}^d$ is the control of particle $l$ at time $t^j$ and scale $\sigma \in \mathcal{S}$. The stacked control cost $\sum_{\sigma\in\mathcal{S}} \langle a_{\sigma}, K_{\sigma} a_{\sigma}\rangle$, with $K_{\sigma} = h_t \operatorname{diag}(K_{\sigma}^1, \ldots, K_{\sigma}^{n_t})$, $K_{\sigma}^j = (K_{kl,\sigma}^jI_d)_{k,l=1}^{n_s,n_s}$, $K_{kl,\sigma}^j = \operatorname{ker}_{\sigma}(z_k^j,z_l^j)$ and $a_{\sigma}$ stacking the controls $a_{l,\sigma}^j$ over all particles and time steps, is the squared norm on the direct sum $\bigoplus_{\sigma}\mathcal{H}_{\sigma}$. By Aronszajn's sum-of-kernels theorem it dominates the single-kernel energy, so \Cref{asm:regularity} and well-posedness still cover the multiscale model. A large $\sigma$ captures coarse, global deformations while a small $\sigma$ captures fine, localized structure.

Because each $K_\sigma\succ0$, the stacked kinetic term is strictly convex and selects a unique minimum-norm splitting; but this splitting spreads the deformation across all scales rather than activating only those it needs. To favor scale-parsimonious representations we add a nonsmooth group-log sparsity penalty~\cite{candes2008:enhancing}
\begin{equation}\label{eq:scale-sparsity}
\textstyle \operatorname{reg}(a) = h_t\sum_{\sigma\in\mathcal{S}} \lambda_{\scr{sp}}^{\sigma}
\sum_{j=1}^{n_t}\sum_{l=1}^{n_s}
\log(1 + \varepsilon^{-1}\|a_{l,\sigma}^{j}\|_2),
\end{equation}

\noindent with per-scale weight $\lambda_{\scr{sp}}^{\sigma}>0$ and smoothing parameter $\varepsilon>0$. The time-step $h_t$ makes the penalty a consistent quadrature of a running cost, mirroring the kinetic term $\sum_\sigma\langle a_\sigma,K_\sigma a_\sigma\rangle$, so that it is consistent under refinement of $n_t$ and $\lambda_{\scr{sp}}^\sigma$ need not be retuned with the temporal resolution. The penalty is group-separable over the per-particle, per-scale, per-step blocks $a_{l,\sigma}^j$, which we exploit in the proximal solver described in \Cref{app:varpro-variants} (see \Cref{alg:varpro-proxgn}). This construction works for the classification and CNF formulation; see \Cref{app:multiscale_cnf} for the multiscale CNF objective.

\subsection{Baseline Optimization Methods}\label{sec:baseline-opt}

The schemes described here optimize the controls $a$ through the shared adjoint reduced gradient~\cref{eq:reduced-gradient} and therefore apply to both the supervised problem~\cref{eq:varopt-learning-discrete} and the CNF problem~\cref{eq:cnf_discrete}. For concreteness we present them on the supervised problem, where the iterate is $\Theta = (a,\Xi)$ and each iteration updates $a$ together with the classifier parameters $\Xi$ (the joint solver of \Cref{alg:iterative-scheme}); for the CNF instance $\Theta = a$ when the reference $\pi_1$ is fixed and $\Theta = (a,\Phi)$ when the mixture parameters $\Phi$ are trained. We consider several baseline algorithms, in part used in~\cite{younes2020:diffeomorphic} and summarized with their globalization and constants in \Cref{tab:optimization-methods}. All except L-BFGS admit both the Euclidean~\cref{eq:reduced-gradient} and the natural~\cref{eq:reduced-gradient-natural} gradient, the latter only for the full-expansion classification problem.
 We refer to~\cite{boyd2004:convex, nocedal2006:numerical} for details.

The dominant cost of all methods is the forward and adjoint solve, each $O(n_s^2 d n_t)$ due to the $n_s \times n_s$ kernel matrix evaluations of the full expansion; the construction in \Cref{sec:control-points} lowers this to $O(n_s m d n_t)$, linear in $n_s$ for $m \ll n_s$ control points.

\begin{table}
\tabadjust
\caption{Baseline algorithms for numerical optimization; those implemented in~\cite{younes2020:diffeomorphic} are marked ($\ast$). GD, NGD, and their Anderson-accelerated variants use Armijo backtracking ($c_1=\mathrm{1e}\mathord{-}4$, factor $0.6$); the CG and quasi-Newton schemes use a strong-Wolfe line search ($c_2=0.2$ for the Nesterov and Fletcher--Reeves variants, $0.9$ for Polak--Ribi\`ere; L-BFGS uses the \texttt{SciPy} routine's internal line search).}
\label{tab:optimization-methods}
\begin{tabular}{lll}
\toprule
\bf acronym & \bf method & \bf globalization \\
\midrule
GD     & (Euclidean) Gradient Descent ($\ast$)                       & Armijo \\
NGD    & Natural (preconditioned) GD ($\ast$)                        & Armijo \\
AAGD   & Anderson-Accelerated GD                                     & Armijo \\
AANGD  & Anderson-Accelerated NGD                                    & Armijo \\
NAGD   & Nesterov-Accelerated GD                                     & Wolfe \\
NANGD  & Nesterov-Accelerated NGD                                    & Wolfe \\
L-BFGS & Limited-memory Broyden--Fletcher--Goldfarb--Shanno ($\ast$) & Wolfe \\
FRCG   & Fletcher--Reeves Conjugate Gradient ($\ast$)                & Wolfe \\
NFRCG  & FRCG (w/ natural gradient)  ($\ast$)                        & Wolfe \\
PRCG   & Polak--Ribi\`ere Conjugate Gradient  ($\ast$)               & Wolfe \\
NPRCG  & PRCG (w/ natural gradient)  ($\ast$)                        & Wolfe \\
\bottomrule
\end{tabular}
\end{table}

\Cref{alg:iterative-scheme} summarizes the joint solver, which updates the controls $a$ and the classifier parameters $\Xi$ simultaneously, in contrast to the VarPro scheme (see \Cref{sec:varpro}), which eliminates $\Xi$ by an inner solve. All baseline schemes (and, with minor modifications, the VarPro outer step of \Cref{sec:outer-step}) share one template: from a gradient $g^{(k)}$ (the joint gradient $(g_a^{(k)}, g_\Xi^{(k)})$ here, or the reduced gradient $g_a^{(k)} = \partial_a \widehat{f}(a^{(k)})$ for VarPro) they form a search direction $s^{(k)}$ and update $\Theta^{(k+1)} = \Theta^{(k)} + \gamma^{(k)} s^{(k)}$ with a line-search step length $\gamma^{(k)} > 0$. The variants are detailed in \Cref{app:baseline-directions}.

\subsection{Variable Projection (VarPro)}\label{sec:varpro}

The optimization problem in~\cref{eq:varopt-learning-discrete} couples the controls $a$ and the classifier parameters $\Xi$. The practical advantage of VarPro~\cite{golub1973:differentiation, kaufman1975:variable, ruhe1980:algorithms, oleary2013:variable, golub2003:separable} is not a reduction in per-iteration cost but in the quality of the outer gradient: by evaluating the reduced gradient at the optimal $\Xi^\star(a)$, the descent direction for $a$ reflects the true geometry of the reduced problem, potentially leading to significantly fewer outer iterations.

The key observation underlying VarPro is that, for a fixed control sequence $a$, the forward trajectory $z^1,\ldots,z^{n_t+1}$ is fully determined by~\cref{eq:euler-step-state}; the objective reduces to a regularized logistic regression in $\Xi$; that is,
\begin{equation}\label{eq:inner}
\Xi^\star(a)
= \operatorname*{arg\,min}_{\Xi\,:\,\xi_1=0}
\left\{
\lambda_{\scr{loss}}\operatorname{loss}(c_{\Xi}(z^{n_t+1}), y) + \lambda_{\scr{cls}}\|\Xi\|_F^2
\right\}.
\end{equation}

Because the loss in~\cref{eq:loss} is convex in $\Xi$ for fixed $z$ and the Tikhonov regularizer is strongly convex, problem~\cref{eq:inner} has a unique solution for any $a$. VarPro eliminates $\Xi$ by substituting the solution map $a \mapsto \Xi^\star(a)$ into the objective in~\cref{eq:varopt-learning-discrete} to obtain the reduced objective
\begin{equation}\label{eq:reduced-obj}
\widehat{f}(a) \defeq f(z^\star(a),a,\Xi^\star(a)).
\end{equation}

\noindent Minimizing $\widehat{f}$ over $a$ alone is equivalent to solving the original joint problem.

VarPro has recently been used for NN training, where the output-layer (linear) parameters are eliminated and a reduced problem in the nonlinear weights remains; the GN VarPro of~\cite{newman2021:varpro} handles the cross-entropy loss of classification, with stochastic and convergence-theoretic developments in~\cite{newman2022:slimtrain, dus2026:grassmannian}. Our setting is analogous; the classifier $\Xi$ is the linear block and the controls $a$ are the nonlinear parameters. It differs in two aspects central to this work. First, the nonlinear dependence enters through the RKHS-parameterized flow, so the reduced gradient is obtained from the discrete adjoint~\cref{eq:adjoint} rather than by backpropagation through a generic network. Second, the eliminated block is specific to the supervised problem: the density-estimation instance has no convex inner block of comparable structure, an asymmetry we return to at the end of this section.

\subsubsection{Reduced Gradient via the Envelope Theorem}

The key computational advantage of VarPro is that the gradient of $\widehat{f}$ can be computed without differentiating through the inner solver. By the chain rule,
\[
\partial_a \widehat{f}(a)
=
\partial_a f\big|_{\Xi = \Xi^\star}\! + \partial_{\bar\Xi} f\big|_{\Xi = \Xi^\star}
\mathrm{d}_a \bar\Xi^\star(a).
\]

\noindent The second term vanishes because $\Xi^\star(a)$ satisfies the first-order optimality condition $\partial_{\bar\Xi} f = 0$ of the inner problem~\cref{eq:inner}. Here and below, $\partial_{\bar\Xi}$ denotes differentiation in the free rows $\bar\Xi = (\xi_2,\ldots,\xi_{n_c})$, i.e.\ along $\mathcal{W}$; the unrestricted gradient $\partial_\Xi f$ does \emph{not} vanish at $\Xi^\star(a)$, since its first row is $\lambda_{\scr{loss}}\sum_k(\widehat{y}_{k,1} - y_{k,1})z^{n_t+1}_k$, which is generically nonzero. The distinction is immaterial in the chain rule above, because $\mathrm{d}_a\Xi^\star$ has a vanishing first row, but it matters wherever the inner residual is measured or the inner Hessian is inverted. Thus, at any $a$ where $\Xi^\star(a)$ is well defined,
\begin{equation}\label{eq:envelope}
\partial_a \widehat{f}(a) = \partial_a f\big|_{\Xi = \Xi^\star}.
\end{equation}

\noindent That is, the gradient of $\widehat{f}$ equals the standard adjoint gradient evaluated at $\Xi^\star(a)$; no differentiation through the inner solver for $\Xi$ is required. The gradients with respect to $a$ and the adjoint equations needed to evaluate~\cref{eq:envelope} are given in~\cref{eq:reduced-gradient,eq:reduced-gradient-natural,eq:adjoint}.

\begin{proposition}[VarPro equivalence and envelope identity]\label[proposition]{prop:varpro}
Assume $\lambda_{\scr{cls}} > 0$ and that the terminal state $z^{n_t+1} = z^\star(a)$ is $C^1$ in $a$. Then:
\begin{enumerate}[label={\it (\roman*)}, wide=0pt, labelindent=0pt]
\item\label{itm:varpro-i}
For every $a$, the inner problem~\cref{eq:inner} admits a unique minimizer $\Xi^\star(a)$, and the map $a \mapsto \Xi^\star(a)$ is continuously differentiable.
\item\label{itm:varpro-ii}
$(a^\dag, \Xi^\dag)$ is a (local) minimizer of the joint problem~\cref{eq:varopt-learning-discrete} if and only if $\Xi^\dag = \Xi^\star(a^\dag)$ and $a^\dag$ is a (local) minimizer of the reduced objective $\widehat{f}$ in~\cref{eq:reduced-obj}.
\item\label{itm:varpro-iii}
The reduced gradient satisfies the envelope identity~\cref{eq:envelope}:
\[
\partial_a \widehat{f}(a) = \partial_a f(z^\star(a), a, \Xi)\big|_{\Xi = \Xi^\star(a)}.
\]
\end{enumerate}
\end{proposition}

\begin{proof}
The three parts follow from strong convexity of the inner problem~\cref{eq:inner}, the standard VarPro stationary-point equivalence~\cite{golub1973:differentiation,golub2003:separable}, and the envelope theorem applied to $\widehat{f}(a) = f(z^\star(a), a, \Xi^\star(a))$; the full argument is in \Cref{app:varpro-proof}.
\end{proof}

\subsubsection{General Numerical Procedure}

\Cref{alg:varpro-step} overviews a single VarPro iteration in four phases:
\begin{enumerate*}[label={\it(\roman*)}]
\item In the \emph{forward pass} we evaluate~\cref{eq:euler-step-state} to get the trajectory $z$.
\item During the \emph{inner solve} we solve the convex subproblem~\cref{eq:inner} for $\Xi^\star(a)$ with terminal state $z^{n_t+1}$ fixed.
\item Next, we evaluate the gradient at $\Xi^\star(a)$; this requires a \emph{backward pass} (adjoint solve) to find $p$ and provides the search direction for $a$.
\item In the \emph{outer step} we perform a line search and update the control $a$.
\end{enumerate*} All variants below share this four-phase structure; only the search-direction construction (3rd phase) and globalization (4th phase) differ.

\subsubsection{Inner Solver}\label{sec:inner-solver}

The inner problem~\cref{eq:inner} is a $2\lambda_{\scr{cls}}$-strongly convex multinomial logistic regression; the $\lambda_{\scr{cls}}$-ridge already makes it strongly convex (hence identifiable); we fix $\xi_1 = 0$ to select class~1 as the reference. We denote the free rows $\xi_2,\ldots,\xi_{n_c}$ by $\bar\Xi \in \mathbb{R}^{(n_c-1)\times d}$ and the inner objective by
\[
f_{\Xi} \defeq \lambda_{\scr{loss}}\operatorname{loss}(c_{\Xi}(z^{n_t+1}), y)
+ \lambda_{\scr{cls}}\|\Xi\|_F^2.
\]

The gradient is $\partial_{\xi_l} f_{\Xi} = 2\lambda_{\scr{cls}}\xi_l + \lambda_{\scr{loss}}\sum_{i=1}^{n_s} (c_{\Xi}(z_i^{n_t+1}) - y_i)_l z_i^{n_t+1}$, $l = 2,\ldots,n_c$, with the per-class Hessian blocks
\[
\textstyle\partial^2_{\xi_l\xi_{l'}} f_{\Xi} = 2\lambda_{\scr{cls}}\delta_{ll'} I_d + \lambda_{\scr{loss}}\sum_{i=1}^{n_s} D_i[l,l'](z_i^{n_t+1}\otimes z_i^{n_t+1}),
\]

\noindent where $D_i \defeq \operatorname{diag}(\hat{y}_i) - \hat{y}_i \otimes \hat{y}_i \in \mathbb{R}^{n_c\times n_c}$ is the per-sample softmax Hessian at the predicted class probabilities $\hat{y}_i \defeq c_{\Xi}(z_i^{n_t+1})$ and $\delta_{ll'}$ the Kronecker delta. The first row of $\Xi$ is set to zero at every iteration, and this full-expansion VarPro inner solve terminates when $\|\partial_{\bar\Xi} f_{\Xi}\|_F \leq \operatorname{tol}_{\scr{in}}=\mathrm{1e}\mathord{-}4$ or at a Newton-step cap; the control-point runs use a tighter Newton solve. Because $z^{n_t+1}$ is held fixed, no forward or adjoint solve is performed and each gradient or Hessian matvec costs $O(n_s d n_c)$.

Since the inner problem is strongly convex, we solve it (essentially exactly) by Newton's method: for a small free block $d_{\scr{free}} \defeq d(n_c-1)$ we assemble the Hessian and take one Cholesky factorization per iteration, at $O(n_s d_{\scr{free}}^2 + d_{\scr{free}}^3)$ cost, and switch to matrix-free Newton-CG for larger $d_{\scr{free}}$ (threshold $d_{\scr{free}} > 200$). On directions $\tilde{\xi} = (\tilde{\xi}_2,\ldots,\tilde{\xi}_{n_c})$, $\tilde{\xi}_l \in \mathbb{R}^d$, the Hessian matvec reads
\[
\textstyle(H\tilde{\xi})_l
= 2\lambda_{\scr{cls}}\tilde{\xi}_l
+ \lambda_{\scr{loss}}\sum_{i=1}^{n_s} (D_iu_i)_lz_i^{n_t+1},
\quad
(u_i)_{l'} = \langle z_i^{n_t+1}, \tilde{\xi}_{l'}\rangle,
\quad l = 2,\ldots,n_c,
\]

\noindent where $u_i \in \mathbb{R}^{n_c}$ collects the per-class projections of $z_i^{n_t+1}$ onto the blocks $\tilde{\xi}_{l'}$ of the direction $\tilde{\xi}$, with $(u_i)_1 = 0$ since $\tilde{\xi}_1 = 0$; the costs are $O(dn_sn_c)$ per multiply. Strong convexity yields $H \succ 0$.

\subsubsection{Outer Step and Globalization}\label{sec:outer-step}

Given the reduced gradient $g_a^{(k)} = \partial_a \widehat{f}(a^{(k)})$, the outer step computes a search direction $s_a^{(k)}$ and an updated iterate $a^{(k+1)} = a^{(k)} + \gamma^{(k)} s_a^{(k)}$. The step length $\gamma^{(k)} > 0$ is determined by a line search: the GN and reduced-gradient directions use Armijo backtracking, and the CG and L-BFGS variants use a strong-Wolfe line search. During the line search the classifier is held at $\Xi^\star = \Xi^\star(a^{(k)})$; each trial evaluation computes $f(a^{(k)} + \gamma^{(k)} s_a^{(k)}, \Xi^\star)$. Since $f(a^{(k)} + \gamma^{(k)} s_a^{(k)}, \Xi^\star)\ge\widehat{f}(a^{(k)} + \gamma^{(k)} s_a^{(k)})$ for all $\gamma^{(k)}$, with equality at the base point $\gamma^{(k)}=0$ because $\Xi^\star=\Xi^\star(a^{(k)})$ is the inner minimizer there, Armijo decrease of this frozen-$\Xi$ objective already implies decrease of the reduced objective $\widehat{f}$~\cite{golub2003:separable}. The trust-region (TR) GN variant (see \Cref{alg:varpro-trgn}) replaces the line search by an explicit step-quality test on a quadratic model.

We consider several VarPro outer-direction variants, all acting on the reduced gradient $g_a$ (or its natural counterpart $g_a^{\scr{nat}}$). Reduced gradient descent (GD) and the Euclidean two-loop L-BFGS (window $m_{\scr{out}}=10$) reuse the generic directions of \Cref{sec:baseline-opt}; NPRCG is the reduced-space counterpart of the scheme used in~\cite{younes2020:diffeomorphic}; and an Anderson mixing extrapolates the GD fixed-point map~\cite{anderson1965:iterative,walker2011:anderson}.

The two schemes we advocate are reduced GN solvers built on a \emph{matrix-free} GN Hessian; we never form this GN approximation of the reduced Hessian,
\[
  H_{\scr{gn}}(a) = 2K(a) + \lambda_{\scr{loss}} J_a^\mathsf{T}H_{\scr{loss}}J_a,
\]

\noindent with $K(a)$ the block-diagonal control-cost Gram matrix, $J_a=\partial_a z^{n_t+1}$ the terminal sensitivity, and $H_{\scr{loss}}$ the block-diagonal terminal-loss Hessian in $z$ (given in \Cref{app:varpro-variants}); it is positive definite whenever $K(a)\succ0$, i.e.\ for pairwise-distinct particles at every time level (cf.\ \Cref{app:representer}). We evaluate its action through one incremental forward and one incremental adjoint solve per matvec, with all metric-dependent quantities taken in the natural inner product (for the full expansion; the control-point solver runs CG in the Euclidean inner product; see \Cref{app:varpro-variants}). GN-CG solves $H_{\scr{gn}}s_a = -g_a$ inexactly by CG (fixed inexact-Newton forcing 0.5 with at most $m_{\scr{cg}}=10$ iterations), globalized by Armijo backtracking (see \Cref{alg:varpro-gncg}); TR-GN replaces the line search by a Steihaug--Toint TR model in the natural norm (see \Cref{alg:varpro-trgn}). See \Cref{app:varpro-variants} for details.

\begin{remark}[Applicability to the density-estimation problem]\label[remark]{rem:varpro-density}
VarPro hinges on a second parameter block whose inner problem is strongly convex, so that $a\mapsto\Xi^\star(a)$ is unique and smooth and the envelope identity of \Cref{prop:varpro} applies. This holds for the supervised inner problem~\cref{eq:inner}. The density-estimation problem~\cref{eq:cnf_discrete} lacks such a block: with a fixed reference $\pi_1$ one simply optimizes $a$. Training the mixture parameters $\Phi$ of \Cref{sec:gmm-reference} gives a non-convex inner log-likelihood whose solution map is neither unique nor smooth. We nonetheless eliminate $\Phi$ by variable projection, fitting the mixture with expectation--maximization (EM)~\cite{bishop2006:pattern}; the resulting scheme lies outside the scope of \Cref{prop:varpro} but is stable in practice.
\end{remark}

\subsection{Bandwidth Selection}\label{sec:bandwidth}

The kernel bandwidth $\sigma > 0$ defines the spatial scale of $v$: small $\sigma$ produces localized deformations, large $\sigma$ smoother, more global ones. We either fix $\sigma$ heuristically or set it adaptively by the median heuristic of~\cite{liu2016:stein}, $\sigma^2 = \operatorname{median}\bigl(\{\|x_i - x_j\|_2^2 : i < j\}\bigr) / \bigl(2\log(n_s + 1)\bigr)$, which estimates the typical inter-point spacing from the training samples $\{x_i\}_{i=1}^{n_s}$.

\subsection{Implementation Details}

The reduced gradient and the GN Hessian matvecs are computed by the discrete adjoint of \Cref{sec:optimality-conditions}, so they are exact to machine precision at frozen $\Xi$, at working memory independent of the computational-graph length. The dominant remaining cost is the kernel evaluation. For the full expansion, rather than forming the dense $n_s\times n_s$ kernel matrix we use the matrix-free \texttt{PyKeOps} backend~\cite{charlier2021:keops}, which carries out the forward integration, adjoint solve, and gradient evaluation as symbolic, JIT-compiled kernel reductions on CPU or GPU at $O(n_s)$ memory. The control-point solver of \Cref{sec:control-points} instead evaluates the smaller $n_s\times m$ cross-kernel in dense \texttt{PyTorch} and streams the $n_s$ particles in mini-batches through the analytic adjoint, so its working memory stays independent of $n_s$.

\section{Results}\label{sec:results}

We present an empirical study of the proposed framework; all experiments use the computing environment of \Cref{app:hardware}. The appendix collects four supporting studies not reported here: the reduced gradient attains the convergence order of the time integrator (see \Cref{app:adjoint-order}); the discrete flow is inverse-consistent, with a round-trip error that shrinks with the time step (see \Cref{app:inv-consistency}); the closed-form RKHS divergence matches an automatic-differentiation trace to machine precision at a fraction of the cost of stochastic estimators (see \Cref{app:res-diveval}); and higher-order (RK2/RK4) integration gives no accuracy gain over Euler at increased cost (see \Cref{app:integration-order}).

\subsection{Classification}

In the main text we compare the VarPro outer-solver variants, benchmark them against joint-update baseline optimizers, study the computational scalability of the reduced solve, and demonstrate the method on two real-data problems (Pendigits and MNIST). Per-problem flow visualizations and convergence data underlying the main-text figures are reported in \Cref{app:additional-class-results}.

\subsubsection{Considered Problems}

\begin{figure}
\centering
\includegraphics[width=0.75\textwidth]{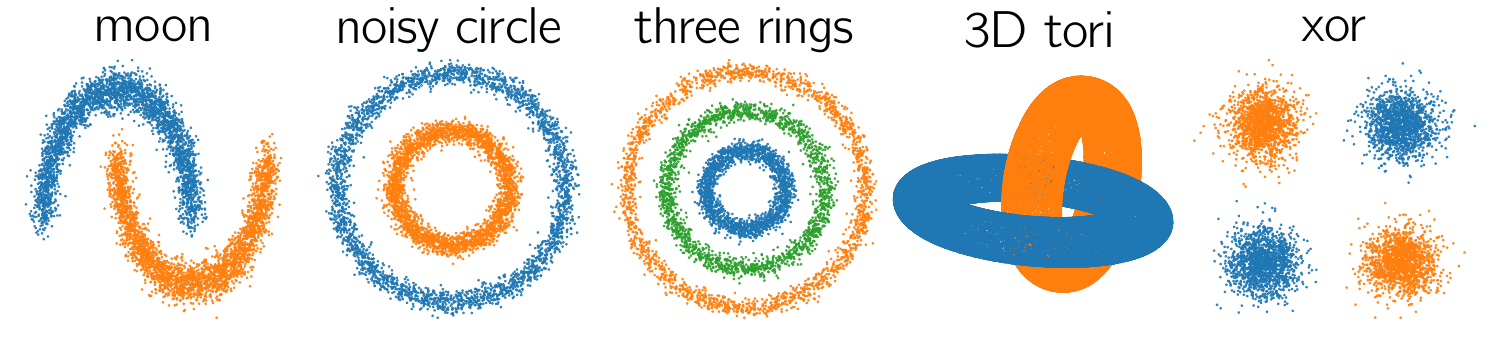}
\caption{The five synthetic classification problems, with points colored by class. From left to right: two moons, noisy circles, three rings, 3D tori, and XOR. Two moons, noisy circles, and XOR are planar binary problems; three rings is a three-class problem; 3D tori is a spatial binary problem of two interlocking tori.}
\label{fig:toy-problems}
\end{figure}

We assess the classifier on five synthetic benchmarks of increasing geometric difficulty (see \Cref{fig:toy-problems}). Each is constructed so that the classes are not linearly separable and a curved decision boundary must be learned. All inputs are standardized to zero mean and unit variance and lifted by one auxiliary affine coordinate before entering the solver.

\subsubsection{Experimental Setup and Metrics}
\label{sec:class-setup}

Unless noted, the classification runs use an Euler flow with the Newton inner $\Xi$-solve of \Cref{sec:inner-solver} (tolerance $\operatorname{tol}_{\scr{in}}=\mathrm{1e}\mathord{-}4$), start from the cold start $a\equiv0$, and report the per-iteration median over five data-instance seeds with a $25/75$ inter-quartile band. The per-problem sample counts, dimensions, and bandwidths, together with the shared solver settings ($n_t=20$, $\lambda_{\scr{loss}}=1$, $\lambda_{\scr{cls}}=\mathrm{1e}\mathord{-}2$), are listed in \Cref{tab:class-settings}. We track the reduced objective $\widehat{f}$, the reduced-gradient norm $\|\partial_a\widehat{f}\|$ (and its relative reduction), the outer-iteration count, and the test accuracy, and measure cost by the hardware-independent flow work
\begin{equation}\label{eq:varpro-work}
n_t\bigl(2 + n_{\scr{fev}} + 2n_{\scr{gev}} + 2n_{\scr{cg}}\bigr)
\end{equation}

\noindent accumulated over outer steps. The leading 2 is the reduced gradient (a forward and an adjoint solve); $n_{\scr{fev}}$, $n_{\scr{gev}}$, and $n_{\scr{cg}}$ are the line-search objective evaluations, gradient evaluations, and GN matvecs; the inner $\Xi$-solve performs no flow solves and is excluded from this count (see \Cref{sec:inner-solver}). It is incurred only by the VarPro schemes, which converge in far fewer outer steps (17--28 versus 58--462 for the joint baselines; see \Cref{tab:varpro-full-baseline}), so counting it would not change the ordering. We summarize each run by the cost to reach a common reference level $\widehat{f}\le\widehat{f}^\star(1+\varepsilon)$, $\varepsilon=5\%$, and stop when the relative reduced-gradient norm falls below $\mathrm{5e}\mathord{-}2$ or the objective stagnates (5 iterations below $\mathrm{1e}\mathord{-}7$ relative change); the variant study of \Cref{sec:varpro-variants} instead runs a fixed 500-iteration budget.

\subsubsection{Comparison of VarPro Outer Solver Variants}
\label{sec:varpro-variants}

We compare optimization methods for the outer control step inside the VarPro loop on the three rings problem, following the protocol of \Cref{sec:class-setup}. For this study each variant runs for a fixed budget of 500 outer iterations with the stopping rules disabled. For the two GN outer solvers the reduced system is solved inexactly by CG with a fixed inexact-Newton forcing and a ten-iteration cap.

\begin{figure}
\centering
\includegraphics[width=\textwidth]{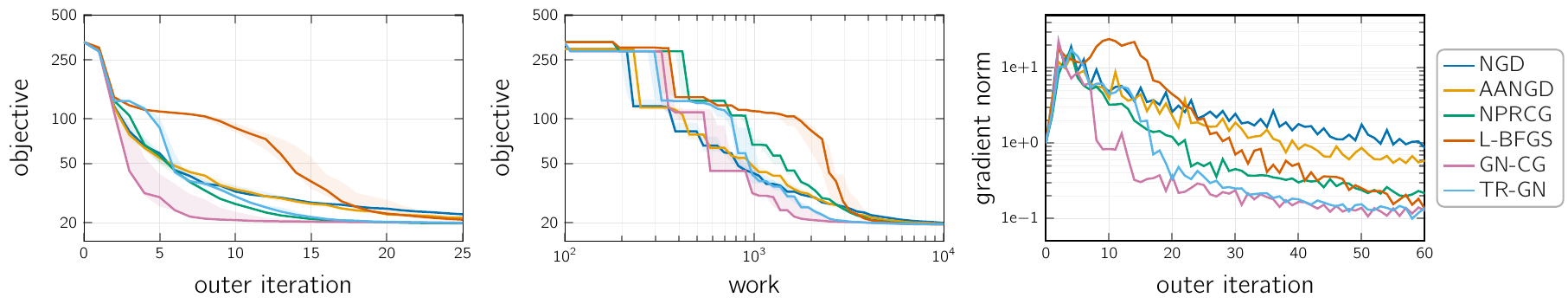}
\caption{VarPro outer-step variants on the three rings problem, median over five data-instance seeds. Reduced objective $\widehat{f}$ against outer-iteration count (left) and computational work (middle), and the relative reduced-gradient norm against the outer iteration (right). The objective panels also shade the $25/75$ inter-quartile band; the gradient panel shows the median only.}
\label{fig:varpro-variants}
\end{figure}

\Cref{fig:varpro-variants} shows the results, with the reported values in \Cref{tab:varpro-variants}. All variants reach the same reduced objective value $\widehat{f}_{\scr{opt}}\approx19.4$. Per outer iteration, the two reduced GN schemes are fastest: GN-CG and TR-GN reach the common reference level in 17 and 21 iterations and reduce the relative reduced-gradient norm by about two orders of magnitude (see the right panel of \Cref{fig:varpro-variants}), whereas the first-order variants need roughly 35--51 iterations and reduce the gradient by only a little over one order of magnitude. On the hardware-independent work axis GN-CG remains most efficient ($\approx$2900), followed by TR-GN ($\approx$3360), the VarPro L-BFGS step ($\approx$4700), and NPRCG ($\approx$4900); this order differs from the iteration ranking, as NPRCG reaches the reference level in only 19 iterations but at a high per-iteration CG cost. The first-order variants (AANGD, NGD) require 6260 and 8360, i.e., 2.2 and 2.9$\times$ the work of GN-CG.

On the reduced problem the RKHS-natural and Euclidean metrics are comparable for the first-order variants. All variants reach $100\%$ test accuracy, up from the near-chance ($\sim$31\%) accuracy of a linear classifier on the raw data.

\subsubsection{Comparison of VarPro Against Joint-Update Baseline Methods}
\label{sec:varpro-vs-baseline}

We explore whether eliminating $\Xi$ and exploiting the reduced GN structure beats optimizing the full vector $(a,\Xi)$ \emph{jointly} with standard methods. We compare GN-CG and TR-GN against a panel of joint optimizers used in~\cite{younes2020:diffeomorphic} on the same three rings instance and protocol as \Cref{sec:varpro-variants}: GD, NAGD, AAGD, FRCG and PRCG forms, and L-BFGS. Each first-order and CG method is run in both the Euclidean and the RKHS-natural (conditioned) inner product; an additional ``N'' marks the natural variant (NGD, NANGD, AANGD, NFRCG, and NPRCG).

\begin{figure}
\centering
\includegraphics[width=\textwidth]{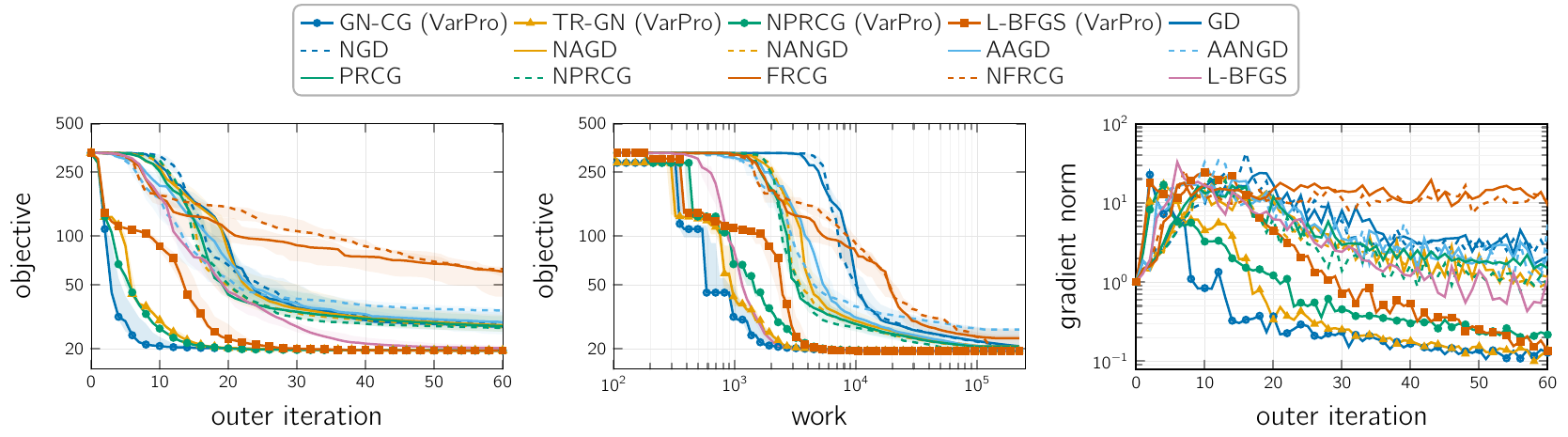}
\caption{Convergence on the three rings problem: the best VarPro variants versus joint $(a,\Xi)$ optimizers, median over five data-instance seeds with $25/75$ inter-quartile band (faint color; objective only). Left: reduced objective $\widehat{f}$ against outer-iteration count; middle: $\widehat{f}$ against cumulative work (log axis); right: relative Euclidean gradient reduction against the outer iteration. VarPro schemes are drawn with markers, joint optimizers without; natural-gradient variants are dashed.}
\label{fig:varpro-full-baseline}
\end{figure}

The comparison (see \Cref{fig:varpro-full-baseline}, with the full baseline numbers in \Cref{tab:varpro-full-baseline}) splits curvature-exploiting methods from the rest: the reduced GN schemes (GN-CG, TR-GN) and the joint quasi-Newton L-BFGS all reach the common reference level at work $\approx$3000, driving the relative gradient down by about two orders of magnitude. Every unconditioned CG method fails to reach this level within the 500-iteration budget (FRCG is unstable and its gradient grows), as does every first-order method except the two natural accelerated variants, AANGD and NANGD, which reach it after 316 and 462 iterations (in 2/5 and 3/5 seeds, against 17--21 iterations in 5/5 seeds for the GN schemes).

VarPro also rescues the classical conditioned CG: optimized jointly over $(a,\Xi)$, NPRCG converges but the coupling to the Euclidean classifier block makes it more than an order of magnitude less efficient than the curvature methods, whereas the \emph{same} scheme on the reduced problem reaches the target in 19 iterations, an $11\times$ work and $17\times$ iteration reduction from decoupling $\Xi$ alone. Conditioning the control block alone or adding acceleration does not substitute for eliminating it.

All methods that fully converge reach the same objective and $100\%$ test accuracy; the reduced GN schemes stand out by reaching it in the fewest iterations (17--21, versus 58 for joint L-BFGS). On the work axis joint L-BFGS is marginally cheaper (2760 versus 2900 for GN-CG) and is the only joint baseline in this range, but it reaches the reference level in only 4/5 seeds against 5/5 for the GN schemes.

We then apply the best methods to the remaining toy examples with the stopping rules of \Cref{sec:class-setup} enabled (see \Cref{fig:varpro-baseline}, \Cref{tab:varpro-baseline}). The VarPro GN schemes outperform joint L-BFGS. On two moons joint L-BFGS terminates early ($\widehat{f}\approx4.8$); relaxing the stopping rule lets it reach the VarPro objective, but only after roughly two orders of magnitude more iterations. All methods classify essentially perfectly, this early-terminated run being the only entry below 100\%.

\subsubsection{Scalability and the Control-Point Reduction}\label{sec:res-scalability}

The dominant cost of our framework is the kernel evaluation. Because $v$ is a kernel expansion over the $n_s$ particles, every flow solve couples all particles through the $n_s\times n_s$ kernel matrix at $O(n_s^2 d n_t)$ cost. The matrix-free \texttt{PyKeOps} backend leaves this unchanged while lowering the memory footprint to $O(n_s)$. We measure both dependencies of the reduced VarPro solve on a single V100 GPU (see \Cref{fig:scalability}) in two experiments:
\begin{enumerate*}[label={\it(\roman*)}]
\item We fix the 3D tori instance and increase $n_s$ from 200 to 3200.
\item We fix the noisy circles instance at $n_s=300$ while isometrically embedding it into $\mathbb{R}^d$ for $d=2,\ldots,256$.
\end{enumerate*}
Each run is executed to convergence over three data-instance seeds. Least-squares fits give wall-time slopes of 1.92 (GN-CG) and 2.00 (TR-GN) against $n_s$, matching the expected $O(n_s^2)$ kernel cost. The dependence on $d$ is sub-linear over the tested range (slope $\approx$0.5, rising toward one at large $d$). The GN outer-iteration count shows no upward trend with $n_s$ or $d$, staying within 7--21 iterations at 100\% accuracy throughout; the growth is entirely in the per-iteration kernel cost.

\begin{figure}
\centering
\includegraphics[width=0.8\textwidth]{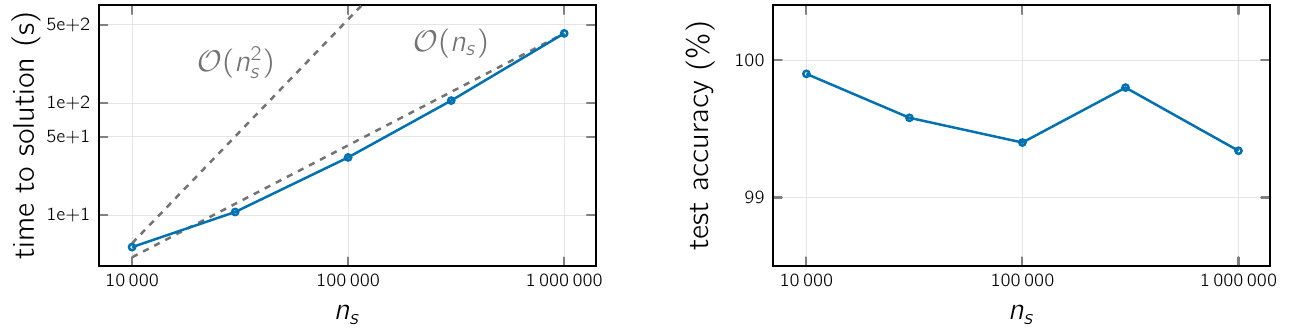}
\caption{Scaling of the reduced GN-CG control-point solver on the three rings problem (resampled at each $n_s$) at a fixed control-point budget $m=128$, timed on a single V100 GPU. Left: time to solution versus the number of samples $n_s$ (log--log), with $\mathcal{O}(n_s)$ and $\mathcal{O}(n_s^2)$ references (dashed). Right: test accuracy versus $n_s$.}
\label{fig:cp-scaling}
\end{figure}

The control-point reduction removes this quadratic term, lowering the per-solve cost to $O(n_s m)$; here, we verify that it is near-lossless. We rerun the reduced GN-CG solver on the five synthetic benchmarks, varying the control-point count $m$ at fixed data size, where $m=n_s$ recovers the full model~\cref{eq:rkhs-flow}; the full accuracy sweep is reported in \Cref{tab:cp-ablation}. A small fraction of control points recovers the full-model accuracy on four of the five problems and comes within 0.8 percentage points on three rings: the simple deformations (two moons, noisy circles, XOR) match exactly at $m=8$ ($\approx$3\% of $n_s$), 3D tori matches at $m=32$, and three rings plateaus at 99.2\% (versus 100\% for the full model) from $m=16$ onward (5--16\% of $n_s$). The reduction thus matches or nearly matches the full model at well under a fifth of the control degrees of freedom; being linear in $n_s$, it extends the solver to data sizes for which the $O(n_s^2)$ full expansion is intractable. \Cref{fig:cp-scaling} makes this concrete: holding $m$ fixed, the streamed reduced GN-CG solver scales linearly (measured slope $\approx$1.0 in $n_s$, versus 1.9 for the full model of \Cref{fig:scalability}) out to $n_s=\mathrm{1e}\mathord{+}6$. This is about 2.5 orders of magnitude beyond the largest full-model problem we could solve in comparable time ($n_s=3200$); the binding constraint is the $O(n_s^2)$ computational cost rather than memory, while the iteration count stays constant at 25 and test accuracy holds near 99\%.

The reduction is equally near-lossless on real, higher-dimensional data. On the VAE-latent MNIST classification of \Cref{sec:res-mnist} (\inum{20000} samples, identical objective and encoder as the full model), the control-point solver recovers the full-model test accuracy to within 0.56 percentage points using $m\le 2048$ ($\approx$10\% of $n_s$) at every latent dimension $d_z\in\{16,32,64\}$, matching it at $d_z=16$ and exceeding it at $d_z=64$ (see \Cref{tab:cp-mnist}). Holding $m=512$ fixed, its test accuracy climbs from 95.2\% to 97.0\% as $n_s$ grows from \inum{5000} to \inum{40000} at a constant iteration count (individual points not reported), reproducing the linear scaling of \Cref{fig:cp-scaling} on real data.

\subsubsection{Pendigits Classification}
\label{sec:res-pendigits}

We apply the reduced GN-CG solver to a real low-dimensional classification task. The UCI Pen-Based Handwritten Digits (Pendigits) dataset comprises \inum{10992} samples in $n_c=10$ classes, each a pen trajectory summarized by 16 features (eight resampled $(x,y)$ pen positions). The inputs are standardized and lifted by one affine coordinate ($d=17$); the velocity field uses a Gaussian kernel with the median-heuristic bandwidth ($\sigma\approx1.35$), an $n_t=20$ Euler flow, loss weight $\lambda_{\scr{loss}}=1$, and classifier regularization $\lambda_{\scr{cls}}=\mathrm{1e}\mathord{-}2$, and the outer step is the GN-CG VarPro solver of \Cref{sec:class-setup} on the matrix-free \texttt{PyKeOps} backend.

On the full stratified split (\inum{7474} training / \inum{3518} test samples), GN-CG converges in 9 outer iterations (relative reduced-gradient norm below tolerance), reducing the objective from $\mathrm{1.7e}\mathord{+}4$ to 69 and reaching 100\% training and 99.0\% test accuracy in 41 seconds on a single V100 GPU. A scaling study (sweeping $n_s$ from \inum{500} to \inum{7000} at a fixed test set, three seeds; see \Cref{app:pendigits-scaling}) confirms the synthetic trends on real data: bounded outer iterations (8--10), 100\% training accuracy throughout, and median test accuracy rising from 93.9\% to 98.9\% at a wall time of under 35 seconds.

\subsubsection{MNIST Classification}
\label{sec:res-mnist}

To test the reduced GN-CG solver on a real, higher-dimensional, multi-class problem, we apply it to MNIST digit classification: $n_s = \inum{20000}$ training and $\inum{5000}$ test images over all $n_c = 10$ classes, each downsampled to $14\times14$ and lifted by one affine coordinate ($d = 197$). The velocity field uses a Gaussian kernel at the median-heuristic bandwidth ($\sigma \approx 1.01$) with $n_t = 20$ Euler steps, loss weight $\lambda_{\scr{loss}} = 1$, and $\lambda_{\scr{cls}} = \mathrm{1e}\mathord{-}2$, and the outer step is the GN-CG VarPro solver, executed on a single V100 GPU with the matrix-free \texttt{PyKeOps} backend.

The solver converges in 12 outer iterations, reducing the objective from $\mathrm{4.6e}\mathord{+}4$ to 186 and reaching 100\% training and 95.2\% test accuracy in $\approx$8 wall-clock hours (dominated by the kernel evaluations). The 100\%/95.2\% train/test gap reflects the generalization limit of the particle-kernel flow in this 197-dimensional raw-pixel setting; classifying in a learned low-dimensional embedding is the natural remedy.

Embedding each image into the latent space of a convolutional VAE (see \Cref{app:vae}), trained here to 16-, 32-, and 64-dimensional latents, and classifying the standardized codes with the identical GN-CG solver (\inum{20000} train / \inum{5000} test), a 32-dimensional latent ($d=33$) reaches 96.8\% test accuracy in 18 iterations and $\approx$8 minutes. This matches or exceeds the raw-pixel 95.2\% at roughly 60$\times$ less cost. The 16-dimensional latent reaches 95.9\% and the 64-dimensional one drops to 95.0\% (below even $d_z=16$), as the kernel flow re-enters the high-dimensional regime the embedding was meant to avoid. Because the control-point reduction makes each solve linear in $n_s$, we train the $d_z=32$ classifier on the full MNIST set (\inum{54000} latents, \inum{6000} held out for validation), a scale the $O(n_s^2)$ full model cannot reach, and evaluate on the standard \inum{10000}-image test set: test accuracy rises to 97.3\% (data not shown). The remaining gap to end-to-end convolutional classifiers reflects the linear classifier on a frozen, unsupervised latent; it does not indicate a limitation of the flow.

\subsection{Generative Modeling}

We report density-estimation results across four synthetic targets, a representative transport, and a real-data MNIST latent-space example; per-problem results are in \Cref{app:genmod-details}.

\subsubsection{Considered Problems}

We evaluate the CNF density-estimation on four 2D target distributions of increasing structural complexity (two concentric circles, a mixture of eight Gaussians (M8G), a pinwheel, and a checkerboard) and, as a real-data example, on MNIST digits mapped into a VAE latent space. The toy problems combine curved support (circles, pinwheel) with multimodal structure (M8G, checkerboard), which stresses both the transport map and the choice of reference distribution.

\subsubsection{Experimental Setup and Metrics}\label{sec:res-genmod-setup}

For all experiments we use $n_s = \inum{10000}$ training samples, a reduction with $m=256$ control points, a forward-Euler flow with $n_t=100$ steps, and $\lambda_{\scr{loss}}=200$. The GMM reference is eliminated by variable projection, its non-convex inner mixture re-fit to the terminal states $z^{n_t+1}=\phi \circ x$ by EM ($50$ sweeps) at each outer step (see \Cref{alg:cnf-varpro-gmm}). The single-scale objective is smooth and is minimized by the reduced GN-CG scheme. The multiscale objective adds the nonsmooth, group-separable scale-sparsity penalty and is minimized by the proximal GN scheme (see \Cref{alg:varpro-proxgn}). The algorithms are executed on a single V100 GPU. In either case we terminate the outer iteration when the dimension-free (max-norm) reduced gradient (respectively, the proximal residual) falls below $\mathrm{1e}\mathord{-}4$, and otherwise cap the count at 200 outer steps;
the multiscale runs also stop on ten consecutive outer steps of stagnant objective (see \Cref{app:varpro-variants}). In practice these flows run to the 200-step cap, except the multiscale M8G/GMM fit, which stopped at 171 steps on that stagnation rule.

We compare single-scale ($\sigma=0.35$) and multiscale ($\sigma\in\{0.35,0.45,0.90\}$) kernels, with bandwidths chosen heuristically; the median heuristic of~\cite{liu2016:stein} gave worse likelihoods. The multiscale scale-sparsity weights are $\mathrm{1e}\mathord{-}4$, $\mathrm{5e}\mathord{-}5$, and $\mathrm{1e}\mathord{-}5$, with smoothing $\varepsilon = \mathrm{1e}\mathord{-}3$.

We report the NLL on the training set and on an independent test set of $n_{\scr{test}} = \inum{10000}$ samples. We also report via the round trip $x_i \mapsto z_i = \phi(x_i) \mapsto \hat{x}_i = \phi^{-1}(z_i)$ on the test set, the mean pointwise error $e_2 = (1/n_{\scr{test}})\sum_{i=1}^{n_{\scr{test}}} \|\hat{x}_i-x_i\|_2$ and the normalized Frobenius error $e_F = \|\hat{X}-X\|_F / n_{\scr{test}}$, where $X,\hat{X}\in\mathbb{R}^{n_{\scr{test}}\times d}$ stack the original and reconstructed test samples ($d=2$ for the toy datasets).

\subsubsection{Synthetic Data Distributions}
\label{sec:syn-data-dist}

For each target we compare a standard UGM against a trainable GMM reference eliminated by variable projection, each with a single-scale and a multiscale kernel. \Cref{tab:genmod-summary} summarizes the multiscale results; the per-problem tables (both parameterizations, with gradient-norm diagnostics) and the transport visualizations are in \Cref{app:genmod-details}.

\begin{figure}
\centering
\includegraphics[width=\textwidth]{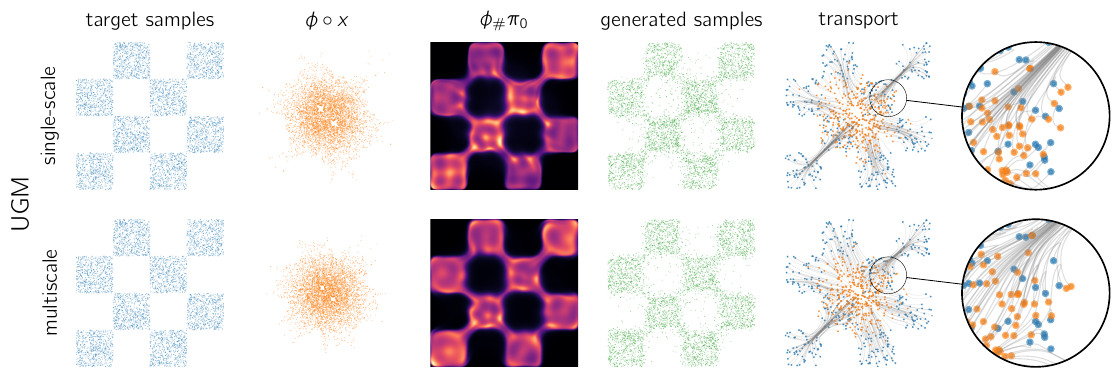}
\includegraphics[width=\textwidth]{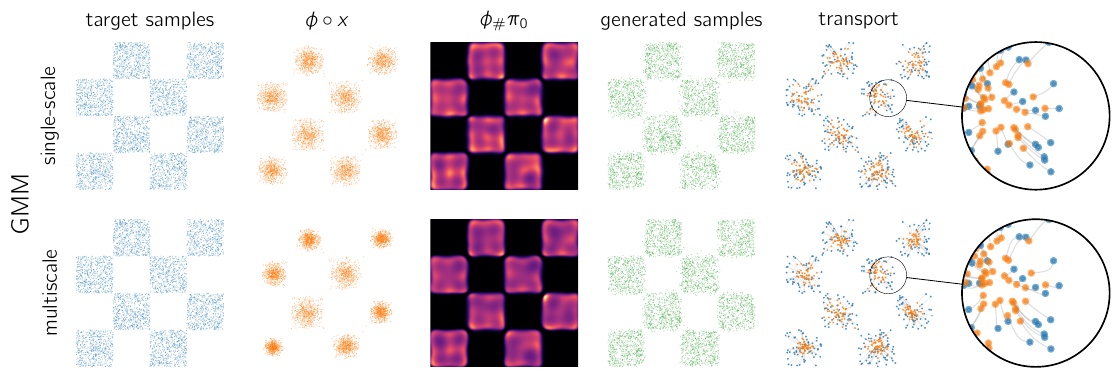}
\caption{Checkerboard test problem, UGM (top block) and GMM (bottom block) references; the two rows in each block are the single-scale and multiscale parameterizations. Columns (left to right): the target distribution, the forward-mapped samples in latent space, the learned data-space density, the generated samples (drawn from the reference and mapped by $\phi^{-1}$), and the particle transport (samples at $t=0$; trajectories and samples at $t=1$). Metrics are given in \Cref{tab:genmod-summary,tab:checkerboard-genmod}.}
\label{fig:checkerboard-varpro-gmm}
\end{figure}

The GMM reference lowers the round-trip reconstruction error $e_2$ over the UGM on the multimodal targets (on the checkerboard from $\mathrm{6.5e}\mathord{-}2$ to $\mathrm{5.2e}\mathord{-}3$, on M8G from $\mathrm{1.6e}\mathord{-}2$ to $\mathrm{5.4e}\mathord{-}5$): the trainable mixture adapts to the modes of the target, so the flow needs to deform it only slightly. At the 200-step budget the test NLL is comparable across the two references, so on these low-dimensional targets the GMM's benefit is less deformation (lower $e_2$) rather than a better density fit. We make no claim about relative convergence \emph{rates} (only endpoint objectives were logged). Solving the inverse Euler steps by fixed-point iteration (see \Cref{app:genmod-details}) gives round-trip errors below $\mathrm{1e}\mathord{-}10$ on the tested samples in every run, whereas the reported $e_2$ uses the explicit reverse-Euler inverse and so measures integration accuracy rather than invertibility. The multiscale parameterization lowers the NLL on the harder targets, and on the checkerboard the test NLL varies by at most $\mathrm{1.7e}\mathord{-}2$ for $256\le m\le1024$ (see \Cref{tab:genmod-cp-ablation}), at a per-solve cost linear in $n_s$. \Cref{fig:checkerboard-varpro-gmm} shows the checkerboard result as a representative solution.

\begin{table}
\tabadjust
\caption{Generative modeling across the four toy targets, multiscale control-point solver ($m=256$). For each target and each of the UGM and GMM references we report the training and test NLL, the final objective, the test-set reconstruction errors $e_2$ and $e_F$, and the time to solution (in seconds). Single-scale results and per-problem gradient diagnostics are in \Cref{app:genmod-details}.}
\label{tab:genmod-summary}
\begin{tabular}{llcccccc}
\toprule
 & & \multicolumn{2}{c}{\bf nll} & & \multicolumn{2}{c}{\bf reconstruction error} & \\
\cmidrule(lr){3-4}\cmidrule(lr){6-7}
\bf target & \bf reference & \bf train & \bf test & \bf objective & $e_2$ & $e_F$ & \bf time \\
\midrule
circles
& UGM & \fnum{3.260455} & \fnum{3.298113} & \fnum{663.261782} & \snum{7.874217e-03} & \snum{1.168554e-04} & \inum{6320.499757} \\
& GMM & \fnum{3.255153} & \fnum{3.296970} & \fnum{656.177509} & \snum{2.267168e-03} & \snum{3.133939e-05} & \inum{6046.077659} \\
\midrule
M8G
& UGM & \fnum{2.807694} & \fnum{2.837518} & \fnum{581.795983} & \snum{1.581085e-02} & \snum{2.573368e-04} & \inum{6065.835989} \\
& GMM & \fnum{2.818989} & \fnum{2.850337} & \fnum{564.112479} & \snum{5.401352e-05} & \snum{7.252023e-07} & \inum{5069.625399} \\
\midrule
pinwheel
& UGM & \fnum{2.311491} & \fnum{2.347938} & \fnum{478.329493} & \snum{1.401059e-02} & \snum{1.842021e-04} & \inum{6159.024467} \\
& GMM & \fnum{2.319794} & \fnum{2.354848} & \fnum{466.687110} & \snum{2.793385e-03} & \snum{4.397049e-05} & \inum{6003.647178} \\
\midrule
checkerboard
& UGM & \fnum{3.511773} & \fnum{3.533090} & \fnum{744.755502} & \snum{6.529109e-02} & \snum{1.004729e-03} & \inum{6066.399942} \\
& GMM & \fnum{3.514202} & \fnum{3.544360} & \fnum{711.237460} & \snum{5.240001e-03} & \snum{7.068113e-05} & \inum{6004.114405} \\
\bottomrule
\end{tabular}
\end{table}

\subsubsection{MNIST}
\label{sec:res-genmod-mnist}

As a real-data example, we train the control-point CNF on MNIST digits represented in a four-dimensional VAE latent space. The VAE is trained separately (its architecture and training are detailed in \Cref{app:vae}), and the CNF is fit to standardized encoder means with $m=256$ control points. We optimize on $n_s = \inum{10000}$ training latents and evaluate the NLL on \inum{10000} validation and \inum{10000} test latents held out from the fit. We use $n_t=100$, $\lambda_{\scr{loss}}=200$ (i.e., $\lambda_{\scr{loss}}/n_s=0.02$), and a multiscale Gaussian kernel whose base bandwidth is the median heuristic (0.639), giving scales $\sigma\in\{0.320,0.639,1.278\}$. We omit the scale-sparsity penalty for this run, so the objective is smooth and the CNF controls are optimized with the reduced GN-CG solver. We compare the UGM against an $8$-component GMM reference.

\begin{figure}
\centering
\includegraphics[width=\textwidth]{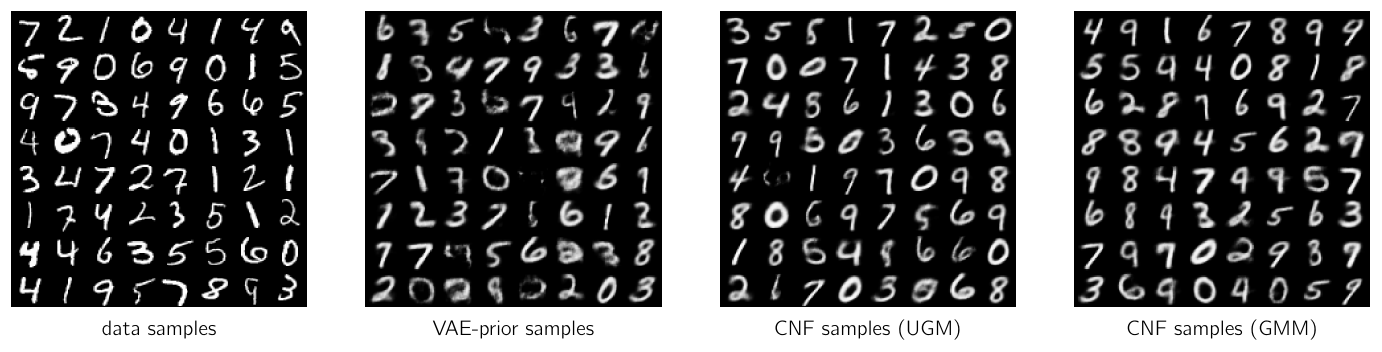}
\caption{MNIST generative modeling in a four-dimensional VAE latent space. From left to right: real MNIST digits, samples from the VAE prior, and control-point CNF samples with the UGM and VarPro GMM references. All panels are decoded by the same frozen VAE decoder; the CNF samples are drawn from the reference, mapped by the learned inverse CNF, and decoded.}
\label{fig:genmod-mnist}
\end{figure}

Solving the inverse Euler steps by fixed-point iteration recovers the tested latents with round-trip error $\sim$$\mathrm{1e}\mathord{-}12$. The UGM reaches train/validation/test NLLs of 4.94/4.96/4.86 and the GMM reference 4.79/4.83/4.72; unlike on the low-dimensional toy targets, here the GMM reference also lowers the test NLL over the UGM, and the near-equal train and test values indicate the fit generalizes without overfitting. These NLLs are computed in the standardized four-dimensional VAE latent and are not comparable to pixel-space likelihoods. Both models improve sample quality over the raw VAE prior: the energy distance~\cite{szekely2013:energy} between generated and real test latents drops from 0.068 (VAE prior) to 0.046 (UGM) and 0.044 (GMM). \Cref{fig:genmod-mnist} compares real MNIST samples, samples generated directly from the VAE prior, and CNF samples for both references. Since the VAE-prior and CNF samples share the decoder, the visual difference reflects the learned latent distribution, which sharpens the digits relative to the raw prior.

\section{Conclusions}\label{sec:conclusions}

We developed a unified optimal-control framework that solves supervised classification and likelihood-based density estimation with one RKHS-parameterized diffeomorphic flow: both tasks share the particle dynamics, control cost, and adjoint machinery, differing in the objective: a terminal classification loss, or a likelihood that also accumulates the flow's divergence along the trajectory. Because the velocity lies in an RKHS, its divergence is exact and closed-form, so density estimation avoids the stochastic trace estimators and potential-flow restriction of existing CNFs. Variable projection removes the terminal model, leaving a reduced problem whose gradient comes from the discrete adjoint and keeps its structure across tasks (density estimation only adds a divergence source term), and the learned map is a true diffeomorphism with explicitly controlled regularity.

We make the following key observations: \begin{enumerate*}[label={\it(\roman*)}]\item Under a common criterion, the reduced VarPro GN schemes need the fewest iterations of any optimizer tried (17--21 versus 58 for the closest joint baseline) and, together with joint L-BFGS, the least work; joint L-BFGS is marginally cheaper on the work axis but reaches the reference level in only 4 of 5 seeds, against 5 of 5 for the GN schemes. All converging methods reach the same solution, so the distinction is efficiency and reliability, not quality. \item For the generative task, a GMM reference markedly improves reconstruction (round-trip $e_2$) on the multimodal targets at comparable likelihood (and lowers the likelihood itself on real data), while the multiscale kernel lowers the NLL on the harder targets. \item The outer-iteration count stays bounded and the accuracy does not degrade as the problem grows, so the optimizer is never the bottleneck; the remaining per-solve cost, quadratic in the particle count, drops to linear under the control-point reduction\end{enumerate*}.

\subsection*{Disclosure of AI-assisted tools} {\itshape We used the generative AI tools Claude (Anthropic) and ChatGPT (OpenAI) to
\begin{enumerate*}[label=(\roman*)]
\item broaden our search for related literature;
\item expand unit test coverage;
\item help design and audit reproducible scripts for generating
the data underlying the figures and tables presented in this paper;
\item refactor and optimize our prototype implementations;
\item edit the manuscript for grammar, clarity, and consistency
of mathematical notation; and
\item check mathematical statements and derivations for potential errors.
\end{enumerate*}
We reviewed and verified the AI-assisted contributions and take
full responsibility for all content in this paper.}

\section*{Acknowledgments}

The work of the authors was partially supported by the National Science Foundation (NSF) under the award DMS-2145845. This material is based in part upon work supported by the NSF under Grant DMS-1929284 while AM was in residence at the Institute for Computational and Experimental Research in Mathematics in Providence, RI, during the ``Stochastic and Randomized Algorithms in Scientific Computing: Foundations and Applications'' semester program.

\appendix
\crefalias{section}{appendix}
\crefalias{subsection}{subappendix}

\section{Computing Environment}
\label{app:hardware}

The experiments reported in this article were carried out on two systems: a MacBook Pro and the Sabine cluster of the Research Computing Data Core (RCDC) at the University of Houston. The implementation uses the \texttt{PyTorch} and \texttt{PyKeOps} software stack described in \Cref{sec:methods}: the full-model kernel reductions run through the symbolic \texttt{PyKeOps} backend, and the control-point solver evaluates its cross-kernel in streamed \texttt{PyTorch}. On the GPU nodes of Sabine both run on CUDA. Concretely, the reported GPU results were produced with Python 3.14.6, \texttt{PyTorch} 2.10.0 (CUDA 12.8 build), \texttt{PyKeOps} 2.3, \texttt{NumPy} 2.5.1, and \texttt{SciPy} 1.18.0 on a single NVIDIA Tesla V100-PCIE-16GB (16\,GB) GPU; \texttt{PyKeOps} compiles its CUDA kernels with the cluster's CUDA 12.4.0 toolkit.

The laptop is an Apple MacBook Pro (Mac16,8) with an Apple M4 Pro chip (a 12-core CPU (8 performance $+$ 4 efficiency cores) and a 16-core GPU (Metal 4)) and 24\,GB of unified memory. The computational environment is Python 3.14.6, \texttt{PyTorch} 2.10.0, \texttt{PyKeOps} 2.3, \texttt{NumPy} 2.4.2, and \texttt{SciPy} 1.17.0, with the kernel reductions running CPU-only (no CUDA device).

The GPU-bound and wall-clock-timed experiments were run on the Sabine V100 nodes. All density-estimation results (the synthetic targets in single- and multiscale form, the VAE-latent MNIST example, and the control-point ablations), the real-data classification benchmarks (Pendigits and MNIST), the full-model and control-point scalability studies (see \Cref{sec:res-scalability,app:full-scalability}), and the supplementary numerical studies (see \Cref{app:adjoint-order,app:inv-consistency,app:res-diveval,app:integration-order}). The MacBook Pro was used in the toy-data classification comparisons of \Cref{sec:varpro-variants,sec:varpro-vs-baseline}. Those comparisons are reported in the hardware-independent work metric~\cref{eq:varpro-work} rather than in wall-clock time, so their conclusions do not depend on the machine they were run on.

Every experiment is fully specified by the model and solver settings reported with it: the number of Euler steps $n_t$, the kernel bandwidth(s), the control-cost and regularization weights, the CG and inner-solve tolerances, and the data-instance seeds. All datasets are publicly available: MNIST, the UCI Pen-Based Handwritten Digits, and the synthetic classification and density-estimation generators. Together with the code release, these settings regenerate every reported figure and table.

\section{Kernel Definition}\label{app:kernel-def}

We consider isotropic scalar kernels $\operatorname{ker}_\sigma(x,y)$, $x,y \in \mathbb{R}^d$, and construct the corresponding matrix-valued kernel $\kappa(x,y) = \operatorname{ker}_\sigma(x,y)I_d$, where $I_d = \operatorname{diag}(1,\ldots,1)$ is the $d\times d$ identity matrix. This isotropic structure is used throughout and enables simplified expressions in the optimality conditions. We select a Gaussian kernel~\cref{eq:gaussiankernel}. Its derivative with respect to the first argument is $\partial_1 \operatorname{ker}_\sigma(x,y) = -(x-y)\operatorname{ker}_\sigma(x,y)/\sigma^2$, which vanishes at $x=y$; this radial self-term $(\partial_1\operatorname{ker}_\sigma)(z_k,z_k)=0$ means a particle's own momentum does not contribute to the divergence~\cref{eq:exact-div} at its position. The Gaussian's rapid Fourier decay embeds the induced RKHS into $C^\infty_0(\mathcal{X},\mathbb{R}^d)$, with $\|\partial^k v_t\|_\infty \le c_k\,\sigma^{-|k|}\|v_t\|_{\mathcal{H}}$, so \Cref{asm:regularity} holds for every $q\ge1$.

For isotropic kernels, we will repeatedly use the identity
\begin{equation}
\label{eq:kernel_identity}
\langle u, \partial_1 \kappa (x,y) v \rangle
= \langle u, v\rangle\partial_1 \operatorname{ker}_\sigma(x,y),
\end{equation}

\noindent which follows directly from the structure $\kappa(x,y)=\operatorname{ker}_\sigma(x,y)I_d$.

\section{Representer Parameterization of the Velocity Field}\label{app:representer}

The velocity field $v_t$ lives in a vector-valued RKHS $\mathcal{H}$ that is continuously embedded in $C^q_0(\mathcal{X}, \mathbb{R}^d)$ for some $q \geq 1$, so that $\|v_t\|_{q,\infty} \leq c\,\|v_t\|_{\mathcal{H}}$ for a constant $c>0$. Since $\mathcal{H}$ is a space of $\mathbb{R}^d$-valued fields, point evaluation $v_t\mapsto v_t(x)$ is a bounded operator into $\mathbb{R}^d$, so by the Riesz representation theorem there is a matrix-valued reproducing kernel $\kappa : \mathcal{X}\times\mathcal{X} \to \mathbb{R}^{d\times d}$ such that $\kappa(\,\cdot\,,x)u \in \mathcal{H}$ and $\langle \kappa(\,\cdot\,,x)u, v_t\rangle_{\mathcal{H}} = \langle u, v_t(x)\rangle$ for all $x,u \in \mathbb{R}^d$ and $v_t \in \mathcal{H}$; see also~\cite[\S8.1.2]{younes2019:shapes} and~\cite{micheli2014:matrix}. We use the standard Cartesian-product construction in which the matrix-valued kernel is diagonal. Restricting further to the isotropic scalar kernel gives $\kappa(\,\cdot\,,\,\cdot\,) = \operatorname{ker}_\sigma(\,\cdot\,,\,\cdot\,) I_d$. The finite kernel expansion of the velocity over the training particles is a trajectory-dependent representer property. Fix the trajectories $\{z_i^\star(t)\}_{i=1}^{n_s}$ generated by an optimal $v^\star$ of minimal energy, i.e., a minimizer of~\cref{eq:varopt-training} with $\operatorname{dist}^2(\operatorname{id}_{\mathcal{X}},z(1))$ replaced by the kinetic energy $\int_0^1\|v_t\|^2_{\mathcal{H}}\,\mathrm{d}t$; the two forms have the same optimal value and the same optimal maps, and it is this form that we discretize. Whenever the objective depends on $v$ only through these trajectories and the control cost, as it does for the supervised loss, projecting each $v^\star_t$ onto
\[
\mathcal{V}_t = \operatorname{span}\{\kappa(\,\cdot\,,z_i^\star(t))u : u \in \mathbb{R}^d,\, i=1,\ldots,n_s\}
\]

\noindent leaves its values at the particles unchanged, since the reproducing property gives
\[
  \langle \kappa(\,\cdot\,,z_i^\star(t))u, w^\perp_t\rangle_{\mathcal{H}} = \langle u, w^\perp_t(z_i^\star(t))\rangle = 0
\]

\noindent for the orthogonal component. Write $v^\star_t = P_{\mathcal{V}_t}v^\star_t + w^\perp_t$. Two steps then close the argument. First, the projected field agrees with $v^\star$ at every particle for a.e.\ $t$, so both fields drive the same initial-value problem along the particles; since each is admissible, uniqueness of the solution of~\cref{eq:flow} makes the trajectories, and hence the data loss, identical; only the control cost differs between the two fields. Second, orthogonality gives $\|v^\star_t\|^2_{\mathcal{H}} = \|P_{\mathcal{V}_t}v^\star_t\|^2_{\mathcal{H}} + \|w^\perp_t\|^2_{\mathcal{H}}$, so the integrated control cost strictly decreases whenever $w^\perp_t \not= 0$ on a set of positive measure. A minimizer therefore cannot carry a nonzero orthogonal component, which forces
\[
v^\star_t(x) = \sum_{i=1}^{n_s} \kappa(x, z_i^\star(t))\,a_i(t) \quad \text{for a.e.}\,\, t.
\]

This finite representation is classical in diffeomorphic matching and computational anatomy, where the optimal velocity is a kernel expansion supported on the moving landmarks (a momentum representation) derived from the first-order optimality conditions~\cite{joshi2000:landmark, miller2006:geodesic, azencott2010:diffeomorphic, younes2019:shapes, mang2023:operator, zhang2021:diffeomorphic}. The orthogonal-projection form used above is a time-dependent analogue of the classical representer theorem~\cite{aronszajn1950:theory, scholkopf2001:generalized}; it is stated in this form for kernel neural ODEs in~\cite{raviola2022:kernelode}, and the two steps above give the argument directly. It requires admissible controls $v \in L^2([0,1];\mathcal{H})$ (so that the embedding $\mathcal{H}\hookrightarrow C^1_0(\mathcal{X},\mathbb{R}^d)$ gives $\operatorname{lip}(v_t)\le c\,\|v_t\|_{\mathcal{H}}\in L^1(0,1)$ and~\cref{eq:flow} has a unique solution) and the strong measurability of $t\mapsto P_{\mathcal{V}_t}v^\star_t$, which follows, for distinct particles, from the strict positive definiteness of $\operatorname{ker}_\sigma$ (which makes the Gram matrix $G(t)$ invertible, hence $t\mapsto G(t)^{-1}$ continuous) together with the continuity of $t\mapsto z_i^\star(t)$ and of the kernel sections $x\mapsto\kappa(\,\cdot\,,x)$. Pairwise distinctness of the particles $z^j_1,\ldots,z^j_{n_s}$ at every step $j=1,\ldots,n_t$, which makes $K(z)\succ0$, is not automatic after discretization: the Euler step~\cref{eq:euler-step-state} is guaranteed injective (indeed a $C^q$ diffeomorphism) when $h_t\operatorname{lip}(v^j)<1$, but need not be otherwise.

For density estimation the objective also contains the running divergence
\[
  \sum_{i=1}^{n_s}\int_0^1 \operatorname{div} v_t(z_i^\star(t))\,\mathrm{d}t,
\]

\noindent which depends on the derivatives of $v$ at the particles and is therefore not covered by the projection argument above. The divergence at a point is itself a bounded linear functional on the $C^1$ RKHS, but its representer is not a kernel section: differentiating the reproducing identity in the second argument gives $\operatorname{div} v(x) = \langle v, \eta_x\rangle_{\mathcal{H}}$ with $\eta_x(\,\cdot\,) = \partial_2\operatorname{ker}_\sigma(\,\cdot\,,x)$. Running the same projection argument on the enlarged span that also contains these derivative sections (for $q\ge1$ and $\operatorname{ker}_\sigma\in C^2$, both automatic for the Gaussian) yields an exact representer with additional radial terms, $v^\star_t(x) = \sum_i \operatorname{ker}_\sigma(x, z_i^\star(t))\bigl[a_i(t) + b_i(t)(x-z_i^\star(t))/\sigma^2\bigr]$; for the Gaussian kernel the dipole sections $\eta_{z_i^\star(t)}$ are linearly independent of the kernel sections $\operatorname{ker}_\sigma(\,\cdot\,,z_i^\star(t))$ at distinct particles (the Gram matrix of the combined family is nonsingular),
so the terms carrying $b_i$ lie outside the span $\mathcal{V}_t$ defined above.

We do not include these terms: for density estimation we adopt the finite kernel expansion~\cref{eq:rkhs-flow} as a modeling choice rather than as a consequence of optimality, so it restricts the model class, since the exact minimizer is not in it. The expansion remains a standard parameterization, and the underlying space is universal: the Gaussian RKHS is dense in $C_0(\mathcal{X},\mathbb{R}^d)$ under the relevant topology~\cite{micchelli2006:universal, wendland2005:scattered}, and multiscale mixtures inherit this universality while spanning several spatial scales.

\section{Optimality Conditions for Diffeomorphic Learning} \label{app:kkt_diffeo}

We derive the first-order optimality conditions associated with the Lagrangian~\cref{eq:lagrangian}. Throughout, we use the notation introduced there, with $n_s$ data points, $n_t$ time steps, and step size $h_t = 1/n_t$.

\subsection*{State Equations}

The state equations follow from stationarity of~\cref{eq:lagrangian} with respect to the multipliers $p = (p^j)_{j=1}^{n_t+1}$. Differentiating with respect to $p_k^{j+1}$ recovers~\cref{eq:euler-step-state} for $j=1,\ldots,n_t$, while differentiating with respect to $p_k^1$ recovers the initial condition $z_k^1 = x_k^0$.

\subsection*{Adjoint Equations}

The adjoint variables $p = (p^j)_{j=1}^{n_t+1}$ satisfy
\[
\partial_{z_k^j} \ell(z,p,a) = 0.
\]

We first derive the condition for the terminal state ($j = n_t+1$); the state $z_k^{n_t+1}$ appears only in the terminal loss and in the interval-$n_t$ constraint, whose multiplier is $p_k^{n_t+1}$. Hence
\[
p_k^{n_t+1} = - \lambda_{\scr{loss}}\,\partial_{z_k^{n_t+1}} \operatorname{loss}(c_{\Xi}(z^{n_t+1}), y).
\]

Next, we derive the optimality conditions for the intermediate states ($2 \leq j \leq n_t$). For a fixed $k$ and $j$, the state $z_k^j$ appears in:
\begin{enumerate*}[label={\it(\roman*)}]
\item the quadratic control cost at time $j$,
\item the constraint at time $j$, and
\item the constraint at time $j-1$\end{enumerate*}.

We treat these terms separately. The derivatives of the running costs are obtained by differentiating the quadratic term at time $j$ to obtain
\[
\begin{aligned}
&h_t\partial_{z_k^j}
\sum_{m=1}^{n_s}
\sum_{l=1}^{n_s}
\big\langle a_m^j, \kappa(z_m^j,z_l^j)a_l^j\big\rangle\\
&\qquad =
h_t\sum_{l=1}^{n_s}
\left(
\big\langle a_k^j,\partial_1\kappa(z_k^j,z_l^j) a_l^j\big\rangle +
\big\langle a_l^j,\partial_2\kappa(z_l^j,z_k^j)a_k^j\big\rangle
\right).
\end{aligned}
\]

The derivatives of the dynamical system at time $j$ are given by
\begin{equation*}
\begin{aligned}
&\partial_{z_k^j}
\sum_{m=1}^{n_s}
\big\langle
p_m^{j+1}, z_m^{j+1}-z_m^j-h_t\sum_{l=1}^{n_s}\kappa(z_m^j,z_l^j)a_l^j
\big\rangle
\\
&\qquad=
-p_k^{j+1}
- h_t\sum_{l=1}^{n_s}
  \big\langle p_k^{j+1},
 \partial_1\kappa(z_k^j,z_l^j)a_l^j
  \big\rangle
- h_t\sum_{m=1}^{n_s}
  \big\langle p_m^{j+1},
  \partial_2 \kappa(z_m^j,z_k^j)a_k^j
  \big\rangle.
\end{aligned}
\end{equation*}

The state $z_k^j$ also appears as the next state in the interval-$(j-1)$ constraint, whose multiplier is $p_k^j$, and therefore contributes $p_k^j$. Collecting all terms and using \cref{eq:kernel_identity} together with the symmetry of the isotropic kernel, $\partial_2\operatorname{ker}_\sigma(z_l^j,z_k^j) = \partial_1\operatorname{ker}_\sigma(z_k^j,z_l^j)$, we obtain, for $j = 2,\ldots,n_t$,
\begin{equation}
\label{eq:adjoint_recursion_full}
p_k^j
= p_k^{j+1} +
h_t \sum_{l=1}^{n_s}
\left(
\big\langle p_k^{j+1} - 2 a_k^j, a_l^j \big\rangle
+
\big\langle p_l^{j+1}, a_k^j\big\rangle
\right)
\partial_1\operatorname{ker}_\sigma(z_k^j,z_l^j).
\end{equation}

Relabeling $j \to j-1$, this is exactly the backward recursion~\cref{eq:adjoint} of the main text.

Finally, stationarity with respect to $z^1_k$, which carries the multiplier $p^1_k$ in~\cref{eq:lagrangian}, continues the recursion to $p^1_k$. Since $p^1$ does not enter the reduced gradient~\cref{eq:reduced-gradient}, the backward sweep terminates at $p_k^2$, and the adjoints $p_k^2,\ldots,p_k^{n_t+1}$ so obtained are exactly those entering~\cref{eq:reduced-gradient}.

\subsection*{Control Equation}

The optimality condition for the controls $a$ is
\[
\partial_{a_l^j}\ell(z,p,a) = 0.
\]

Differentiating the quadratic term $\langle a,Ka\rangle$ and the constraint term yields, for every $l=1,\ldots,n_s$,
\[
h_t\sum_{k=1}^{n_s}
\operatorname{ker}_\sigma(z_l^j,z_k^j)\,
\big(2 a_k^j - p_k^{j+1}\big) = 0.
\]

\noindent That is, $h_t K^j(2 a^j - p^{j+1}) = 0$. If the particles $z^j_1,\ldots,z^j_{n_s}$ are pairwise distinct, $K^j$ is symmetric positive definite (strict positive definiteness of the Gaussian kernel), and since $h_t>0$ the condition reduces to the pointwise relation $2 a_k^j - p_k^{j+1} = 0$; this is the natural reduced gradient~\cref{eq:reduced-gradient-natural} set to zero. Distinctness is not automatic after discretization; see \Cref{app:representer}.

\subsection*{Terminal Loss and Gradients}

The terminal states $z_k^{n_t+1}$ with $k=1,\ldots,n_s$ are mapped to class probabilities by the classifier $c_{\Xi}$ of~\cref{eq:loss}, with parameters $\Xi = (\xi_1,\ldots,\xi_{n_c}) \in \mathbb{R}^{n_c\times d}$ and $\xi_1 = 0$ fixing class~1 as the reference,
\begin{equation*}
c_{\Xi}(z_k^{n_t+1})
=
\frac{\exp(\Xi z_k^{n_t+1})}{\langle e_{n_c}, \exp(\Xi z_k^{n_t+1})\rangle}
\in \mathbb{R}^{n_c},
\qquad
e_{n_c} = (1,\ldots,1) \in \mathbb{R}^{n_c}.
\end{equation*}

With the one-hot labels $y_k \in \mathbb{R}^{n_c}$, the terminal cross-entropy loss and its gradient with respect to $z_k^{n_t+1}$ are
\[
\operatorname{loss}(c_{\Xi}(z^{n_t+1}), y)
= -\sum_{k=1}^{n_s} \big\langle y_k, \log c_{\Xi}(z_k^{n_t+1})\big\rangle,
\]
\noindent and
\[
\partial_{z_k^{n_t+1}}
\operatorname{loss}(c_{\Xi}(z^{n_t+1}), y)
=
\Xi^\mathsf{T}\big(c_{\Xi}(z_k^{n_t+1}) - y_k\big).
\]

With the terminal condition
\[
p_k^{n_t+1} = -\lambda_{\scr{loss}}\partial_{z_k^{n_t+1}}\operatorname{loss}(c_{\Xi}(z^{n_t+1}), y)
\]

\noindent derived above, the terminal adjoint variable is therefore
\begin{equation*}
p_k^{n_t+1}
= -\lambda_{\scr{loss}}\,\Xi^\mathsf{T}\big(c_{\Xi}(z_k^{n_t+1}) - y_k\big).
\end{equation*}

The classifier reduced gradient follows by differentiating the terminal loss and the regularizer $\lambda_{\scr{cls}}\|\Xi\|_F^2$ with respect to $\Xi$ directly, since the terminal state $z^{n_t+1}$ does not depend on $\Xi$: for $l = 2,\ldots,n_c$ (with $\xi_1 = 0$ fixed),
\begin{equation*}
\partial_{\xi_l} f
= 2\lambda_{\scr{cls}}\,\xi_l
+ \lambda_{\scr{loss}}\sum_{k=1}^{n_s}\big(c_{\Xi}(z_k^{n_t+1}) - y_k\big)_l\, z_k^{n_t+1}.
\end{equation*}

\noindent This is the second reduced gradient of \Cref{prop:varpro}; it vanishes at the VarPro inner solution $\Xi^\star(a)$.

\section{Reduced Optimality System for the Control-Point Model}\label{app:cp-kkt}

The control-point reduction of \Cref{sec:control-points} replaces the single particle block by two state blocks: the $m$ control points $c=(c_i^j)_{i,j=1}^{m,n_t}$, which carry the momentum $a$ and flow under their own $m \times m$ dynamics, and the $n_s$ data particles $z=(z_k^j)$, which are advected by the shared field~\cref{eq:cp-flow}. The control cost is the RKHS norm of $v$ on the control points, $\operatorname{cost}(a)=h_t\sum_j\langle a^j, K_{cc}^j a^j\rangle$ with $K_{cc}^j=(\operatorname{ker}_\sigma(c_i^j,c_l^j)I_d)_{i,l=1}^{m}$, while the terminal loss acts on the data particles $z^{n_t+1}$. The Lagrangian therefore carries two adjoint blocks, $q=(q_i^j)$ on the control points and $p=(p_k^j)$ on the data particles, and the same stationarity calculus that produced~\cref{eq:adjoint,eq:reduced-gradient} yields the following reduced system.

\paragraph*{Data-particle adjoint} Because $z_k$ is only advected, its adjoint recursion involves only the advection Jacobian, with the terminal condition of~\cref{eq:adjoint} on the data particles:
\begin{equation}\label{eq:cp-adjoint-p}
\begin{aligned}
p_k^{j-1} &= p_k^{j} + h_t\sum_{i=1}^{m}
\big\langle a_i^{j-1}, p_k^{j}\big\rangle\partial_1\operatorname{ker}_\sigma(z_k^{j-1}, c_i^{j-1}),
\\
p_k^{n_t+1} &= -\lambda_{\scr{loss}}\Xi^\mathsf{T}\big(c_{\Xi}(z_k^{n_t+1}) - y_k\big).
\end{aligned}
\end{equation}

\paragraph*{Control-point adjoint} The control points flow under their own $m$-point dynamics and carry the control cost, so $q$ satisfies the full-model recursion~\cref{eq:adjoint} on the $m$ control points, augmented by the cross-kernel contribution from the data particles. In addition, each control point drives the data-particle advection, adding a cross-kernel term driven by $p$:
\begin{equation}\label{eq:cp-adjoint-q}
\begin{aligned}
q_i^{j-1} &= q_i^{j} + h_t\sum_{l=1}^{m} \Big(\big\langle q_i^{j} - 2a_i^{j-1}, a_l^{j-1}\big\rangle + \big\langle q_l^{j}, a_i^{j-1}\big\rangle\Big)\,\partial_1\operatorname{ker}_\sigma(c_i^{j-1}, c_l^{j-1})\\
&\qquad + h_t\sum_{k=1}^{n_s} \big\langle a_i^{j-1}, p_k^{j}\big\rangle\,\partial_2\operatorname{ker}_\sigma(z_k^{j-1}, c_i^{j-1}),
\end{aligned}
\end{equation}

\noindent with $q^{n_t+1}_i=0$ because the control points do not enter the terminal loss. As in~\cref{eq:adjoint}, $\partial_1$ and $\partial_2$ denote differentiation in the first and second kernel argument.

\paragraph*{Reduced gradient} Collecting the $a$-derivatives of the control cost, the control flow, and the data-particle advection gives the reduced gradient in the controls,
\begin{equation}\label{eq:cp-reduced-gradient}
g_i^j = h_t\Big[\sum_{l=1}^{m}\operatorname{ker}_\sigma(c_i^j,c_l^j)\big(2a_l^j - q_l^{j+1}\big) - \sum_{k=1}^{n_s}\operatorname{ker}_\sigma(z_k^j, c_i^j)\,p_k^{j+1}\Big].
\end{equation}

\noindent This is the control-basis analogue of~\cref{eq:reduced-gradient}: the first sum reproduces the full-model $2a-p$ structure on the $m$ control points (with $q$ in place of $p$), and the second is the cross-kernel term through which the $n_s$ data particles enter, at $O(n_s m d)$ cost per step. For $m=n_s$ with $c=z$ the two adjoint blocks superpose and~\cref{eq:cp-reduced-gradient} returns the full-model gradient~\cref{eq:reduced-gradient}. The reduced GN Hessian matvec and the VarPro classifier elimination carry over unchanged, now acting on $a\in\mathbb{R}^{dmn_t}$ with the data particles entering only through the cross-kernel sums. The control-point classifier solve runs Newton's method in double precision with a tolerance of $\mathrm{1e}\mathord{-}11$ and at most $50$ steps; like the full-expansion solve of \Cref{sec:inner-solver}, it fixes $\xi_1=0$ and penalizes the remaining rows, so both solvers minimize the same inner problem~\cref{eq:inner}.

\paragraph*{Density estimation on control points} In the density-estimation instance the running divergence of \Cref{app:kkt_cnf} depends on $z_k^j$, $c_i^j$ and $a_i^j$, so, unlike in the full model, it enters all three relations: through $\partial_{z_k}$ it adds a source term to the data-particle adjoint~\cref{eq:cp-adjoint-p}, through $\partial_{c_i}$ one to the control-point adjoint~\cref{eq:cp-adjoint-q}, and through $\partial_{a_i}$ one directly to the reduced gradient~\cref{eq:cp-reduced-gradient}. Since this is the configuration used for all density-estimation experiments, we record the three terms explicitly. With the control-point velocity~\cref{eq:cp-velocity} the divergence at a data particle is
\[
\operatorname{div} v_{t^j}(z^j_k)
= \sum_{i=1}^{m}\big\langle \partial_1\operatorname{ker}_\sigma(z^j_k,c^j_i),\,a^j_i\big\rangle
= -\sum_{i=1}^{m}\operatorname{ker}_\sigma(z^j_k,c^j_i)\big\langle \Sigma^{-1}(z^j_k-c^j_i),\,a^j_i\big\rangle,
\]

\noindent so, writing $\operatorname{div}_\sigma^{\scr{cp}}(z,c,a) = -({\lambda_{\scr{loss}}h_t}/{n_s})\sum_{j,k}\operatorname{div} v_{t^j}(z^j_k)$ for its contribution to the objective, the three source terms are $-\partial_{z^{j-1}_k}\operatorname{div}_\sigma^{\scr{cp}}$ added to~\cref{eq:cp-adjoint-p}, $-\partial_{c^{j-1}_i}\operatorname{div}_\sigma^{\scr{cp}}$ added to~\cref{eq:cp-adjoint-q} (both recursions are written for level $j-1$), and $+\partial_{a^j_i}\operatorname{div}_\sigma^{\scr{cp}}$ added to~\cref{eq:cp-reduced-gradient}, with
\[
\partial_{a^j_i}\operatorname{div}_\sigma^{\scr{cp}}
= \frac{\lambda_{\scr{loss}}h_t}{n_s}\sum_{k=1}^{n_s}
\operatorname{ker}_\sigma(z^j_k,c^j_i)\,\Sigma^{-1}\big(z^j_k-c^j_i\big).
\]

\noindent With $r^j_{ki}\defeq\langle\Sigma^{-1}(z^j_k-c^j_i),a^j_i\rangle$, the two adjoint sources are
\[
\begin{aligned}
\partial_{z^j_k}\operatorname{div}_\sigma^{\scr{cp}}
&= \frac{\lambda_{\scr{loss}}h_t}{n_s}\sum_{i=1}^{m}\Big[\partial_1\operatorname{ker}_\sigma(z^j_k,c^j_i)\,r^j_{ki}+\operatorname{ker}_\sigma(z^j_k,c^j_i)\,\Sigma^{-1}a^j_i\Big],\\
\partial_{c^j_i}\operatorname{div}_\sigma^{\scr{cp}}
&= -\frac{\lambda_{\scr{loss}}h_t}{n_s}\sum_{k=1}^{n_s}\Big[\partial_1\operatorname{ker}_\sigma(z^j_k,c^j_i)\,r^j_{ki}+\operatorname{ker}_\sigma(z^j_k,c^j_i)\,\Sigma^{-1}a^j_i\Big],
\end{aligned}
\]

\noindent so that $\sum_k\partial_{z^j_k}+\sum_i\partial_{c^j_i}=0$, as translation invariance requires; the $\partial_{z_k}$ source has no $a_k$ term, since the data particles carry no momentum. The GMM-reference terms, by contrast, act only on the terminal states and so enter only through the terminal condition of~\cref{eq:cp-adjoint-p}.

\section{Optimality Conditions for CNF} \label{app:kkt_cnf}

The training objective~\cref{eq:varopt-cnf} is a KL divergence, and
\[
  \operatorname{dist}_{\scr{kl}}(\phi_\#\pi_0\,\|\,\pi_1) = \operatorname{dist}_{\scr{kl}}(\pi_0\,\|\,\phi^{\#}\pi_1) = -H(\pi_0) - \mathbb{E}_{\pi_0}[\log \phi^{\#}\pi_1];
\]

\noindent minimizing it equals minimizing the expected NLL up to the $v$-independent data-entropy constant $-H(\pi_0)$, which is dropped in~\cref{eq:cnf-continuous-objective}.

Introducing adjoint variables $p\in\mathbb{R}^{(n_t+1)dn_s}$, the discrete CNF problem~\cref{eq:cnf_discrete} has the Lagrangian
\begin{equation}
\label{eq:cnf_lagrangian}
\begin{aligned}
\ell(z,p,a) &=
f(z,a) +
\sum_{k=1}^{n_s}\big\langle p_k^1,z_k^1-x_k^0\big\rangle
\\
&\qquad +
\sum_{j=1}^{n_t}\sum_{k=1}^{n_s}
\big\langle
p_k^{j+1},
z_k^{j+1}-z_k^j-h_t\sum_{l=1}^{n_s}
\operatorname{ker}_\sigma(z_k^j,z_l^j)a_l^j
\big\rangle
\end{aligned}
\end{equation}

\noindent with $f$ defined in~\cref{eq:cnf-discrete-objective}. We derive the first-order optimality conditions for this Lagrangian, using the isotropic kernel $\operatorname{ker}_\sigma$, and the kernel identity~\cref{eq:kernel_identity} of \Cref{app:kkt_diffeo}.

The contribution of the ODE constraint in~\cref{eq:cnf_lagrangian} is identical to that of the Lagrangian for the supervised classification problem~\cref{eq:lagrangian}, and the two problems share the same control cost. Consequently, the state equation~\cref{eq:euler-step-state} and the contributions for the kinetic energy and the constraint to both the adjoint and control equations carry over from \Cref{app:kkt_diffeo}. Only two ingredients of the CNF objective~\cref{eq:cnf_discrete} differ from the classification case: the terminal cost $-(\lambda_{\scr{loss}}/n_s)\log\pi_1(z_k^{n_t+1})$ and the running divergence term
\begin{equation}
\label{eq:cnf_divergence}
\operatorname{div}_\sigma(z,a)
= h_t\frac{\lambda_{\scr{loss}}}{n_s}\sum_{j=1}^{n_t}\sum_{k=1}^{n_s}\sum_{l=1}^{n_s}
\operatorname{ker}_\sigma(z_k^j,z_l^j)\,
\big\langle \Sigma^{-1}(z_k^j-z_l^j),\, a_l^j\big\rangle.
\end{equation}

\noindent We state only the differences these two terms yield.

\subsection*{Terminal Condition}

The terminal state $z_k^{n_t+1}$ enters the objective through $-(\lambda_{\scr{loss}}/n_s)\log\pi_1(z_k^{n_t+1})$. By the terminal stationarity condition of \Cref{app:kkt_diffeo}, the terminal adjoint is
\[
p_k^{n_t+1} = (\lambda_{\scr{loss}}/n_s)\partial_z\log\pi_1(z_k^{n_t+1});
\]

\noindent that is, it is driven by the score $\partial_z\log\pi_1$ instead of the softmax residual $c_\Xi(z_k^{n_t+1}) - y_k$.

\subsection*{Divergence Source Terms}

The divergence term~\cref{eq:cnf_divergence} depends on both the states and the controls, and therefore contributes a source term to both the adjoint equation and the control equation. Adding $-\partial_{z_k^j}\operatorname{div}_\sigma(z,a)$ to~\cref{eq:adjoint_recursion_full} gives, for $j = 2,\ldots,n_t$,
\[
\begin{aligned}
p_k^j
&=
p_k^{j+1} +
h_t\sum_{l=1}^{n_s}
\big(
\langle p_k^{j+1} - 2 a_k^j, a_l^j\rangle
+ \langle p_l^{j+1}, a_k^j\rangle
\big)
\partial_1\operatorname{ker}_\sigma(z_k^j,z_l^j)
\\
&\qquad
- h_t\frac{\lambda_{\scr{loss}}}{n_s}
\partial_{z_k^j}\sum_{m=1}^{n_s}\sum_{l=1}^{n_s}
\operatorname{ker}_\sigma(z_m^j,z_l^j)
\big\langle \Sigma^{-1}(z_m^j-z_l^j), a_l^j\big\rangle.
\end{aligned}
\]

\noindent Carrying out the differentiation, the source term reads
\begin{equation}\label{eq:cnf-div-source}
\begin{aligned}
&\partial_{z_k^j}\!\!\sum_{m,l=1}^{n_s}\!\!\operatorname{ker}_\sigma(z_m^j,z_l^j)
\big\langle \Sigma^{-1}(z_m^j\!-\!z_l^j), a_l^j\big\rangle
\\
&\qquad= \sum_{l=1}^{n_s}\partial_1\operatorname{ker}_\sigma(z^j_k,z^j_l)
\Big[
\big\langle \Sigma^{-1}(z^j_k\!-\!z^j_l), a^j_l\big\rangle
+ \big\langle \Sigma^{-1}(z^j_l\!-\!z^j_k), a^j_k\big\rangle
\Big]
\\[-0.2em]
&\qquad\qquad
+ \sum_{l=1}^{n_s}\operatorname{ker}_\sigma(z^j_k,z^j_l)\,\Sigma^{-1}\big(a^j_l - a^j_k\big).
\end{aligned}
\end{equation}

Adding $\partial_{a_l^j}\operatorname{div}_\sigma(z,a)$ to the classification control equation yields
\[
h_t\sum_{k=1}^{n_s}\operatorname{ker}_\sigma(z_l^j,z_k^j)\big(2 a_k^j - p_k^{j+1}\big)
= -h_t\frac{\lambda_{\scr{loss}}}{n_s}\sum_{k=1}^{n_s}
\operatorname{ker}_\sigma(z_k^j,z_l^j)\Sigma^{-1}(z_k^j-z_l^j),
\]

\noindent for every $l = 1,\ldots,n_s$. Unlike the classification case, the right-hand side does not vanish: the divergence acts as a source term.

\subsection*{GMM Reference Gradients}

Recall the trainable GMM reference $\pi_1(\,\cdot\,;\Phi)$ with parameters $\Phi = (\eta,\mu,\zeta)$ and the terminal NLL $\operatorname{nllike}(\Phi,z)$ of \Cref{sec:gmm-reference}. Its gradients with respect to the mixture parameters are as follows: Let $\gamma_{ki} = w_i \mathcal{N}(z_k^{n_t+1}; \mu_i, \Sigma_i) \big/ \sum_{l=1}^{n_{\scr{mod}}} w_l \mathcal{N}(z_k^{n_t+1}; \mu_l, \Sigma_l)$ denote the responsibilities evaluated at the terminal states. Then the gradients are,
\begin{gather*}
\partial_{\mu_i} \operatorname{nllike}(\Phi,z)
= \frac{1}{n_s}
\sum_{k=1}^{n_s}
\gamma_{ki}\,\Sigma_i^{-1}(\mu_i - z_k^{n_t+1}),
\\
\partial_{\eta_i}\operatorname{nllike}(\Phi,z)
=
w_i - \frac{1}{n_s}\sum_{k=1}^{n_s} \gamma_{ki},
\\
\partial_{\zeta_{i,j}} \operatorname{nllike}(\Phi,z)
=
\frac{1}{2n_s}
\sum_{k=1}^{n_s}
\gamma_{ki}
\left[
1 - \frac{(z_{k,j}^{n_t+1} - \mu_{i,j})^2}{s_{i,j}}
\right]
\end{gather*}

\noindent with $s_{i,j} = \exp(\zeta_{i,j})$. These gradients can be used in a block-coordinate or joint optimization scheme together with the control variables.

\section{Multiscale Extension}
\label{app:multiscale_cnf}

The multiscale velocity model and the scale-wise sparsity penalty are introduced in \Cref{sec:multiscale}; here we record the multiscale CNF objective referenced there. With the scale-dependent controls $a_{l,\sigma}^j \in \mathbb{R}^d$ and the multiscale velocity field
\[
  v_{t^j}(x) = \sum_{\sigma\in\mathcal{S}}\sum_{l=1}^{n_s}\operatorname{ker}_\sigma(x,z_l^j)a_{l,\sigma}^j,
\]

\noindent the particle dynamics read
\begin{equation*}
z_k^{j+1}
=
z_k^j
+
h_t \sum_{\sigma\in\mathcal{S}} \sum_{l=1}^{n_s}
\operatorname{ker}_{\sigma}(z_k^j,z_l^j)a_{l,\sigma}^j,
\qquad
z_k^1 = x_k^0.
\end{equation*}

\noindent With this, the discrete CNF objective~\cref{eq:cnf_discrete} generalizes to
\[
\frac{\lambda_{\scr{loss}}}{n_s}\sum_{k=1}^{n_s}
\Big[
-\log\pi_1(z_k^{n_t+1})
-h_t\sum_{j=1}^{n_t}
\sum_{\sigma\in\mathcal{S}}\sum_{l=1}^{n_s}
\big\langle \partial_1\operatorname{ker}_\sigma(z_k^j,z_l^j), a_{l,\sigma}^j\big\rangle
\Big]
 + \operatorname{reg}(a).
\]

\noindent The regularization model is given by
\[
\operatorname{reg}(a)
= \sum_{\sigma\in\mathcal{S}}\langle a_\sigma, K_\sigma a_\sigma\rangle
+ h_t\sum_{\sigma\in\mathcal{S}}\lambda_{\scr{sp}}^{\sigma}
\sum_{j=1}^{n_t}\sum_{l=1}^{n_s}
\log\big(1+\varepsilon^{-1}\|a_{l,\sigma}^j\|_2\big),
\]

\noindent where $K_\sigma$ is the scale-dependent kernel operator of \Cref{sec:multiscale}. Notice that by \Cref{lem:exact-div}, the closed-form divergence is additive across scales. The state and adjoint equations of \Cref{app:kkt_cnf} carry over unchanged apart from the summation over scales. The control equation also acquires the derivative of the group-log penalty. Write $g^j_{l,\sigma}$ for the gradient of the smooth part (data misfit, divergence, and kinetic energy) with respect to $a^j_{l,\sigma}$, and $c_\sigma\defeq h_t\lambda_{\scr{sp}}^\sigma$. At a block with $a^j_{l,\sigma}\neq0$ stationarity reads
\[
g^j_{l,\sigma} + c_\sigma\,\frac{a^j_{l,\sigma}}{\|a^j_{l,\sigma}\|_2\,\bigl(\varepsilon+\|a^j_{l,\sigma}\|_2\bigr)} = 0.
\]
At a block with $a^j_{l,\sigma}=0$, where the penalty is not differentiable, it is the inclusion $0\in g^j_{l,\sigma}+(c_\sigma/\varepsilon)\,\mathbb{B}$ in the (Fr\'echet) subdifferential, with $\mathbb{B}$ the closed unit ball of $\mathbb{R}^d$; that is, $\|g^j_{l,\sigma}\|_2\le c_\sigma/\varepsilon$.

\section{Algorithms}\label{app:algorithms}

Here, we collect pseudocode for the optimization schemes of \Cref{sec:baseline-opt,sec:varpro}.

\Cref{alg:iterative-scheme} is the joint solver underlying the baseline optimizers of \Cref{tab:optimization-methods}: each iteration solves the state~\cref{eq:euler-step-state} and adjoint~\cref{eq:adjoint} equations, evaluates the reduced gradient in the controls together with the classifier gradient, forms a joint search direction $(s_a, s_\Xi)$, and updates the controls $a$ and the classifier $\Xi$ simultaneously after a line search; the baselines differ only in how this direction is built.

\Cref{alg:varpro-step} is the generic VarPro iteration. It eliminates the classifier by an inner solve for $\Xi^\star(a)$ of the convex subproblem~\cref{eq:inner} at the current terminal state, evaluates the adjoint and the reduced gradient at $\Xi^\star(a)$ and updates only the controls $a$, setting $\Xi \gets \Xi^\star(a)$.

\Cref{alg:varpro-gncg,alg:varpro-trgn} give the two reduced GN outer steps we advocate, which share the forward/inner/adjoint solves and the matrix-free GN Hessian action and differ only in their globalization. \Cref{alg:varpro-gncg} (GN-CG) solves the reduced system $H_{\scr{gn}} s_a = -g_a$ inexactly by CG. It falls back to GD if $s_a$ is not a descent direction. The scheme is globalized using an Armijo line-search. \Cref{alg:varpro-trgn} (TR-GN) instead solves the trust-region subproblem~\cref{eq:tr-subproblem} by Steihaug--Toint truncated CG, and accepts or rejects the step, adjusting the radius $\Delta$, according to the ratio $\rho$ of actual to predicted decrease.

\Cref{alg:varpro-proxgn} is the outer step for the composite multiscale objective of \Cref{app:varpro-variants}: it keeps the smooth GN model but treats the nonsmooth group-log penalty through its proximal map, solving the composite GN model~\cref{eq:prox-gn-model} by an accelerated proximal-gradient inner iteration~\cite{beck2009:fast}. The method is globalized using an Armijo line search that falls back to a proximal-gradient step. The block-wise group-log proximal map and the constants are given in \Cref{app:varpro-variants}.

\Cref{alg:cnf-varpro-gmm} instantiates the VarPro iteration for the density-estimation problem~\cref{eq:cnf_discrete} with a trainable GMM reference eliminated by variable projection. It mirrors \Cref{alg:varpro-step}: the terminal loss is replaced by the CNF likelihood and the inner classifier solve by an EM re-fit $\Phi^\star(a)$ of the mixture to the current terminal states $\{z_k^{n_t+1}\}_{k=1}^{n_s}$, after which the outer control step is taken by GN-CG (single-scale) or proximal GN (multiscale). The outer step uses \Cref{alg:varpro-gncg} or \Cref{alg:varpro-proxgn}, executed in the Euclidean inner product on the Euclidean GN operator (the full-expansion solver uses the natural inner product).

\begin{algorithm}
\algadjust
\caption{Joint optimization of the controls $a$ and classifier $\Xi$ for the discrete supervised problem~\cref{eq:varopt-learning-discrete}.}
\label{alg:iterative-scheme}
\begin{algorithmic}[1]
\State \textbf{input:}
  $n_{\scr{iter}} > 0$; initial guesses $\Xi_{\scr{init}}$,
  $a_{\scr{init}}$; training data $\Psi = (x^0,y)$
  \State \textbf{initialize}: $\Xi^{(1)} \gets \Xi_{\scr{init}}$, $a^{(1)}
  \gets a_{\scr{init}}$, $k \gets 1$, $\text{stop} \gets 0$
  \While{$k < n_{\scr{iter}} \land \neg$ stop}
  \State $z^{(k)} \gets$ solve state equation~\cref{eq:euler-step-state} with initial condition $x^0$ given candidate $a^{(k)}$
  \State $p^{(k)} \gets$ solve adjoint equation~\cref{eq:adjoint} given $\{\Xi^{(k)}, a^{(k)}, z^{(k)}\}$
  \State $g_a^{(k)} \gets$ evaluate reduced gradient~\cref{eq:reduced-gradient} given $\{a^{(k)}, z^{(k)}, p^{(k)}\}$
  \State $g_{\Xi}^{(k)} \gets$ evaluate classifier gradient given $\{\Xi^{(k)}, z^{(k)}, y\}$
  \State $(s_a^{(k)}, s_\Xi^{(k)}) \gets$ compute search direction given $(g_a^{(k)}, g_{\Xi}^{(k)})$
  \State $\gamma^{(k)} \gets$ determine step size via line search
  \State $a^{(k+1)} \gets a^{(k)} + \gamma^{(k)} s_a^{(k)}, \quad \Xi^{(k+1)} \gets \Xi^{(k)} + \gamma^{(k)} s_\Xi^{(k)}$
  \State $k \gets k + 1$
  \State stop $\gets$ check convergence
  \EndWhile
\State \Return $\Xi_{\scr{opt}} \gets \Xi^{(k)} $, $a_{\scr{opt}} \gets a^{(k)}$
\end{algorithmic}
\end{algorithm}

\begin{algorithm}
\algadjust
\caption{VarPro iteration for solving \cref{eq:varopt-learning-discrete}. The classifier is not carried between outer steps: it is a function of the control alone and is recomputed as $\Xi^\star(a^{(k)})$ in line~4 of every call, so the outer state is $a$ only and the returned solution pair is $(a_{\scr{opt}},\Xi^\star(a_{\scr{opt}}))$.}
\label{alg:varpro-step}
\begin{algorithmic}[1]
\State \textbf{input}: iterate $a^{(k)}$; training data $\Psi = (x^0, y)$
\State $z^{(k)} \gets$ solve state equation \cref{eq:euler-step-state} with initial condition $x^0$ given candidate $a^{(k)}$
\State evaluate control cost $\operatorname{cost}(z^{(k)},a^{(k)}) = \langle a^{(k)}, K(z^{(k)})\,a^{(k)}\rangle$
\State $\Xi^\star(a^{(k)}) \gets$ solve~\cref{eq:inner} with $z^{n_t+1}$ fixed
\State $p^{(k)} \gets$ solve adjoint equation~\cref{eq:adjoint} given $\{\Xi^\star(a^{(k)}), a^{(k)}, z^{(k)}\}$
\State $g_a^{(k)} \gets$ evaluate reduced gradient~\cref{eq:reduced-gradient} (or \cref{eq:reduced-gradient-natural}) given $\{a^{(k)}, z^{(k)}, p^{(k)}\}$
\State $s_a^{(k)} \gets$ compute search direction given $g_a^{(k)}$
\State $\gamma^{(k)} \gets$ perform line search
\State $a^{(k+1)} \gets a^{(k)} + \gamma^{(k)} s_a^{(k)}$
\State stop $\gets$ check for convergence
\State $k \gets k + 1$
\State \Return $\{a^{(k)}, \text{stop}\}$
\end{algorithmic}
\end{algorithm}

\begin{algorithm}
\algadjust
\caption{VarPro outer step with GN-CG as outer solver. The parameters of the algorithm are the tolerance $\operatorname{tol}_{\scr{cg}} > 0$ (forcing term) and maximum number of iterations $m_{\scr{cg}} \in \mathbb{N}$.}
\label{alg:varpro-gncg}
\begin{algorithmic}[1]
\State \textbf{input:} iterate $a^{(k)}$; data $\Psi = (x^0, y)$
\State $z^{(k)}\gets$ forward solve; $\Xi^\star(a^{(k)})\gets$ inner solve; $p^{(k)}\gets$ adjoint solve
\State $g_a^{(k)}\gets$ evaluate reduced gradient
\State build matrix-free matvec $\tilde{a}\mapsto H_{\scr{gn}}^{\scr{nat}}(a^{(k)})\tilde{a}$ at the frozen trajectory
\State $s_a^{(k)} \gets$ solve $H_{\scr{gn}}^{\scr{nat}}(a^{(k)})s_a = -g_a^{\scr{nat},(k)}$
\If{$\langle g_a^{(k)}, s_a^{(k)}\rangle \geq 0$}
  \State $s_a^{(k)} \gets -g_a^{\scr{nat},(k)}$
\EndIf
\State $\gamma^{(k)} \gets$ perform line search on $f(a^{(k)} + \gamma s_a^{(k)},\Xi^\star(a^{(k)}))$
\State $a^{(k+1)} \gets a^{(k)} + \gamma^{(k)} s_a^{(k)}$
\end{algorithmic}
\end{algorithm}

\begin{algorithm}
\algadjust
\caption{VarPro outer step with TR-GN as outer solver. The parameters of the algorithm are the thresholds $\eta_1 > 0$, $\eta_2>0$, the shrink/expand factors $\gamma_1>0$, $\gamma_2 >0$ and the bounds $\Delta_{\scr{min}}>0$, $\Delta_{\scr{max}}>0$.}
\label{alg:varpro-trgn}
\begin{algorithmic}[1]
\State \textbf{input:} iterate $a^{(k)}$; data $\Psi = (x^0, y)$; current radius $\Delta^{(k)}$
\State $z^{(k)}\gets$ forward solve; $\Xi^\star(a^{(k)})\gets$ inner solve; $p^{(k)}\gets$ adjoint solve
\State $g_a^{(k)}\gets$ evaluate reduced gradient
\State build matrix-free matvec $\tilde{a}\mapsto H_{\scr{gn}}^{\scr{nat}}(a^{(k)})\tilde{a}$ at the frozen trajectory
\State $s_a^{(k)},\,\text{on\_boundary} \gets$ solve~\cref{eq:tr-subproblem} with step truncated at TR boundary
\State $\rho^{(k)} \gets (\widehat{f}(a^{(k)}) - \widehat{f}(a^{(k)} + s_a^{(k)})) / (-m_k(s_a^{(k)}))$
\If{$\rho^{(k)} \geq \eta_1$}
  \State $a^{(k+1)} \gets a^{(k)} + s_a^{(k)}$
  \If{$\rho^{(k)} \geq \eta_2 \land \text{on\_boundary}$}
    \State $\Delta^{(k+1)} \gets \min(\gamma_2 \Delta^{(k)},\Delta_{\scr{max}})$
  \Else
    \State $\Delta^{(k+1)} \gets \Delta^{(k)}$
  \EndIf
\Else
  \State $a^{(k+1)} \gets a^{(k)}$
  \State $\Delta^{(k+1)} \gets \max(\gamma_1 \|s_a^{(k)}\|_{\scr{nat}},\Delta_{\scr{min}})$
\EndIf
\end{algorithmic}
\end{algorithm}

\begin{algorithm}
\algadjust
\caption{VarPro outer step with proximal GN as outer solver, for the composite (multiscale) objective $\widetilde{f}=\nu+\psi$~\cref{eq:prox-gn-model}. The parameters are the number of inner iterations $m_{\scr{fista}}\in\mathbb{N}$, the Armijo constant $c_1>0$, and the backtracking factor $\beta\in(0,1)$. The proximal map $\operatorname{prox}_{\psi/L}$ decouples into the block-wise group-log map of \Cref{app:varpro-variants}.}
\label{alg:varpro-proxgn}
\begin{algorithmic}[1]
\State \textbf{input:} iterate $a^{(k)}$; reference $\pi_1$; data $\Psi = (x^0)$
\State $z^{(k)}\gets$ forward solve
\State $\Phi^\star(a^{(k)})\gets$ inner solve (if the reference is trained)
\State $p^{(k)}\gets$ adjoint solve
\State $g_a^{(k)}\gets \partial_a\nu(a^{(k)})$
\State build matrix-free matvec $\tilde{a}\mapsto H_{\scr{gn}}^{(k)}\tilde{a}$ at the frozen trajectory
\State $L\gets$ estimate $\|H_{\scr{gn}}^{(k)}\|$ by power iteration
\State $s\gets 0$;\quad $y\gets 0$;\quad $\tau\gets 1$
\For{$i=1,\dots,m_{\scr{fista}}$}
  \State $w\gets y - L^{-1}\big(g_a^{(k)} + H_{\scr{gn}}^{(k)} y\big)$
  \State $s^{+}\gets \operatorname{prox}_{\psi/L}\!\big(a^{(k)}+w\big) - a^{(k)}$
  \State $\tau^{+}\gets \tfrac12\big(1+\sqrt{1+4\tau^2}\big)$;\quad $y\gets s^{+} + \tfrac{\tau-1}{\tau^{+}}\big(s^{+}-s\big)$
  \State $s\gets s^{+}$;\quad $\tau\gets \tau^{+}$
\EndFor
\State $s_a^{(k)}\gets s$
\State $\gamma^{(k)}\gets$ Armijo backtracking on $\widetilde{f}$
\If{no sufficient decrease}
  \State $s_a^{(k)}\gets \operatorname{prox}_{\psi/L}\!\big(a^{(k)}-L^{-1}g_a^{(k)}\big) - a^{(k)}$;\quad $\gamma^{(k)}\gets$ backtrack on $\widetilde{f}$ (else $0$)
\EndIf
\State $a^{(k+1)}\gets a^{(k)} + \gamma^{(k)} s_a^{(k)}$
\end{algorithmic}
\end{algorithm}

\begin{algorithm}
\algadjust
\caption{VarPro CNF iteration for the density-estimation problem~\cref{eq:cnf_discrete} with a GMM reference eliminated by an inner EM re-fit. The outer control step is GN-CG (see \Cref{alg:varpro-gncg}) for the single-scale objective and proximal GN (see \Cref{alg:varpro-proxgn}) for the multiscale objective. The parameters are the number of EM sweeps $m_{\scr{em}}\in\mathbb{N}$ and the variance floor $s_{\min}>0$ ($s_{\min}=\mathrm{1e}\mathord{-}6$). The covariances are diagonal, $\Sigma_i=\operatorname{diag}(s_{i,1},\ldots,s_{i,d})$ with $s_{i,r}=\exp(\zeta_{i,r})$, as in \Cref{sec:gmm-reference}.}
\label{alg:cnf-varpro-gmm}
\begin{algorithmic}[1]
\State \textbf{input:} controls $a^{(k)}$; data $\Psi = (x^0)$
\State $z^{(k)}\gets$ solve~\cref{eq:euler-step-state} from $x^0$ given $a^{(k)}$, accumulating the running divergence for the log-likelihood
\State $\Phi^\star(a^{(k)})\gets \Call{EM-REFIT}{\{z_k^{n_t+1}\}_{k=1}^{n_s}}$
\State $p^{(k)}\gets$ solve the CNF adjoint with terminal GMM score $\partial_z\log\pi_1(\,\cdot\,;\Phi^\star)$
\State $g_a^{(k)}\gets$ evaluate the CNF reduced gradient given $\{a^{(k)}, z^{(k)}, p^{(k)}\}$
\State $a^{(k+1)}\gets$ outer step from $\{a^{(k)}, g_a^{(k)}, H_{\scr{gn}}^{(k)}\}$
\State $k\gets k+1$;\quad stop $\gets$ check convergence
\State \Return $\{a^{(k)}, \text{stop}\}$ (form $\Phi^\star(a_{\scr{opt}})$ by a final re-fit at the returned controls)
\Statex
\Procedure{EM-REFIT}{$\{z_k\}_{k=1}^{n_s}$}
  \State $\mu\gets k\text{-means}(\{z_k\})$ (fixed seed)
  \State $w_i\gets 1/n_{\scr{mod}}$
  \State $s_{i,r}\gets\max\bigl(\operatorname{Var}_r\{z_k\},s_{\min}\bigr)$, $r=1,\ldots,d$
  \For{$t=1,\dots,m_{\scr{em}}$}
    \State $\gamma_{ki}\gets w_i\mathcal{N}(z_k;\mu_i,\Sigma_i)\big/\textstyle\sum_{l}w_l\mathcal{N}(z_k;\mu_l,\Sigma_l)$,\quad $\Sigma_i=\operatorname{diag}(s_i)$
    \State $N_i\gets \textstyle\sum_{k}\gamma_{ki}$
    \State $w_i\gets N_i/n_s$;\quad $\mu_i\gets N_i^{-1}\textstyle\sum_k\gamma_{ki}z_k$
    \State $s_{i,r}\gets N_i^{-1}\textstyle\sum_k\gamma_{ki}(z_{k,r}-\mu_{i,r})^2+s_{\min}$, $r=1,\ldots,d$
  \EndFor
  \State \Return $\Phi=(w,\mu,s)$
\EndProcedure
\end{algorithmic}
\end{algorithm}

Because the re-fit is cold-started from a fixed-seed $k$-means~\cite{bishop2006:pattern} initialization rather than warm-started from $\Phi^{(k)}$, the solution map $\Phi^\star(a)$ need not be continuous in $a$, consistent with the caveat of \Cref{rem:varpro-density} that the GMM inner problem lies outside the smooth-VarPro theory. Outer convergence is nonetheless stable in practice (see \Cref{sec:syn-data-dist}).

\section{Additional Derivations and Algorithmic Details}
\label{app:methods-details}

Here, we provide some additional details regarding the algorithms we consider in the main manuscript.

\subsection{Baseline Search Directions}
\label{app:baseline-directions}

The joint and VarPro schemes share the search-direction template of \Cref{sec:baseline-opt}: given a gradient $g^{(k)}$, a direction $s^{(k)}$ is produced and the iterate is updated with the direction-appropriate line search of \Cref{tab:optimization-methods}. We give here the three direction families listed in \Cref{tab:optimization-methods}.

\subsubsection*{Gradient Descent Variants}
The Euclidean GD direction is $s^{(k)} = -g^{(k)}$. Under an exact line search successive gradients, and hence successive search directions, are orthogonal. The resulting zigzag path slows plain GD on ill-conditioned problems. The natural-gradient variant NGD preconditions by the kernel-induced metric of the RKHS~\cite{younes2019:shapes}. This corresponds to descent in the Riemannian geometry of $\mathcal{H}^d_\sigma$ and improves the conditioning of the linearized problem; see \cref{eq:reduced-gradient-natural}. These directions correspond to the GD/NGD rows of \Cref{tab:optimization-methods}. We also consider Nesterov's accelerated GD method~\cite{nesterov2004:introductory}, which augments the gradient step with a momentum term extrapolated from the two previous iterates; applied to the Euclidean and to NGD it yields the NAGD and NANGD rows of \Cref{tab:optimization-methods}. As a second acceleration strategy we consider Anderson acceleration of the joint steepest-descent fixed-point map $\chi(\Theta^{(k)}) = \Theta^{(k)} + \gamma^{(k)} s^{(k)}$, using the same safeguarded mixing scheme that we describe for the VarPro outer step in \Cref{sec:outer-step}. Here the mixing acts on the concatenated iterate $(a, \Xi)$, so the residual couples the control and classifier blocks in a single Euclidean norm. Applied to the Euclidean and to the natural control gradient this yields the AAGD and AANGD rows of \Cref{tab:optimization-methods}.

\subsubsection*{Nonlinear Conjugate Gradient}
Nonlinear CG generates directions $s^{(k+1)} = -g^{(k+1)} + \beta^{(k+1)} s^{(k)}$ that are conjugate with respect to a quadratic model~\cite[Chap.~5]{nocedal2006:numerical}. We use the Fletcher--Reeves and Polak--Ribi\`ere coefficients $\beta_{\scr{fr}}$, $\beta_{\scr{pr}}$ with the standard $\beta^{(k+1)} = \max(\beta_{\scr{pr}}^{(k+1)}, 0)$ safeguard for global convergence~\cite[Chap.~5]{nocedal2006:numerical}. The natural-gradient variants (NFRCG, NPRCG of \Cref{tab:optimization-methods}) evaluate these coefficients in the kernel-induced inner product $\langle u, v\rangle_{\scr{nat}} = \langle u, K v\rangle$ rather than the Euclidean one, so each pairing couples the Euclidean gradient $g$ with the conditioned gradient $g^{\scr{nat}} = K^{-1}g$ (e.g.\ the Fletcher--Reeves numerator becomes $g^{(k+1)\mathsf{T}} K^{-1} g^{(k+1)}$); this combines the nonlinear CG of~\cite[Sec.~5.2]{nocedal2006:numerical} with the preconditioning of~\cite[Sec.~5.1]{nocedal2006:numerical}. Using the Euclidean pairing of the natural gradients would instead misweight the directions and degrade conjugacy.

\subsubsection*{Limited-Memory BFGS (L-BFGS)}

We use the two-loop L-BFGS recursion of~\cite[Alg.~7.4]{nocedal2006:numerical} on the $m=10$ most recent secant pairs $(s^{(i)}, y^{(i)}) = (\Theta^{(i+1)} - \Theta^{(i)},\, g^{(i+1)} - g^{(i)})$, with the standard self-scaled initial inverse-Hessian and curvature filtering (pairs with $\langle y^{(i)}, s^{(i)}\rangle > 0$ are kept); a non-descent direction clears the history and falls back to GD. This is the outer step of the VarPro L-BFGS variant of \Cref{sec:outer-step}. As a baseline joint solver we use the compiled L-BFGS-B routine of \texttt{SciPy} ($m=10$), which realizes the same recursion with its own curvature skipping and line search.

\subsection{Proof of the VarPro Equivalence (\Cref{prop:varpro})}
\label{app:varpro-proof}

The inner objective in~\cref{eq:inner} is the sum of a convex function of $\Xi$ (the softmax cross-entropy~\cref{eq:loss} with fixed $z^{n_t+1}$) and the $2\lambda_{\scr{cls}}$-strongly convex regularizer $\lambda_{\scr{cls}}\|\Xi\|_F^2$; it is therefore $2\lambda_{\scr{cls}}$-strongly convex with a unique minimizer $\Xi^\star(a)$, giving part \ref{itm:varpro-i}. Smoothness of $a \mapsto \Xi^\star(a)$ follows from the implicit function theorem applied to the optimality condition $\partial_{\bar\Xi} f(z^\star(a), a, \Xi) = 0$, whose Jacobian in $\bar\Xi$ is the positive-definite regularized Hessian of the inner objective. Throughout this proof, $\partial_{\bar\Xi}$ denotes differentiation in the free rows $\bar\Xi = (\xi_2,\ldots,\xi_{n_c})$, i.e., along the feasible set $\mathcal{W} = \{\Xi : \xi_1 = 0\}$ of~\cref{eq:inner}, over which the joint problem~\cref{eq:varopt-learning-discrete} is posed. Correspondingly, $[\partial^2_{\bar\Xi\bar\Xi}f]^{-1}$ is the inverse of the principal submatrix of the full $\Xi\Xi$ Hessian, which in general differs from the corresponding submatrix of that Hessian's inverse. For part \ref{itm:varpro-ii}, write $F(a,\Xi) \defeq f(z^\star(a), a, \Xi)$ for the objective with the state eliminated, so that $\widehat{f}(a) = F(a, \Xi^\star(a)) = \min_{\Xi\,\in\,\mathcal{W}} F(a,\Xi)$ by part \ref{itm:varpro-i}. The standard variable-projection results give the stationary-point equivalence~\cite{golub1973:differentiation,golub2003:separable}; the local-minimizer equivalence asserted here also uses the strong convexity of the inner problem and the continuity of $\Xi^\star$ from part \ref{itm:varpro-i}.
($\Leftarrow$) Suppose $\Xi^\dag = \Xi^\star(a^\dag)$ and $a^\dag$ minimizes $\widehat{f}$ on a neighborhood $U$. Then for every $a \in U$ and every $\Xi\in\mathcal{W}$,
\[
F(a,\Xi) \ge \min_{\Xi'\in\,\mathcal{W}} F(a,\Xi') = \widehat{f}(a) \ge \widehat{f}(a^\dag) = F(a^\dag,\Xi^\dag),
\]

\noindent so $(a^\dag,\Xi^\dag)$ is a local minimizer of the joint problem (global in $\Xi$ over $\mathcal{W}$); no continuity is needed for this direction.
($\Rightarrow$) Suppose $(a^\dag,\Xi^\dag)$ minimizes the joint problem on a neighborhood $V = U \times N$, with $N$ a relative neighborhood of $\Xi^\dag$ in $\mathcal{W}$. Fixing $a = a^\dag$ makes $\Xi^\dag$ a local minimizer of the strongly convex map $\Xi \mapsto F(a^\dag,\Xi)$, hence its unique global minimizer on $\mathcal{W}$, so $\Xi^\dag = \Xi^\star(a^\dag)$. Because $\Xi^\star$ is continuous by part \ref{itm:varpro-i}, there is a neighborhood $U' \subseteq U$ of $a^\dag$ with $\Xi^\star(U') \subseteq N$; for $a \in U'$ we then have $\widehat{f}(a) = F(a,\Xi^\star(a)) \ge F(a^\dag,\Xi^\dag) = \widehat{f}(a^\dag)$, so $a^\dag$ locally minimizes $\widehat{f}$. Part \ref{itm:varpro-iii} is the envelope theorem~\cite{golub1973:differentiation,golub2003:separable}: differentiating $\widehat{f}(a) = f(z^\star(a), a, \Xi^\star(a))$ by the chain rule,
\[
\partial_a \widehat{f}(a)
= \partial_a f\big|_{\Xi = \Xi^\star(a)}
+ \partial_{\bar\Xi} f\big|_{\Xi = \Xi^\star(a)}\,\mathrm{d}_a \bar\Xi^\star(a),
\]

\noindent and the second term vanishes because $\partial_{\bar\Xi} f|_{\Xi = \Xi^\star(a)} = 0$ by the first-order optimality condition of \ref{itm:varpro-i}. Note that $\partial_a f$ denotes the partial derivative through both the explicit $a$-dependence and through $z^\star(a)$; the latter contribution is exactly what is computed by the discrete adjoint reduced gradient~\cref{eq:reduced-gradient}. This completes the proof.

\begin{remark}[Relation to alternating minimization]\label[remark]{rem:altmin}
Alternating minimization updates $a$ with $\Xi$ frozen and then $\Xi$ with $a$ frozen. Solving the inner problem to optimality makes VarPro's reduced gradient coincide with this frozen-$\Xi$ gradient by the envelope identity~\cref{eq:envelope}. Since our line search and GN model also hold $\Xi$ at $\Xi^\star(a^{(k)})$, one VarPro GN-CG step is one sweep of alternating minimization with a GN update in $a$ and an exact update in $\Xi$. Eliminating the classifier and choosing the outer curvature model are therefore separate decisions, and both reduce the outer-iteration count: on the reduced problem GN-CG needs 17 iterations against 51 for NGD (see \Cref{tab:varpro-variants}), and L-BFGS needs 28 on the reduced problem against 58 jointly (see \Cref{tab:varpro-full-baseline}). A single inner iteration loses the envelope identity and yields a plain alternating scheme.
\end{remark}

\subsection{VarPro Outer-Direction Variants}
\label{app:varpro-variants}

We provide the full instantiation of the VarPro outer-direction variants summarized in \Cref{sec:varpro}, including the reduced GN Hessian, its matrix-free action, and the TR model.

We consider several VarPro outer-direction variants. Steepest descent and L-BFGS share the generic descriptions in \Cref{sec:baseline-opt}; we record only the VarPro-specific choices here, summarizing for each variant its search direction, cost, and globalization strategy.

\subsubsection*{Gradient Descent}
For gradient descent we set $s_a^{(k)} = -g_a^{(k)}$ (Euclidean) or $s_a^{(k)} = -g_a^{\scr{nat},(k)}$ (natural). The natural-gradient variant is invariant under the kernel-induced metric of $\mathcal{H}^d_\sigma$ and has no extra per-iteration cost.

\subsubsection*{L-BFGS}

We apply the two-loop L-BFGS direction to the reduced gradient $g_a$ with sliding window $m_{\scr{out}} = 10$. With the natural gradient, $y^{(k)}=g^{\scr{nat},(k+1)}-g^{\scr{nat},(k)}=K(a^{(k+1)})^{-1}g^{(k+1)}-K(a^{(k)})^{-1}g^{(k)}$ is a difference of gradients taken in two different metrics, since $K$ depends on $a$ through the trajectory, so the secant pairs carry no consistent curvature information about $\widehat f$. We therefore use the Euclidean variant.

\subsubsection*{Polak--Ribi\`ere Nonlinear CG}

We apply the conditioned Polak--Ribi\`ere nonlinear CG ($s_a^{(k+1)} = -g_a^{\scr{nat},(k+1)} + \beta^{(k+1)} s_a^{(k)}$, with the safeguarded natural-metric coefficient of \Cref{app:baseline-directions}) to the reduced gradient, restarting to $-g_a^{\scr{nat},(k+1)}$ on a non-descent direction or a failed strong-Wolfe line search. This is the conditioned CG originally used for diffeomorphic learning~\cite{younes2020:diffeomorphic}, applied here to the reduced objective $\widehat{f}$ after variable projection has eliminated $\Xi$ rather than to the joint $(a,\Xi)$ iterate; \Cref{sec:varpro-vs-baseline} compares the two. Beyond the reduced gradient it incurs no additional per-iteration cost.

\subsubsection*{Anderson Acceleration}

Anderson acceleration (AA)~\cite{anderson1965:iterative,walker2011:anderson,he2025:ngmres} extrapolates from a sliding window (size $m_{\scr{aa}} = 5$) of residuals $r^{(k)} = \chi(a^{(k)}) - a^{(k)}$ of the GD fixed-point map $\chi(a) = a + \gamma s_a$, obtaining the mixing coefficients from the standard Tikhonov-regularized least-squares over pairwise residual differences ($\lambda_{\scr{aa}} = \mathrm{1e}\mathord{-}8$)~\cite{walker2011:anderson}. Because AA is not a monotone-descent method, we accept the extrapolated iterate $a_{\scr{aa}}^{(k)}$ only under two safeguards~\cite{zhang2020:globally}:
\begin{enumerate*}[label={\it(\roman*)}]
\item \emph{step-length pre-filter:} $\|a_{\scr{aa}}^{(k)} - a^{(k)}\|_2 \leq c_{\scr{aa}}\|r^{(k)}\|_2$ with $c_{\scr{aa}} = 3$;
\item \emph{descent check:} one extra forward pass to verify $\widehat{f}(a_{\scr{aa}}^{(k)}) \leq \widehat{f}(\chi(a^{(k)}))$.
\end{enumerate*}
On failure we fall back to $\chi(a^{(k)})$. The window is reset whenever the inner line search fails.

\subsubsection*{GN-CG Method}

We consider a GN-CG method. At the current inner optimum $\Xi^\star(a)$, the Hessian of the state-eliminated objective $F(a,\Xi)$ (with $\Xi$ frozen at $\Xi^\star(a)$) in the controls $a$ is
\[
\partial^2_{aa} F(a,\Xi^\star(a)) = 2 K(a) + \lambda_{\scr{loss}} J_a^\mathsf{T}H_{\scr{loss}}J_a + \tilde{H}(a),
\]

\noindent where $K(a) = h_t \operatorname{diag}(K^1, \ldots, K^{n_t})$ is the block-diagonal Gram matrix of the control cost~\cref{eq:discrete_cost}, with blocks $K^j = (\operatorname{ker}_\sigma(z_k^j, z_l^j))_{k,l=1}^{n_s}$ evaluated on the current trajectory $z = z(a)$; $J_a = \partial_a z^{n_t+1}$ is the linear terminal sensitivity, $H_{\scr{loss}}$ is the per-sample terminal-loss Hessian in $z$ (with $D_i$ the softmax Hessian)
\[
H_{\scr{loss}}(z_i^{n_t+1}) = \Xi^\mathsf{T} D_i\Xi,
\qquad D_i = \operatorname{diag}(\hat{y}_i) - \hat{y}_i \otimes \hat{y}_i,
\qquad \hat{y}_i = c_{\Xi^\star\!(a)}(z_i^{n_t+1}),
\]

\noindent and $\tilde{H}(a)$ collects the second-order ODE sensitivity $\partial_{aa} z^{n_t+1}$ together with the indirect dependence of $K$ on $a$ through the state. The exact reduced Hessian $\partial^2_a\widehat{f}(a)$ also carries a variable-projection coupling term $-\,\partial^2_{a\bar\Xi}f\,(\partial^2_{\bar\Xi\bar\Xi}f)^{-1}\,\partial^2_{\bar\Xi a}f \preceq 0$ from the sensitivity $\partial_a\bar\Xi^\star(a)$ of the eliminated classifier (the barred derivatives are taken along $\mathcal{W}$; see \Cref{app:varpro-proof}). Dropping both this coupling term and $\tilde{H}(a)$ gives the GN approximation $H_{\scr{gn}}(a)$, which is positive definite whenever $K(a)\succ0$, i.e.\ for pairwise-distinct particles at every time level. The search direction $s_a \in \mathbb{R}^{dn_sn_t}$ is obtained by approximately solving $H_{\scr{gn}}(a)s_a = - g_a$ via a CG method with relative tolerance $\operatorname{tol}_{\scr{cg}} > 0$ (a fixed inexact-Newton forcing~\cite[\S7.1]{nocedal2006:numerical} with default $0.5$) and maximum iteration count $m_{\scr{cg}} \in \mathbb{N}$ (default $10$). A tolerance sweep confirmed that this loose forcing already attains the achievable optimum at minimal cost and that neither a tighter constant tolerance nor an adaptive Eisenstat--Walker forcing sequence~\cite[\S7.1]{nocedal2006:numerical} improves it.

We do not form or store the Hessian; we compute its action on a vector. Given a candidate $\tilde{a} = (\tilde{a}^j)_{j=1}^{n_t}$, the forward sensitivities are obtained via
\[
\tilde{z}^{j+1}_k = \tilde{z}^j_k + h_t
\sum_{l=1}^{n_s}
\Bigl[
\operatorname{ker}_{\sigma}(z_k^j, z_l^j) \tilde{a}^j_l
+ \bigl(
\langle \partial_1\!\operatorname{ker}_{\sigma}(z_k^j, z_l^j), \tilde{z}^j_k\rangle
+ \langle \partial_2\!\operatorname{ker}_{\sigma}(z_k^j, z_l^j), \tilde{z}^j_l\rangle
\bigr) a_l^j
\Bigr]
\]

\noindent with $\tilde{z}^1 = 0$ and terminal incremental state $\tilde{z}^{n_t+1} = J_a \tilde{a}$. To compute the action of the softmax Hessian we apply $H_{\scr{loss}}(z_i^{n_t+1})$ to $\tilde{z}_i^{n_t+1}$ for each sample $i$:
\[
w_i = \lambda_{\scr{loss}}\,\Xi^\mathsf{T} D_i\Xi\tilde{z}_i^{n_t+1},
\quad
D_i u = \hat{y}_i \odot u - \langle \hat{y}_i, u\rangle \hat{y}_i,
\]

\noindent where $\odot$ denotes the entrywise (Hadamard) product and $u = \Xi\,\tilde{z}_i^{n_t+1} \in \mathbb{R}^{n_c}$. To compute the adjoint sensitivities, we solve the incremental adjoint equation
\[
\tilde{p}_k^{j-1}
= \tilde{p}_k^j + h_t \sum_{l=1}^{n_s}
\bigl(
\langle \tilde{p}_k^j, a_l^{j-1}\rangle
+ \langle \tilde{p}_l^j, a_k^{j-1}\rangle
\bigr)
\partial_1\!\operatorname{ker}_{\sigma}(z_k^{j-1}, z_l^{j-1})
\]

\noindent with $\tilde{p}_k^{n_t+1} = -w_k$ backward in time.

The natural-coordinate Hessian matvec is $(H_{\scr{gn}}^{\scr{nat}}(a) \tilde{a})_k^j = 2\tilde{a}_k^j - \tilde{p}_k^{\,j+1}$, $j = 1,\ldots,n_t$. For pairwise-distinct particles at every time level, the operator $H_{\scr{gn}}^{\scr{nat}}(a)$ is symmetric positive definite in the natural inner product
\begin{equation}\label{eq:nat-ip}
\langle u, v\rangle_{\scr{nat}}
\defeq h_t \sum_{j=1}^{n_t} \langle u^j, K^j v^j\rangle.
\end{equation}

\noindent All metric-dependent quantities in our CG implementation (residual norms, curvature check, CG coefficients, relative stopping test) use~\cref{eq:nat-ip}. The kernel action $K^j$ is applied matrix-free at the frozen trajectory $\{z^j\}$ of the current outer iterate (the \texttt{PyKeOps} backend never materializes an $n_s\times n_s$ block); each CG iteration then costs one Hessian matvec plus $O(n_t n_s^2 d)$ for the natural inner product evaluations, comparable to one extra forward solve per CG iteration. \Cref{alg:varpro-gncg} summarizes the outer step. This natural-inner-product implementation is that of the full-expansion VarPro solver; the control-point solvers of \Cref{app:cp-kkt} instead run CG in the standard Euclidean inner product on the Euclidean GN operator, which is internally consistent and is the path used for the scalability, ablation, and density-estimation results.

For the CNF density-estimation the terminal objective is the NLL $-\log\pi_1$ rather than the softmax cross-entropy, so the softmax Hessian $H_{\scr{loss}}$ is replaced by the curvature of $-\log\pi_1$ at the terminal state. For the GMM reference of \Cref{sec:gmm-reference} the exact terminal Hessian is
\[
\partial_{zz}\bigl(-\log\pi_1(z)\bigr)
= \sum_{i=1}^{n_{\scr{mod}}}\gamma_i(z)\,\Sigma_i^{-1}
- \operatorname{Cov}_{\gamma(z)}\!\bigl(\Sigma_i^{-1}(z-\mu_i)\bigr),
\]

\noindent where $\operatorname{Cov}_{\gamma(z)}$ denotes the responsibility-weighted covariance of the component scores $\Sigma_i^{-1}(z-\mu_i)$. The covariance term can be sufficiently large in regions between mixture modes that the exact Hessian becomes indefinite. We therefore use the responsibility-weighted GN (Fisher) surrogate
\[
H_{\scr{sur}}(z) = \sum_{i=1}^{n_{\scr{mod}}}\gamma_i(z)\,\Sigma_i^{-1} \succeq 0,
\]

\noindent which retains only the within-component precisions; since $\gamma_i(z) \ge 0$ and $\Sigma_i^{-1}\succ 0$ it is symmetric positive semidefinite, so $H_{\scr{gn}}$ remains symmetric positive definite. The incremental terminal condition of the Hessian--vector product then reads $w_k = (\lambda_{\scr{loss}}/n_s)\,H_{\scr{sur}}(z_k^{n_t+1})\,\tilde{z}_k^{n_t+1}$, in place of the softmax-Hessian action above. The curvature of the running-divergence term, which has no definite sign in general, is omitted from $H_{\scr{gn}}$; the divergence term still enters the reduced gradient exactly through the adjoint of \Cref{app:kkt_cnf}. For the UGM this reduces to $H_{\scr{sur}}=\Sigma^{-1}$, which is the exact terminal Hessian.

\subsubsection*{GN Trust-Region Model}

Here, we replace the line-search globalization with a TR model~\cite{conn2000:trust}. Let
\[
m_k(s_a) \defeq \langle g_a^{\scr{nat},(k)}, s_a\rangle_{\scr{nat}} + \tfrac{1}{2} \langle s_a, H_{\scr{gn}}^{\scr{nat},(k)} s_a\rangle_{\scr{nat}}.
\]

\noindent At outer iteration $k$ we solve the subproblem
\begin{equation}\label{eq:tr-subproblem}
\minopt_{s_a}\;\;
m_k(s_a)
\quad\text{subject to}\quad
\|s_a\|_{\scr{nat}} \leq \Delta^{(k)}.
\end{equation}

\noindent The TR radius $\Delta^{(k)}$ is measured in the natural norm~\cref{eq:nat-ip}. We solve~\cref{eq:tr-subproblem} approximately by the Steihaug--Toint truncated-CG method~\cite{steihaug1983:conjugate,toint1981:towards}, run in that norm: CG iterates until it reaches the trust-region boundary, where it truncates the last step to the boundary. The step is accepted, and the radius contracted or expanded, by the standard ratio test $\rho^{(k)}$ of actual to predicted reduction and its radius-update rule~\cite[Chap.~4]{nocedal2006:numerical} (see \Cref{alg:varpro-trgn}), with thresholds $\eta_1 = \mathrm{1e}\mathord{-}2$, $\eta_2 = 0.75$ and shrink/expand factors $\gamma_1 = 0.5$, $\gamma_2 = 2.0$. The radius floor is $\Delta_{\scr{min}} = \mathrm{1e}\mathord{-}6$, and $(\Delta_{\scr{init}}, \Delta_{\scr{max}}) = (1, 100)$ are scaled by $\|g_a^{(0)}\|_{\scr{nat}}$ so they act as unitless multipliers of the initial natural-norm gradient.

In the ratio $\rho^{(k)}$ (and likewise in the Armijo line search of \Cref{alg:varpro-gncg} and the Anderson descent check), the trial term $\widehat{f}(a^{(k)}+s_a^{(k)})$ is evaluated at the frozen classifier $\Xi^\star(a^{(k)})$, i.e., as $F(a^{(k)}+s_a^{(k)},\Xi^\star(a^{(k)}))$ with $F$ as in \Cref{app:varpro-proof}, rather than by re-solving the inner problem for $\Xi$. For the trust-region ratio and the Armijo test this is conservative: the base point satisfies $F(a^{(k)},\Xi^\star(a^{(k)}))=\widehat{f}(a^{(k)})$ exactly while $F(\,\cdot\,,\Xi^\star(a^{(k)}))\ge\widehat{f}(\,\cdot\,)$ over-estimates the trial objective, so the computed decrease is a lower bound on the true reduced decrease. The Anderson descent check $\widehat{f}(a_{\scr{aa}}^{(k)}) \leq \widehat{f}(\chi(a^{(k)}))$ compares two frozen-$\Xi$ values at different points, neither of which is the base point at which the frozen classifier is exact, so it carries no such guarantee for $\widehat{f}$; it acts as a heuristic filter, and any convergence guarantee is inherited from the fallback step $\chi(a^{(k)})$, under the standard assumptions for gradient descent with Armijo backtracking, rather than from the check itself~\cite{zhang2020:globally}. The classifier is re-solved to $\Xi^\star(a^{(k+1)})$ at the start of the next outer iteration.

\subsubsection*{Proximal GN for the Scale-Sparsity Penalty}

With the multiscale scale-sparsity penalty~\cref{eq:scale-sparsity} the reduced objective is composite,
\[
\widetilde{f}(a) = \nu(a) + \psi(a),
\]

\noindent where $\nu$ collects the smooth terms (the reduced data misfit and the kinetic energy $\sum_\sigma\langle a_\sigma, K_\sigma a_\sigma\rangle$) and $\psi$ is the nonsmooth group-log term~\cref{eq:scale-sparsity}. The GN-CG and TR-GN steps above target a smooth objective and do not apply to $\psi$; we therefore replace the outer step by a proximal GN scheme~\cite{lee2014:proximal} that keeps the smooth curvature model but treats $\psi$ through its proximal operator. At outer iteration $k$ we form the same GN model of $\nu$ used above and minimize the composite quadratic
\begin{equation}\label{eq:prox-gn-model}
\textstyle
q_k(s) = \langle g_a^{(k)}, s\rangle + \tfrac12\langle s, H_{\scr{gn}}^{(k)} s\rangle + \psi(a^{(k)}+s),
\qquad g_a^{(k)} = \partial_a\nu(a^{(k)}),
\end{equation}

\noindent approximately for the increment $s$. Because $H_{\scr{gn}}^{(k)}$ is available only through its matvec, we solve~\cref{eq:prox-gn-model} with an accelerated proximal-gradient inner iteration~\cite{beck2009:fast}: each inner step takes a gradient step on the smooth part of $q_k$ with step length $1/L$, where $L$ is a power-iteration estimate of $\|H_{\scr{gn}}^{(k)}\|$ (with a $1.25$ safety factor), refreshed every fifth outer step, followed by the proximal map of $\psi/L$.

The penalty is separable over the per-particle, per-scale, per-time-step blocks $u=a_{l,\sigma}^{j}\in\mathbb{R}^d$, so its proximal map decouples into independent group problems. For a block with weight $w=\lambda_{\scr{sp}}^{\sigma}h_t/L$ the proximal map $\operatorname{prox}_{w\rho_\varepsilon}(u)$ with $\rho_\varepsilon(u)=\log(1+\varepsilon^{-1}\|u\|_2)$ acts only on the radius $r=\|u\|_2$: writing $s=\|\hat u\|_2$ for the shrunk radius, stationarity of $\tfrac12(s-r)^2 + w\log(1+s/\varepsilon)$ gives the scalar quadratic
\[
s^2 + (\varepsilon-r)s + (w-\varepsilon r) = 0,
\qquad
\Delta = (r+\varepsilon)^2 - 4w,
\]

\noindent whose larger root $s^\star=\tfrac12\big(r-\varepsilon+\sqrt{\Delta}\big)$ is kept when $\Delta\ge0$, $s^\star>0$, and it lowers the group objective below the value $\tfrac12 r^2$ at $s=0$; otherwise the block is set to exactly zero, so the map both shrinks and prunes. The shrunk block is $\hat u=(s^\star/r)u$. The resulting increment $s$ defines the outer search direction, globalized by a composite Armijo backtracking line search on $\widetilde{f}$ (sufficient decrease $\widetilde f(a+\gamma s)\le\widetilde f(a)+c_1\gamma\Delta$ with $\Delta=\langle g_a^{(k)},s\rangle+\psi(a+s)-\psi(a)<0$, as in~\cite{lee2014:proximal}); if the GN increment fails to produce sufficient decrease we fall back to a proximal-gradient step on $\widetilde{f}$. Since $L$ only estimates $\|H_{\scr{gn}}^{(k)}\|$, and not $\operatorname{lip}(\nabla\nu)$, monotone decrease is not guaranteed a priori, so this fallback is itself globalized by the composite backtracking line search and the step is rejected if no sufficient decrease is found. \Cref{alg:varpro-proxgn} summarizes the resulting outer step.

We measure stationarity of the composite problem by the max-norm proximal-gradient residual
\begin{equation}\label{eq:prox-residual}
\textstyle
r_\infty^{(k)} = \big\|a^{(k)} - \operatorname{prox}_{\psi}(a^{(k)} - g_a^{(k)})\big\|_\infty,
\end{equation}

\noindent which reduces to the reduced-gradient max-norm $\|g_a^{(k)}\|_\infty$ in the smooth (single-scale, $\psi\equiv0$) case solved by GN-CG. We terminate the outer iteration when $r_\infty^{(k)}<\mathrm{1e}\mathord{-}4$ or at a cap of 200 outer steps. The multiscale runs also enable an objective-stagnation stop, which triggers after ten consecutive outer steps with a relative objective improvement below $\mathrm{1e}\mathord{-}5$. It is disabled for the single-scale runs: the inner EM re-fit of the GMM reference makes the objective non-monotone across outer steps, so this test can fire before the budget is exhausted. The max-norm makes the tolerance dimension-free, so the same threshold applies across problem sizes; we report both the $\ell_2$ and $\ell_\infty$ residual norms in \Cref{app:genmod-details}.

\section{VAE Latent Embedding for the MNIST Experiments}
\label{app:vae}

Two of our MNIST experiments operate not on the raw $28\times28$ pixels but on a low-dimensional latent code produced by a convolutional VAE~\cite{kingma2014:auto}: the density-estimation example of \Cref{sec:res-genmod-mnist} and the classification example of \Cref{sec:res-mnist}. Embedding sidesteps the degradation of kernel methods in high dimensions noted in \Cref{sec:limitations}: the flow then operates in a $d_z$-dimensional latent space rather than in the $784$-dimensional pixel space. The same VAE architecture is used in both cases, differing only in the latent dimension $d_z$; it is trained separately and kept frozen, and both the CNF and the GN-CG classifier are fit to the standardized encoder means. For the generative task, latent codes are decoded back to images by the VAE decoder.

\paragraph{Architecture} The encoder maps an image $x\in[0,1]^{1 \times 28\times28}$ through four $3\times3$ convolutions of channel widths $(32,64,64,64)$ with ReLU activations (the second convolution has stride 2, halving the spatial resolution to $14\times14$, and the others stride $1$), followed by a fully connected layer to a hidden width $h$ and two linear heads that output the latent mean $\mu\in\mathbb{R}^{d_z}$ and log-variance $\log\sigma^2\in\mathbb{R}^{d_z}$. A latent code is drawn by the reparameterization $z=\mu+\sigma\odot\varepsilon$ with $\varepsilon\sim\mathcal{N}(0,I_{d_z})$. The decoder mirrors the encoder: a fully connected layer maps $z$ to a $64\times14\times14$ feature map, a stride-2 transposed convolution restores the $28\times28$ resolution, and a final $3\times3$ convolution with a sigmoid activation produces the reconstruction $\hat{x}\in[0,1]^{1\times28\times28}$.

\paragraph{Training and embedding}
We train the VAE with the Adam optimizer by minimizing the $\beta$-weighted VAE objective
\[
f_{\scr{vae}}
= \operatorname{dist}_{\scr{bce}}(\hat{x},x) + \beta\operatorname{reg}_{\scr{kl}}\big(\mathcal{N}(\mu,\operatorname{diag}\sigma^2)\,\big\|\,\mathcal{N}(0,I_{d_z})\big),
\]

\noindent the pixelwise binary cross-entropy reconstruction term plus a $\beta$-weighted Kullback--Leibler regularizer toward a standard-normal prior. The latent dimension is chosen per experiment: $d_z=4$ for the density-estimation example of \Cref{sec:res-genmod-mnist}, and $d_z\in\{16,32,64\}$ for the classification example of \Cref{sec:res-mnist}, where $d_z=32$ is the accuracy optimum; the hidden width is $h=32$, $\beta=\mathrm{5e}\mathord{-}4$, and the encoder is trained with Adam (learning rate $\mathrm{1e}\mathord{-}3$, batch size $128$). The reconstruction term is the binary cross-entropy averaged over the $784$ pixels, and the KL term is summed over the latent coordinates; both are averaged over the batch.

After training, the encoder means $\mu$ are extracted and standardized to zero mean and unit variance per coordinate before entering the solver. For density estimation the standardized means of the \inum{10000} training, \inum{10000} validation, and \inum{10000} test images feed the CNF, and new digits are generated by sampling the CNF reference, applying the learned inverse flow, undoing the standardization, and decoding the resulting latent with the VAE decoder. For classification the same standardized means feed the GN-CG classifier, with the training/test split reported there.

\section{Numerical Results}

Here, we report complementary results for those reported in the main text.

\subsection{Order of Convergence of the Discrete Adjoint Gradient}
\label{app:adjoint-order}

We verify that the discrete adjoints of the time integrators produce a reduced gradient whose accuracy matches the order of the underlying time integrator. We fix one classification instance (the three rings data set, with $n_s=100$ particles subsampled from the training set, $d=3$ and $\sigma\approx0.5$) and a time-constant control $\bar a\in\mathbb{R}^{n_s\times d}$ with i.i.d.\ $\mathcal{N}(0,0.3^2)$ entries. Holding the control constant in time makes it independent of the number of time steps $n_t$, so the same parameter can be compared across temporal resolutions. The gradient compared across resolutions is the gradient with respect to this shared parameter, $\nabla_{\bar a}F_{n_t}=\sum_{j=1}^{n_t}\partial_{a^j}F_{n_t}\in\mathbb{R}^{n_s\times d}$. We take an RK4 solution with $n_{\scr{ref}}=2048$ steps as the continuous reference; we then refine $n_t\in\{2,4,8,\dots,256\}$ and measure the Euclidean-norm error of three quantities against the reference: the terminal state $z^{n_t+1}$, the terminal-loss gradient, and the reduced gradient.

\begin{figure}
\centering
\includegraphics[width=\textwidth]{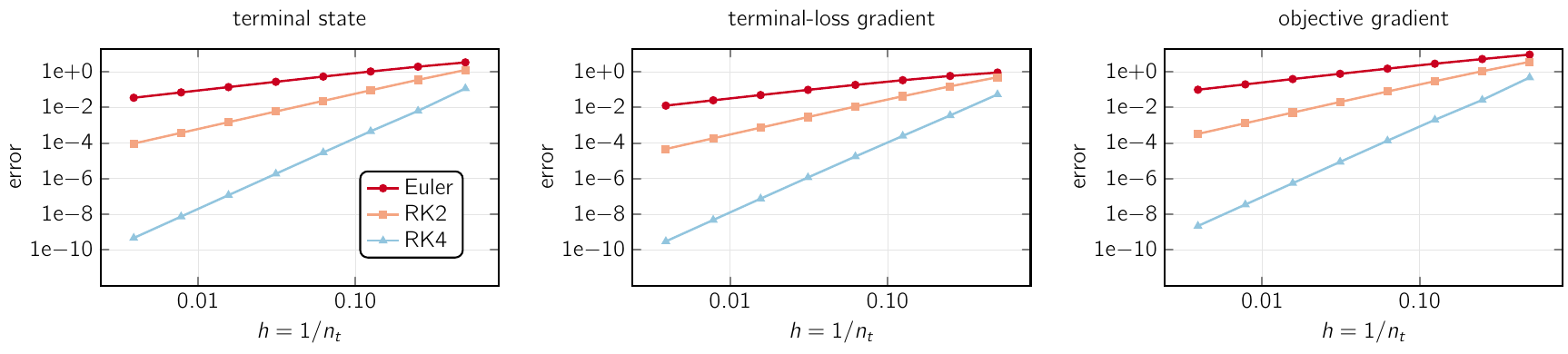}
\caption{Order-of-convergence study of the discrete adjoint gradient on the three rings classification instance with a time-constant control. We report the Euclidean error versus step size $h_t=1/n_t$ for the terminal state (left), the terminal-loss (flow-adjoint) gradient (center), and the reduced gradient (right), against an RK4 reference at $n_{\scr{ref}}=2048$.}
\label{fig:adjoint-order}
\end{figure}

\Cref{fig:adjoint-order} shows the three error curves versus $h_t$ on a log--log scale. All three quantities converge at the order of the integrator; the fitted orders are collected in \Cref{tab:adjoint-order}. The coarsest resolutions are pre-asymptotic (a fit over the full range $n_t=2$--$256$ is as low as 0.894 for the Euler terminal-loss gradient), so we fit every slope over the four finest steps ($n_t=32$--$256$), where all three quantities recover the expected orders ($\approx$1 for Euler, and correspondingly $\approx$2 and $\approx$4 for RK2 and RK4).

\begin{table}
\tabadjust
\caption{Empirical convergence orders (log--log least-squares slopes) of the terminal state and of the discrete adjoint gradient of the terminal loss and of the full objective. All nine slopes are fitted over the same asymptotic window, the four finest resolutions $n_t=32,\ldots,256$; each quantity converges at the order of the corresponding time integrator.}
\label{tab:adjoint-order}
\begin{tabular}{lcccc}
\toprule
 & & & \multicolumn{2}{c}{gradient}\\
integrator & theoretical & terminal state & terminal-loss & objective \\
\midrule
Euler & 1 & \fnum{0.992457400} & \fnum{0.979324469} & \fnum{0.989203769} \\
RK2   & 2 & \fnum{1.993996660} & \fnum{1.987539691} & \fnum{1.985852873} \\
RK4   & 4 & \fnum{3.992949442} & \fnum{3.973871628} & \fnum{3.986546118} \\
\bottomrule
\end{tabular}
\end{table}

The study above compares the discrete gradient against a fine continuous reference; on its own it does not establish that the discrete adjoint is the \emph{exact} gradient of the \emph{discrete} objective. We verify that separately by a Taylor-remainder test. For a random unit direction $\tilde{a}$ we evaluate
\begin{equation}\label{eq:taylor-remainder}
t_0(h) = \big|F(a + h\tilde{a},\Xi) - F(a,\Xi)\big|,
\qquad
t_1(h) = \big|F(a + h\tilde{a},\Xi) - F(a,\Xi) - h\,\langle g_a, \tilde{a}\rangle\big|,
\end{equation}

\noindent at a generic (nonzero, random) control $a$ and frozen $\Xi$, with $g_a$ the discrete adjoint gradient~\cref{eq:reduced-gradient}. Here $t_0 = O(h)$ holds for any $g_a$, whereas $t_1=O(h^2)$ holds if and only if $\langle\nabla_aF(a,\Xi)-g_a,\tilde a\rangle=0$: a wrong gradient leaves a first-order remainder unless its error is orthogonal to $\tilde a$. On a random instance ($n_s=20$ standard-normal particles in $d=4$, $n_c=5$ classes with random one-hot labels, $n_t=20$ Euler steps, $\sigma=3$, $\lambda_{\scr{cls}}=\mathrm{1e}\mathord{-}3$, $\lambda_{\scr{loss}}=6$), in double precision, the fitted slopes are $t_0\sim h^{1.01}$ and $t_1\sim h^{2.00}$. The ratio $t_1/h^2$ is constant ($\approx\mathrm{5.098e}\mathord{-}2$) until $t_1$ reaches the float64 cancellation floor $\approx\epsilon_{\scr{mach}}|F(a,\Xi)|=O(\mathrm{1e}\mathord{-}14)$ near $h=\mathrm{1e}\mathord{-}6$. The same test on the classifier block gives $t_1\sim h^{2.00}$ with ratio $\mathrm{2.214e}\mathord{-}1$ over the same range. The Taylor test confirms the discrete gradient to numerical precision along the tested directions on the tested instances.

\subsection{Inverse Consistency}\label{app:inv-consistency}

\begin{figure}
\centering
\includegraphics[width=0.6\textwidth]{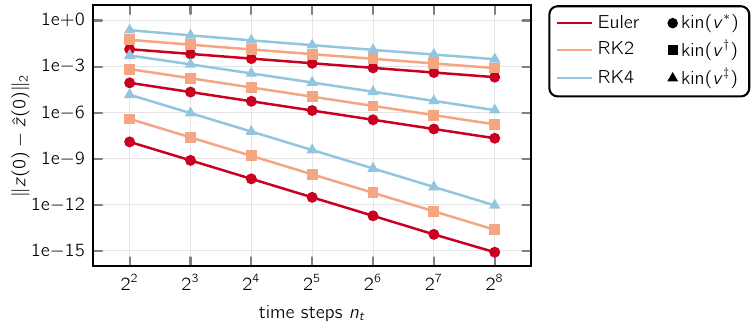}
\caption{Inverse consistency of the discrete flow for a synthetic (autonomous) RKHS velocity field: forward--backward round-trip error $\|z(0)-\hat z(0)\|$ versus the number of forward-Euler time steps $n_t$, for three kinetic-energy levels $\operatorname{kin}(v\ast) < \operatorname{kin}(v^\dag) < \operatorname{kin}(v^\ddag)$.}
\label{fig:inverse-consistency}
\end{figure}

As a consistency check on the discrete flow, we integrate a synthetic (autonomous) RKHS velocity field (nothing is learned here) forward from $t=0$ to $t=1$ and then backward to $t=0$, and measure the inverse consistency error $\|z(0)-\hat z(0)\|_2$ as a function of the number of time steps $n_t$ and the flow's kinetic energy \[\operatorname{kin}(v) = \int_0^1 \|v_t\|_{\mathcal{H}}^2\,\mathrm{d}t.\] We use an autonomous RKHS field ($n_s = 64$ particles in $d = 2$) at three energy levels obtained by scaling the control field. \Cref{fig:inverse-consistency} shows that the inverse consistency error decays at first order in the step size (fitted orders $1.01$--$1.04$), consistent with the explicit forward/backward numerical schemes, and grows monotonically with the kinetic energy.

\subsection{Divergence Evaluation}\label{app:res-diveval}

\begin{figure}
\centering
\includegraphics[width=\textwidth]{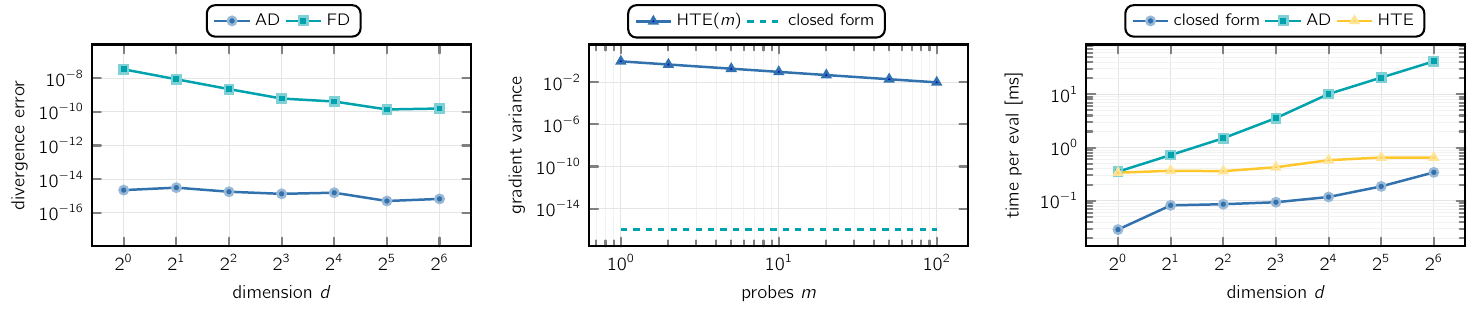}
\caption{Exact divergence for RKHS velocity fields. Left: Agreement of the closed-form divergence~\cref{eq:exact-div-gauss} with an AD trace and an FD trace vs dimension $d$. Middle: HTE / gradient variance versus the number of probes $m$; the exact expression has zero variance. Right: Wall-clock per divergence evaluation versus $d$ for the closed-form, FD-trace, and HTE variants.}
\label{fig:exact-div}
\end{figure}

\Cref{fig:exact-div} verifies the claim that the divergence of the RKHS velocity field is available in exact, deterministic, closed form. We evaluate the closed-form expression~\cref{eq:exact-div-gauss} on random particle configurations $\{z_k\}_{k=1}^{n_s}$ with $n_s = 64$ and controls drawn i.i.d.\ from a standard normal, sweeping the ambient dimension $d \in \{1,2,\ldots,64\}$ and the number of Hutchinson probes $m$. The kernel bandwidth is set per dimension by the median-of-pairwise-distances heuristic so that the field remains non-degenerate as $d$ grows. We report three quantities.

We compare the closed-form divergence against two independent references: a forward-mode automatic-differentiation (AD) trace of the velocity Jacobian and a central finite-difference (FD) trace, both built from the kernel matrix alone and thus independent of the analytic gradient used in~\cref{eq:exact-div-gauss}. The closed form matches the AD trace to machine precision (the maximum absolute discrepancy remains below $\mathrm{4e}\mathord{-}15$ in float64 across all tested dimensions), confirming that~\cref{eq:exact-div-gauss} is exact rather than a low-order approximation. The FD trace, by contrast, is limited by truncation and agrees only to $\mathcal{O}(\mathrm{1e}\mathord{-}8)$.

The standard deviation of the Hutchinson trace estimator (HTE) decays as $m^{-1/2}$ and the induced training-gradient variance as $m^{-1}$, the rates expected of a stochastic trace estimator; the exact expression has identically zero variance and requires no probes.

Consistent with the $\mathcal{O}(n_s^2 d)$ bound of~\Cref{lem:exact-div} (a single kernel matrix--vector product), the closed-form evaluation is the least expensive of the three variants: its wall-clock cost grows only mildly with $d$ (a factor of about $4$ across the $32$-fold increase from $d=2$ to $d=64$); at the small particle count $n_s=64$ used here the kernel reduction is dominated by launch overhead rather than the $d$-linear arithmetic, so the measured cost is sub-linear in $d$. The full AD trace instead scales linearly in $d$ and is $121\times$ slower at $d=64$. Even a single-probe HTE is more expensive than the exact evaluation while incurring the variance quantified in the plot in the middle. In the RKHS setting the stochastic estimator therefore offers no advantage over the exact expression, in either accuracy or cost.

\subsection{Time-Integrator Order}\label{app:integration-order}

We explore if the time-integration fidelity carries over to the downstream task: does a higher-order integrator produce a better \emph{classifier}.

We consider the interlocking tori in $\mathbb{R}^3$ (lifted to $\mathbb{R}^4$), which a linear classifier separates at only 66.7\% test accuracy. For each explicit integrator and each number of time steps $n_t$, we solve the VarPro problem to convergence. We deploy the learned control and classifier through a fine reference flow (RK4 at $n_t=64$) to measure test accuracy. We consider LBFGS as an outer optimizer. We report the hardware-independent work $n_t n_{\scr{stages}} (2 + n_{\scr{fev}} + 2n_{\scr{gev}})$, where $n_{\scr{stages}}\in\{1,2,4\}$ is the number of RK stages, the leading 2 counts the reduced gradient (one forward and one adjoint solve) evaluated per outer step (consistent with~\cref{eq:varpro-work}), and $n_{\scr{fev}}$, $n_{\scr{gev}}$ count the additional forward and gradient evaluations of the line search; the inner classifier solve is integrator-independent, and therefore excluded.

We report results in \Cref{tab:time-integrator-order-vs-classification-accuracy}: the test accuracy is $99.17\%$ (119/120) for every integrator and every choice of $n_t\in\{2,\dots,32\}$, against the $66.7\%$ baseline.

The cost grows by roughly the stage factor $s$: at $n_t=32$, RK4 spends $\sim 4.3\times$ the work of Euler for identical accuracy. For the considered example we conclude that higher-order integration buys no classification accuracy and only costs more.

\begin{table}
\tabadjust
\caption{Test accuracy as a function of the number of time steps used during training and the order of time integration. We consider the two interlocking tori as a test problem. L-BFGS is the outer solver. We report the test accuracy and the work for $n_t \in \{2,32\}$ for Euler, RK2, and RK4.}
\label{tab:time-integrator-order-vs-classification-accuracy}
\begin{tabular}{lrrrr}
\toprule
& $n_t$ & \bf Euler & \bf RK2 & \bf RK4 \\
\midrule
test accuracy &  2 & \fnum{0.991666667} & \fnum{0.991666667} & \fnum{0.991666667} \\
test accuracy & 32 & \fnum{0.991666667} & \fnum{0.991666667} & \fnum{0.991666667} \\
\midrule
work          & 32 & \inum{4576} & \inum{9024} & \inum{19456} \\
\bottomrule
\end{tabular}
\end{table}

\section{Additional Classification Results}\label{app:additional-class-results}

We collect additional material for the classification experiments of \Cref{sec:results}: the problem and solver settings, the control-point ablation, flow visualizations for the toy problems, and the convergence data underlying the main-text figures.

\subsection{Problem and Solver Settings}\label{app:class-settings}

\Cref{tab:class-settings} lists the settings for the five synthetic benchmarks. All five share the solver configuration: $n_t=20$, $\lambda_{\scr{loss}}=1$, $\lambda_{\scr{cls}}=\mathrm{1e}\mathord{-}2$, the cold start $a\equiv0$, and the Newton inner $\Xi$-solve of \Cref{sec:inner-solver}. At this cold start the flow is the identity, so the inner solve already returns the optimal linear-softmax classifier; the flow then improves on that baseline instead of on a degenerate $\Xi=0$ classifier. Each planar problem is standardized and lifted by one auxiliary affine coordinate, so the solver dimension is $d+1$. The bandwidth $\sigma$ is set per problem by the median rule of \Cref{sec:bandwidth} on the lifted training features; the resulting values are reported below. Data instances are regenerated per seed with a fixed base seed; the five reported seeds differ only in the data instance.

\begin{table}
\tabadjust
\caption{Problem and solver settings for the five synthetic classification benchmarks. From left to right: training and test sample counts, solver dimension $d$ (after the affine lift), number of classes $n_c$, and the median-rule bandwidth $\sigma$. All runs share $n_t=20$, $\lambda_{\scr{loss}}=1$, and $\lambda_{\scr{cls}}=\mathrm{1e}\mathord{-}2$.}
\label{tab:class-settings}
\begin{tabular}{lccccc}
\toprule
problem & $n_s$ (train) & $n_s$ (test) & $d$ & $n_c$ & $\sigma$ \\
\midrule
two moons     & \inum{300} & \inum{200} & 3 & 2 & \fnum{0.5149} \\
noisy circles & \inum{300} & \inum{150} & 3 & 2 & \fnum{0.5187} \\
three rings   & \inum{300} & \inum{240} & 3 & 3 & \fnum{0.5170} \\
3D tori       & \inum{200} & \inum{160} & 4 & 2 & \fnum{0.7106} \\
XOR           & \inum{300} & \inum{200} & 3 & 2 & \fnum{0.5838} \\
\bottomrule
\end{tabular}
\end{table}

\subsection{Control-Point Ablation}

\Cref{tab:cp-ablation} reports the full control-point sweep summarized in \Cref{sec:res-scalability} of the main text.

\begin{table}
\tabadjust
\caption{The control-point reduction is near-lossless. Test accuracy (\%) of the reduced GN-CG solver on the five synthetic benchmarks as the control-point count $m$ is varied at fixed data size $n_s$; $m=n_s$ recovers the full model~\cref{eq:rkhs-flow}. A small fraction of control points reaches the full-model accuracy on four of the five problems and comes within 0.8 percentage points on three rings.}
\label{tab:cp-ablation}
\begin{tabular}{lrrrrrr}
\toprule
\bf problem & $n_s$ & $m{=}8$ & $m{=}16$ & $m{=}32$ & $m{=}64$ & $m{=}n_s$ \\
\midrule
two moons     & 300 & 100.0 & 100.0 & 100.0 & 100.0 & 100.0 \\
noisy circles & 300 & 100.0 & 100.0 & 100.0 & 100.0 & 100.0 \\
three rings   & 300 &  77.1 &  99.2 &  99.2 &  99.2 & 100.0 \\
3D tori       & 200 &  96.9 &  98.8 & 100.0 & 100.0 & 100.0 \\
XOR           & 300 & 100.0 & 100.0 & 100.0 & 100.0 & 100.0 \\
\bottomrule
\end{tabular}
\end{table}

\Cref{tab:cp-mnist} reports the corresponding real-data comparison on the VAE-latent MNIST classification of \Cref{sec:res-mnist}, summarized in \Cref{sec:res-scalability} of the main text.

\begin{table}
\tabadjust
\caption{Control-point reduction on real data: VAE-latent MNIST classification (\inum{20000} train / \inum{5000} test, identical objective and encoder as the full model). For each latent dimension $d_z$ we report the full-model test accuracy and the control-point test accuracy at the control count $m$ that attains it ($m\le 2048$, i.e., $\approx 10\%$ of $n_s$). The control-point solver is near-lossless throughout, within 0.56 percentage points of the full model, approximately matching it at $d_z=16$ and exceeding it at $d_z=64$.}
\label{tab:cp-mnist}
\begin{tabular}{ccc}
\toprule
\bf latent dimension $d_z$ ($d$) & \bf full model & \bf reduced model ($m$) \\
\midrule
16 (17) & 95.86 & 95.82 ($\phantom{0}256$)   \\
32 (33) & 96.80 & 96.24 ($2048$)  \\
64 (65) & 94.96 & 95.02 ($1024$)  \\
\bottomrule
\end{tabular}
\end{table}

\subsection{Full-Model Computational Scalability}
\label{app:full-scalability}

\begin{figure}
\centering
\includegraphics[width=\textwidth]{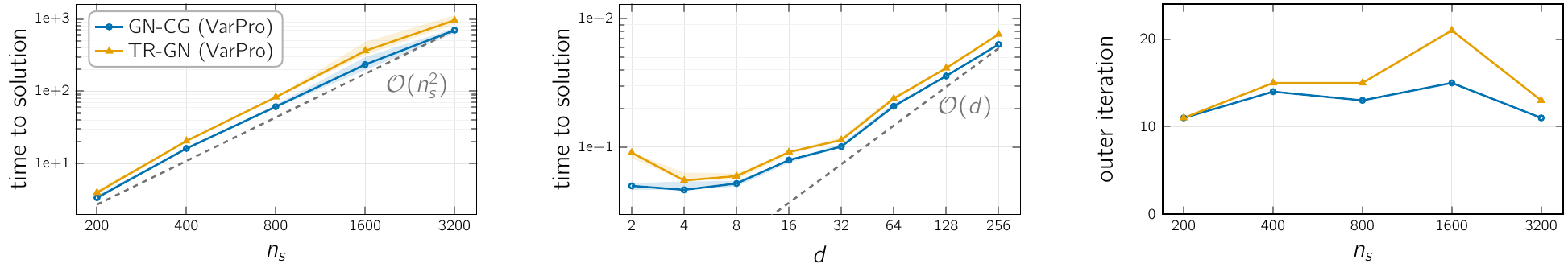}
\caption{Computational scalability of the reduced VarPro solve (GN-CG, TR-GN; median over three data-instance seeds with $25/75$ band). Left: wall time against the number of particles $n_s$ (3D tori, log--log) with the dashed $\mathcal{O}(n_s^2)$ reference. Middle: wall time against the ambient dimension $d$ (noisy circles isometrically embedded into $\mathbb{R}^d$, log--log) with the dashed $\mathcal{O}(d)$ reference. Right: outer-iteration count against $n_s$ (linear axis).}
\label{fig:scalability}
\end{figure}

\Cref{fig:scalability} reports the full-model scalability study summarized in \Cref{sec:results} of the main text. We report the wall time of the reduced VarPro solve against the particle count $n_s$ and the ambient dimension $d$, together with the (bounded) outer-iteration count.

\subsection{Flow Visualizations}

\begin{figure}
\centering
\includegraphics[width=\textwidth]{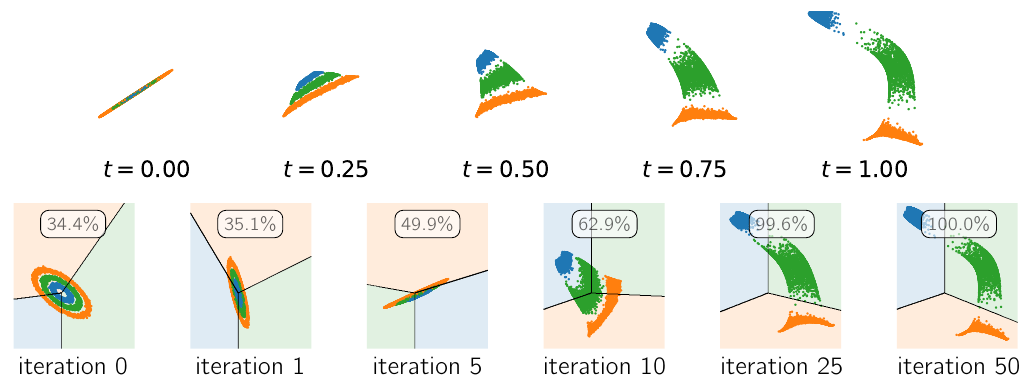}
\caption{Evolution of the three rings problem under the learned diffeomorphic model. Top: the original configuration and the transformed training data at selected flow times using the final controls $a$. Bottom: the transformed data at selected optimization iterations, with the linear classifier $\Xi$ shown in black and the predicted classes indicated by the background shading.}
\label{fig:3rings_flow_optimization_full}
\end{figure}

\begin{figure}
\centering
\includegraphics[width=\textwidth]{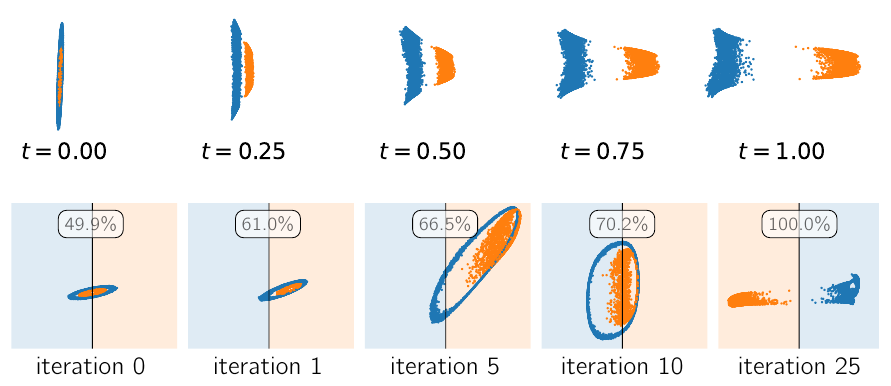}
\caption{Evolution of the noisy circles problem under the learned diffeomorphic model. Top: the original configuration and the transformed training data at selected flow times using the final controls $a$. Bottom: the transformed data at selected optimization iterations, with the linear classifier $\Xi$ shown in black and the predicted classes indicated by the background shading.}
\label{fig:noisy_circles_flow_optimization}
\end{figure}

\begin{figure}
\centering
\includegraphics[width=\textwidth]{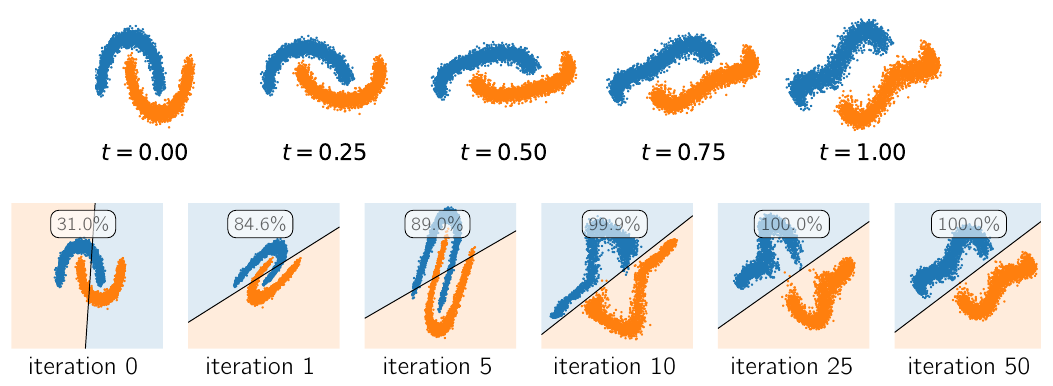}
\caption{Evolution of the two moons problem under the learned diffeomorphic model. Top: the original configuration and the transformed training data at selected flow times using the final controls $a$. Bottom: the transformed data at selected optimization iterations, with the linear classifier $\Xi$ shown in black and the predicted classes indicated by the background shading.}
\label{fig:moons_flow_optimization}
\end{figure}

\begin{figure}
\centering
\includegraphics[width=\textwidth]{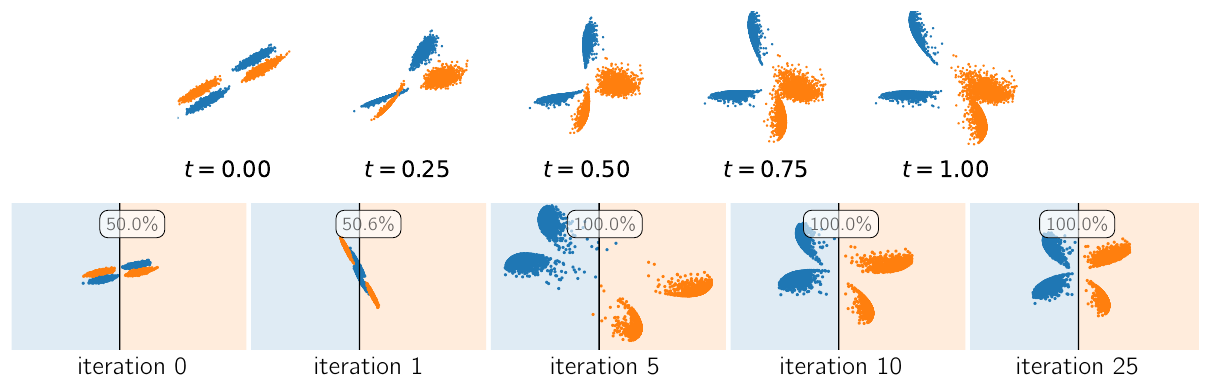}
\caption{Evolution of the XOR problem under the learned diffeomorphic model. Top: the original configuration and the transformed training data at selected flow times using the final controls $a$. Bottom: the transformed data at selected optimization iterations, with the linear classifier $\Xi$ shown in black and the predicted classes indicated by the background shading.}
\label{fig:xor_flow_optimization}
\end{figure}

\begin{figure}
\centering
\includegraphics[width=\textwidth]{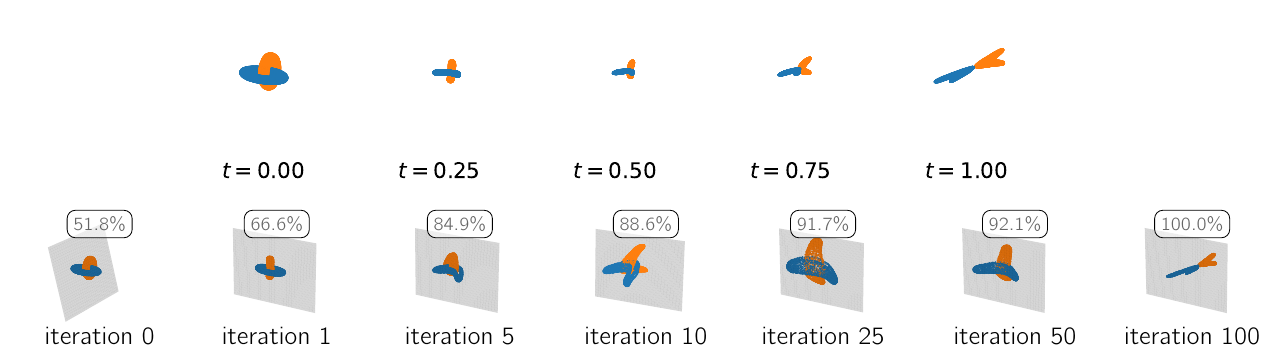}
\caption{Visualization of the learned diffeomorphic flow for the 3D tori dataset. The panels show the original configuration and the transformed training data at selected flow times under the final controls $a$. Points are colored by their true class labels.}
\label{fig:tori_flow}
\end{figure}

\Cref{fig:3rings_flow_optimization_full,fig:noisy_circles_flow_optimization,fig:moons_flow_optimization,fig:xor_flow_optimization,fig:tori_flow} visualize the learned diffeomorphic flow on the five toy problems. We show the transformed training data at selected flow times and optimization iterations, with the linear classifier and the predicted-class background shading. \Cref{fig:3rings_flow_optimization_full} is the three rings illustration referenced in \Cref{sec:results} of the main text.

\subsection{Convergence Results}

The following tables give the values underlying the classification convergence figures of the main text. In particular, \Cref{tab:varpro-full-baseline} for the full baseline comparison, \Cref{tab:varpro-variants} for the outer-solver variants, and \Cref{tab:varpro-baseline} for the four additional toy problems.

\begin{table}
\caption{Convergence study. We compare the top VarPro schemes versus joint $(a,\Xi)$ optimizers on the three rings problem. We report the median over five data-instance seeds. We report the reduced objective $\widehat{f}$ at the 500-iteration budget and the work and outer iterations to reach the common reference level $\widehat{f} \le 20.25$, with the number of seeds (out of five) reaching it. The work and iteration columns are medians over the seeds that reach the reference level only. Entries marked ``---'' do not reach this level within the budget; their reported $\widehat{f}$ is the value at iteration 500. The results correspond to those reported in \Cref{fig:varpro-full-baseline}.}
\label{tab:varpro-full-baseline}
\tabadjust
\begin{tabular}{lrrrc}
\toprule
\bf method     & \bf objective & \bf work & \bf iters & \bf reached\\
\midrule
GN-CG (VarPro)  & \fnum{19.406} &   \inum{2900} & 17 & 5/5 \\
TR-GN (VarPro)  & \fnum{19.385} &   \inum{3360} & 21 & 5/5 \\
L-BFGS (VarPro) & \fnum{19.353} &   \inum{4700} & 28 & 5/5 \\
NPRCG (VarPro)  & \fnum{19.386} &   \inum{4900} & 19 & 5/5 \\
\midrule
L-BFGS          & \fnum{19.421} &   \inum{2760} &  58 & 4/5 \\
NPRCG           & \fnum{19.924} &  \inum{54510} & 318 & 4/5 \\
AANGD           & \fnum{26.306} &  \inum{68940} & 316 & 2/5 \\
NANGD           & \fnum{20.239} &  \inum{98880} & 462 & 3/5 \\
NFRCG           & \fnum{20.192} & \inum{114580} & 410 & 3/5 \\
AAGD            & \fnum{20.332} & --- & --- & 0/5 \\
NGD             & \fnum{20.646} & --- & --- & 0/5 \\
PRCG            & \fnum{20.689} & --- & --- & 0/5 \\
NAGD            & \fnum{20.696} & --- & --- & 0/5 \\
GD              & \fnum{20.740} & --- & --- & 0/5 \\
FRCG            & \fnum{23.269} & --- & --- & 0/5 \\
\bottomrule
\end{tabular}
\end{table}

\begin{table}
\tabadjust
\caption{Comparison of VarPro outer-step variants on the three rings problem. We report the median over five data-instance seeds. We consider several optimization algorithms (rows). We report the reduced objective $\widehat{f}$, the Euclidean reduced-gradient norm $\|\partial_a\widehat{f}\|$, and its reduction relative to the initial norm (all at the 500-iteration limit), together with the work and the number of outer iterations to reach the common reference level $\widehat{f}\le20.25$. The results correspond to those of \Cref{fig:varpro-variants}.}
\label{tab:varpro-variants}
\begin{tabular}{lccccc}
\toprule
\bf variant & $\widehat{f}$ & $\|\partial_a\widehat{f}\|_2$ & $\|\partial_a\widehat{f}\|_{2,\scr{rel}}$ & \bf work & \bf iters\\
\midrule
NGD (VarPro)    & \fnum{19.469} & \snum{1.728e-01} & \snum{6.80e-02} & \inum{8360} & 51 \\
AANGD (VarPro)  & \fnum{19.467} & \snum{1.543e-01} & \snum{4.79e-02} & \inum{6260} & 36 \\
NPRCG (VarPro)  & \fnum{19.386} & \snum{2.931e-02} & \snum{1.15e-02} & \inum{4900} & 19 \\
L-BFGS (VarPro) & \fnum{19.353} & \snum{2.329e-02} & \snum{9.60e-03} & \inum{4700} & 28 \\
GN-CG (VarPro)  & \fnum{19.406} & \snum{1.435e-02} & \snum{5.64e-03} & \inum{2900} & 17 \\
TR-GN (VarPro)  & \fnum{19.385} & \snum{2.979e-02} & \snum{1.17e-02} & \inum{3360} & 21 \\
\bottomrule
\end{tabular}
\end{table}

\begin{figure}
\centering
\includegraphics[width=\textwidth]{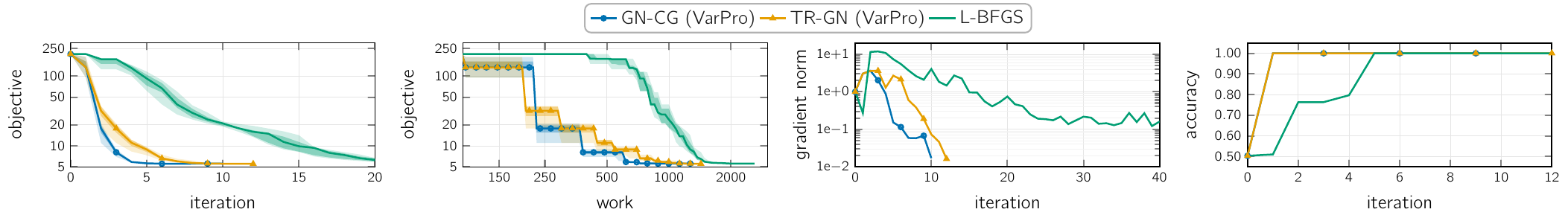}
\includegraphics[width=\textwidth]{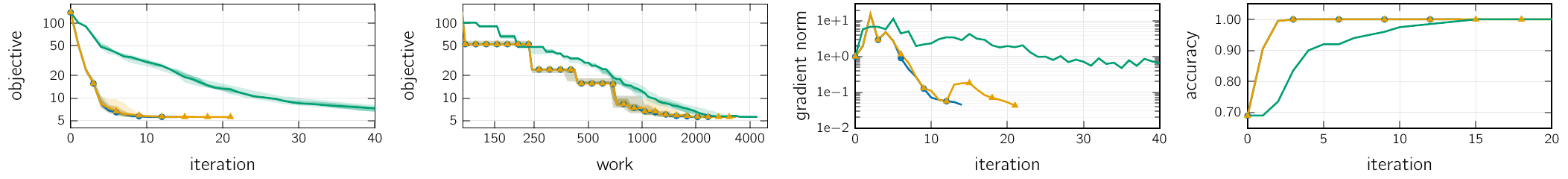}
\includegraphics[width=\textwidth]{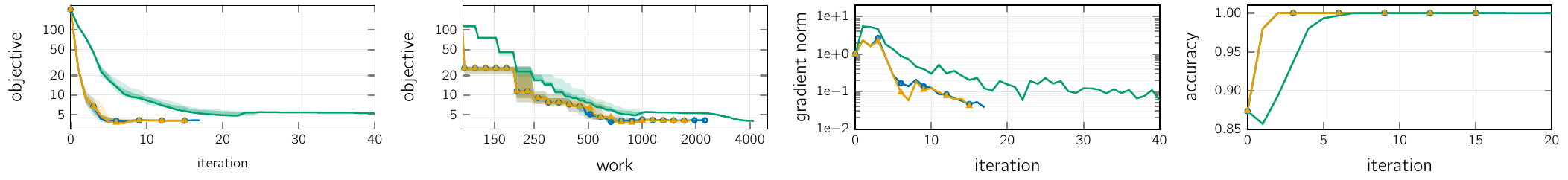}
\includegraphics[width=\textwidth]{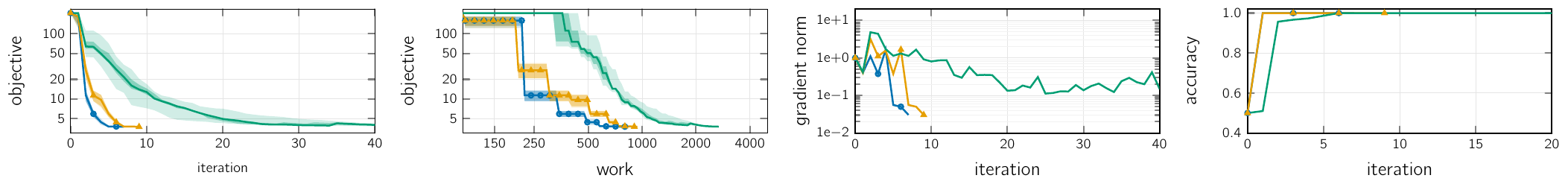}
\caption{Convergence summary for four additional toy problems. We consider the noisy circles (top row), the 3D tori (second row), the two moons (third row), and XOR (bottom row). We compare the VarPro solvers TR-GN, GN-CG versus joint $(a,\Xi)$ L-BFGS optimizer. We report the median over five data-instance seeds with $25/75$ inter-quartile band (faint color) and the min-max band (fainter color). Left: objective against outer-iteration count; middle left: objective against cumulative work (logarithmic axis); middle right: relative Euclidean gradient reduction against outer-iteration count; right: test accuracy against outer iteration count. The reported results correspond to \Cref{tab:varpro-baseline}.}
\label{fig:varpro-baseline}
\end{figure}

\begin{table}
\caption{Convergence summary for four additional toy problems. We report the median over five data-instance seeds. We report the final reduced objective $\widehat{f}$, the work and outer iterations to reach the common reference level (within $5\%$ of the best value attained on each instance, $\widehat{f}\le 5.76$, $5.74$, $3.83$, and $3.87$, respectively), the number of seeds (of five) reaching it, and the median training and test classification accuracy. As in \Cref{tab:varpro-full-baseline}, the work and iteration medians are taken over the seeds that reach the reference level only. Entries marked ``---'' do not reach this level under the common relative-gradient stopping criterion. The reported results correspond to \Cref{fig:varpro-baseline}.}
\label{tab:varpro-baseline}
\tabadjust
\begin{tabular}{llcrrccr}
\toprule
& & & & & & \multicolumn{2}{c}{\bf accuracy (\%)}\\
\cmidrule(lr){7-8}
\bf problem & \bf method & $\widehat{f}$ & \bf work & \bf iters & \bf reached & \bf train & \bf test\\
\midrule
noisy circles & GN-CG (VarPro) & \fnum{5.514} &  \inum{700} &   5 & 5/5 & 100.0 & 100.0 \\
              & TR-GN (VarPro) & \fnum{5.516} & \inum{1000} &   8 & 5/5 & 100.0 & 100.0 \\
              & L-BFGS         & \fnum{5.570} & \inum{1560} &  23 & 5/5 & 100.0 & 100.0 \\
\midrule
3D tori       & GN-CG (VarPro) & \fnum{5.618} & \inum{1410} &   8 & 4/5 & 100.0 & 100.0 \\
              & TR-GN (VarPro) & \fnum{5.618} & \inum{1470} &   9 & 4/5 & 100.0 & 100.0 \\
              & L-BFGS         & \fnum{5.626} & \inum{2360} &  55 & 4/5 & 100.0 & 100.0 \\
\midrule
two moons     & GN-CG (VarPro) & \fnum{3.695} &  \inum{720} &   6 & 3/5 & 100.0 & 100.0 \\
              & TR-GN (VarPro) & \fnum{3.692} &  \inum{720} &   6 & 3/5 & 100.0 & 100.0 \\
              & L-BFGS         & \fnum{4.775} &         --- & --- & 0/5 & 100.0 &  99.5 \\
\midrule
XOR           & GN-CG (VarPro) & \fnum{3.739} &  \inum{580} &   5 & 4/5 & 100.0 & 100.0 \\
              & TR-GN (VarPro) & \fnum{3.739} &  \inum{740} &   7 & 4/5 & 100.0 & 100.0 \\
              & L-BFGS         & \fnum{3.815} & \inum{1840} &  34 & 4/5 & 100.0 & 100.0 \\
\bottomrule
\end{tabular}
\end{table}

\subsection{Pendigits Scaling Study}
\label{app:pendigits-scaling}

\begin{figure}
\centering
\includegraphics[width=\textwidth]{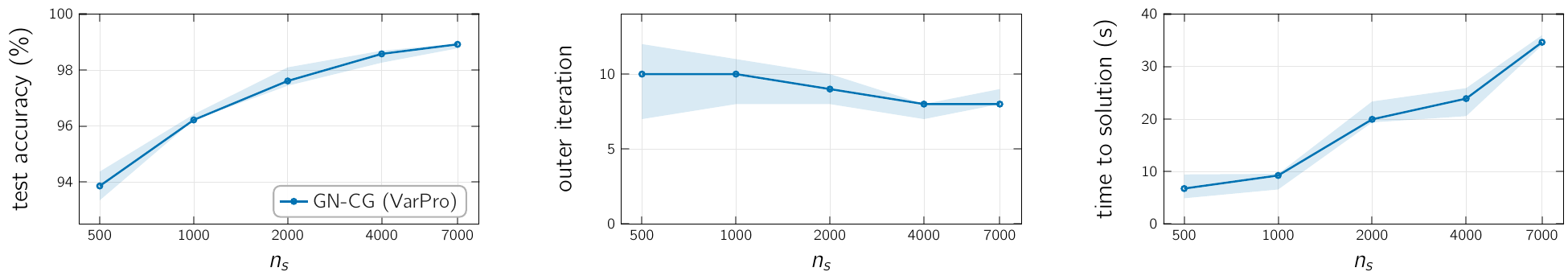}
\caption{Scaling of the reduced GN-CG solver on the real Pendigits data (median over three data-instance seeds with the min--max band), sweeping the training size $n_s$ at a fixed test set. Left: test accuracy against $n_s$. Middle: outer-iteration count against $n_s$. Right: wall-clock time to solution against $n_s$. The iteration count stays bounded while the test accuracy rises with $n_s$; the wall time stays under 35 seconds and grows sub-quadratically because the GPU is not saturated at these sizes.}
\label{fig:pendigits-scaling}
\end{figure}

We repeat the synthetic scalability study of \Cref{fig:scalability} on the real Pendigits data of \Cref{sec:res-pendigits}: holding the test set fixed at $n_{\scr{test}}=\inum{3518}$ samples and increasing the training size $n_s$ from 500 to \inum{7000}, we run the reduced GN-CG solver over three data-instance seeds and record the median (with the min--max band) of the test accuracy, the outer-iteration count, and the wall-clock time to solution. \Cref{fig:pendigits-scaling} and \Cref{tab:pendigits-scaling} report the results. Two trends carry over from the synthetic case. First, the outer-iteration count stays bounded (8--10) across the entire range and every run reaches 100\% training accuracy, so the optimizer effort does not grow with problem size and the wall time stays under 35 seconds; unlike the synthetic study, the wall time grows sub-quadratically here because at these sizes the GPU is not saturated. Second, generalization improves monotonically with $n_s$: the median test accuracy rises from 93.9\% at $n_s=500$ to 98.9\% at $n_s=\inum{7000}$. On this task the reduced solver turns additional data directly into accuracy at nearly fixed optimizer cost.

\begin{table}
\tabadjust
\caption{Pendigits scaling study. We report the median over three data-instance seeds at a fixed test set with $n_{\scr{test}}=\inum{3518}$. We increase the training size $n_s$. We report the test and training accuracy in \%, the outer-iteration count, and the wall-clock time to solution of the reduced GN-CG solver on a single GPU.}
\label{tab:pendigits-scaling}
\begin{tabular}{rrrrr}
\toprule
& \multicolumn{2}{c}{\bf accuracy (\%)}\\
$n_s$ & \bf test  & \bf train & \bf iterations & \bf time (s) \\
\midrule
        500 & 93.9 & 100.0 & 10 &  6.8 \\
\inum{1000} & 96.2 & 100.0 & 10 &  9.3 \\
\inum{2000} & 97.6 & 100.0 &  9 & 20.0 \\
\inum{4000} & 98.6 & 100.0 &  8 & 23.9 \\
\inum{7000} & 98.9 & 100.0 &  8 & 34.6 \\
\bottomrule
\end{tabular}
\end{table}

\section{Generative Modeling: Per-Problem Details}
\label{app:genmod-details}

This appendix reports the complete generative-modeling results summarized in \Cref{tab:genmod-summary} and discussed in \Cref{sec:syn-data-dist}. All fits use the control-point solver with $m=256$ control points. For each target (concentric circles, M8G, pinwheel, checkerboard) and each reference model (UGM; GMM eliminated by VarPro) we give the single-scale and multiscale metrics with gradient-norm diagnostics, together with the corresponding transport visualizations. The reconstruction error $e_2$ reported throughout is the forward--backward round trip $x_i\mapsto z_i=\phi(x_i)\mapsto\hat x_i=\phi^{-1}(z_i)$ computed with the explicit reverse-Euler inverse, so it also reflects time-integration accuracy; with the implicit fixed-point inverse (each backward Euler step solved by fixed-point iteration) the round-trip error is below $\mathrm{1e}\mathord{-}10$ on the tested samples of every run reported here and in \Cref{tab:genmod-cp-ablation}.

\Cref{tab:genmod-cp-ablation} first documents the sensitivity of the control-point CNF to the control-point count $m\le1024\approx10\%$ of $n_s=\inum{10000}$; unlike in \Cref{tab:cp-ablation}, the full expansion $m=n_s$ is not part of this study. On the checkerboard, the multiscale test NLL changes by $\mathrm{4e}\mathord{-}3$ (VarPro GMM) and $\mathrm{8e}\mathord{-}3$ (UGM) from $m=256$ to $m=1024$, with a non-monotone dip of $\mathrm{1.7e}\mathord{-}2$ for the UGM at $m=512$. The round-trip error $e_2$ is nearly unchanged. The reported $m=256$ fits lie within $\mathrm{1.7e}\mathord{-}2$ (UGM) and $\mathrm{4.4e}\mathord{-}3$ (GMM) of the best tested value, at a per-solve cost linear in $n_s$.

\begin{table}
\tabadjust
\caption{Sensitivity of the control-point CNF to the control-point count $m$ (multiscale checkerboard). We report the test NLL and round-trip error $e_2$ as the control-point count $m$ is varied at fixed data size $n_s=\inum{10000}$; $m=n_s$ would recover the full RKHS expansion~\cref{eq:rkhs-flow}, which is not part of this study. From $m=256$ on, the test NLL varies by at most $\mathrm{1.7e}\mathord{-}2$ (UGM) and $\mathrm{4.4e}\mathord{-}3$ (GMM).}
\label{tab:genmod-cp-ablation}
\begin{tabular}{rcccc}
\toprule
& \multicolumn{2}{c}{\bf UGM} & \multicolumn{2}{c}{\bf GMM} \\
\cmidrule(lr){2-3}\cmidrule(lr){4-5}
$m$ & \bf test nll & $e_2$ & \bf test nll & $e_2$ \\
\midrule
\inum{64}   & \fnum{3.5563} & \snum{8.072e-02} & \fnum{3.5871} & \snum{2.670e-03} \\
\inum{128}  & \fnum{3.5458} & \snum{6.379e-02} & \fnum{3.5559} & \snum{4.142e-03} \\
\inum{256}  & \fnum{3.5331} & \snum{6.529e-02} & \fnum{3.5444} & \snum{5.240e-03} \\
\inum{512}  & \fnum{3.5163} & \snum{6.901e-02} & \fnum{3.5419} & \snum{5.183e-03} \\
\inum{1024} & \fnum{3.5247} & \snum{6.551e-02} & \fnum{3.5400} & \snum{5.238e-03} \\
\bottomrule
\end{tabular}
\end{table}

\begin{table}
\tabadjust
\setlength{\tabcolsep}{3pt}
\caption{Concentric Circles test problem: complete control-point generative-modeling metrics ($m=256$). Each block reports one reference model (the UGM and the GMM eliminated by VarPro) for the single-scale (S) and multiscale (M) parameterization. From left to right we report the reference model, the parameterization, the NLL on the training and test data, the final objective, the reduced-gradient (single-scale) or proximal-residual (multiscale) norm in the $\ell_2$ and $\ell_\infty$ norms, the test-set reconstruction error in the $\ell_2$ and Frobenius norms ($x_i \mapsto z_i = \phi(x_i) \mapsto \hat{x}_i = \phi^{-1}(z_i)$), and the time to solution (in seconds). All runs are iterated to the 200-outer-step budget; the dimension-free $\ell_\infty$ norm is the reported convergence measure. The corresponding transport visualizations are shown in \Cref{fig:circles-genmod}.}
\label{tab:circles-genmod}
\begin{tabular}{llcccccccc}
\toprule
& & \multicolumn{2}{c}{\bf nll} & & \multicolumn{2}{c}{\bf gradient norm} & \multicolumn{2}{c}{\bf reconstruction error} & \\
\cmidrule(lr){3-4}\cmidrule(lr){6-7}\cmidrule(lr){8-9}
\bf reference & \bf scale & \bf train & \bf test & \bf objective & $\ell_2$ & $\ell_\infty$ & $e_2$ & $e_F$ & \bf time \\
\midrule
UGM & S & \fnum{3.282436} & \fnum{3.325295} & \fnum{695.014686} & \snum{2.217644e-01} & \snum{7.100547e-03} & \snum{6.087124e-03} & \snum{9.274555e-05} & \inum{404.725580} \\
    & M & \fnum{3.260455} & \fnum{3.298113} & \fnum{663.261782} & \snum{4.025392e-01} & \snum{8.416849e-03} & \snum{7.874217e-03} & \snum{1.168554e-04} & \inum{6320.499757} \\
\midrule
GMM & S & \fnum{3.260100} & \fnum{3.303512} & \fnum{659.127591} & \snum{1.306590e-01} & \snum{5.114718e-03} & \snum{9.933122e-04} & \snum{1.474960e-05} & \inum{736.764142} \\
    & M & \fnum{3.255153} & \fnum{3.296970} & \fnum{656.177509} & \snum{2.180860e-01} & \snum{3.355586e-03} & \snum{2.267168e-03} & \snum{3.133939e-05} & \inum{6046.077659} \\
\bottomrule
\end{tabular}
\end{table}

\begin{figure}
\centering
\includegraphics[width=\textwidth]{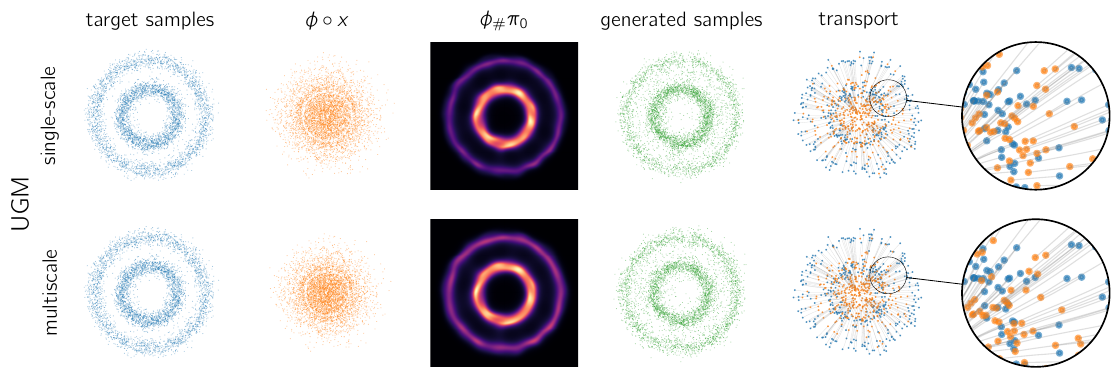}\\
\includegraphics[width=\textwidth]{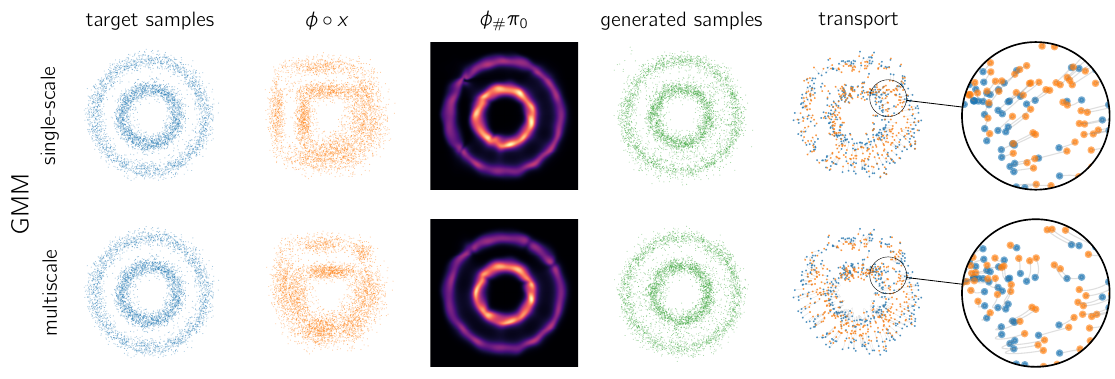}
\caption{Concentric circles test problem: transport visualizations for the two reference models. Top block: UGM; bottom block: GMM (VarPro). Within each panel the two rows compare the single-scale (top) and multiscale (bottom) parameterization, and the columns show (from left to right) the target distribution, the forward-mapped samples in latent space, the learned density in the data space, the generated samples obtained by sampling from the reference and applying the inverse map $\phi^{-1}$, and the particle transport (initial samples at time $t=0$ and the trajectories and samples at the final time $t=1$). The corresponding metrics are reported in \Cref{tab:circles-genmod}.}
\label{fig:circles-genmod}
\end{figure}

\begin{table}
\tabadjust
\setlength{\tabcolsep}{3pt}
\caption{M8G test problem: complete control-point generative-modeling metrics ($m=256$). Each block reports one reference model (the UGM and the GMM eliminated by VarPro) for the single-scale (S) and multiscale (M) parameterization. From left to right we report the reference model, the parameterization, the NLL on the training and test data, the final objective, the reduced-gradient (single-scale) or proximal-residual (multiscale) norm in the $\ell_2$ and $\ell_\infty$ norms, the test-set reconstruction error in the $\ell_2$ and Frobenius norms ($x_i \mapsto z_i = \phi(x_i) \mapsto \hat{x}_i = \phi^{-1}(z_i)$), and the time to solution (in seconds). All runs are iterated to the $200$-outer-step budget except the multiscale GMM fit, which stopped at 171 steps on the objective-stagnation rule of \Cref{app:varpro-variants}; the dimension-free $\ell_\infty$ norm is the reported convergence measure. The corresponding transport visualizations are shown in \Cref{fig:8gaussians-genmod}.}
\label{tab:8gaussians-genmod}
\begin{tabular}{llcccccccc}
\toprule
& & \multicolumn{2}{c}{\bf nll} & & \multicolumn{2}{c}{\bf gradient norm} & \multicolumn{2}{c}{\bf reconstruction error} & \\
\cmidrule(lr){3-4}\cmidrule(lr){6-7}\cmidrule(lr){8-9}
\bf reference & \bf scale & \bf train & \bf test & \bf objective & $\ell_2$ & $\ell_\infty$ & $e_2$ & $e_F$ & \bf time \\
\midrule
UGM & S & \fnum{2.818243} & \fnum{2.852773} & \fnum{635.022587} & \snum{1.772311e-01} & \snum{3.350664e-03} & \snum{1.419527e-02} & \snum{2.279021e-04} & \inum{868.322640} \\
    & M & \fnum{2.807694} & \fnum{2.837518} & \fnum{581.795983} & \snum{2.759640e-01} & \snum{5.431116e-03} & \snum{1.581085e-02} & \snum{2.573368e-04} & \inum{6065.835989} \\
\midrule
GMM & S & \fnum{2.818683} & \fnum{2.850242} & \fnum{564.009728} & \snum{5.293549e-02} & \snum{1.577868e-03} & \snum{2.378991e-05} & \snum{3.241178e-07} & \inum{523.961800} \\
    & M & \fnum{2.818989} & \fnum{2.850337} & \fnum{564.112479} & \snum{1.615249e-01} & \snum{2.621617e-03} & \snum{5.401352e-05} & \snum{7.252023e-07} & \inum{5069.625399} \\
\bottomrule
\end{tabular}
\end{table}

\begin{figure}
\centering
\includegraphics[width=\textwidth]{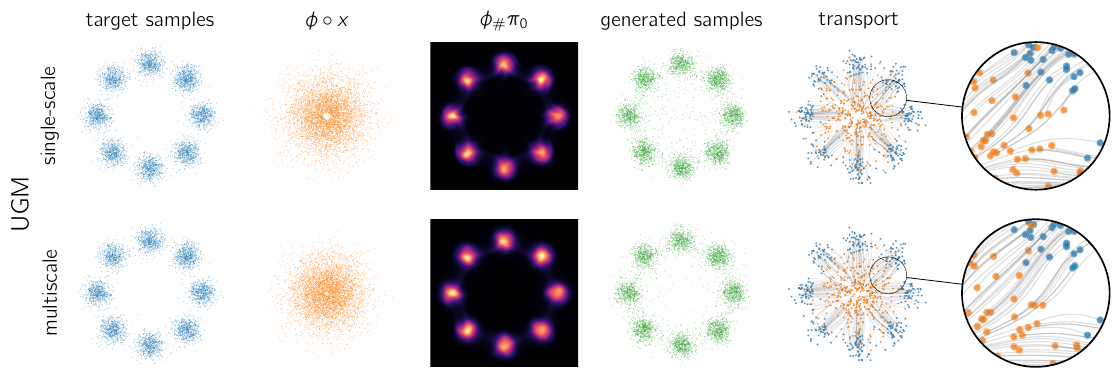}\\
\includegraphics[width=\textwidth]{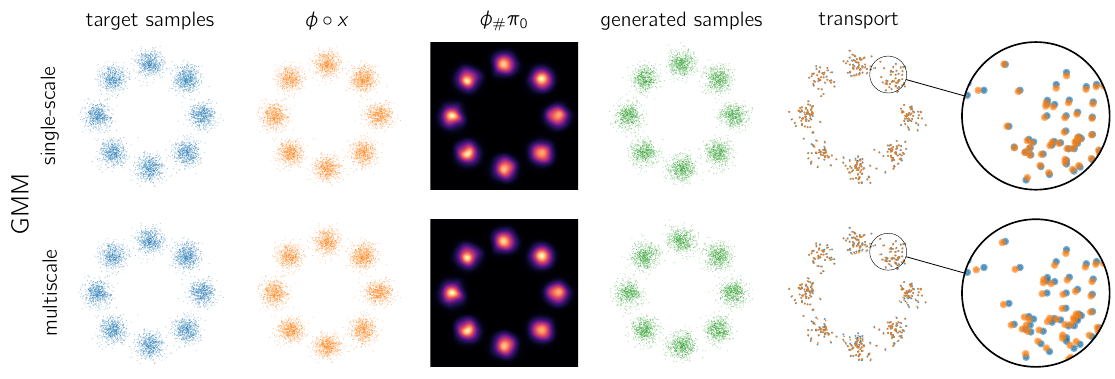}
\caption{M8G test problem: transport visualizations for the two reference models. Top block: UGM; bottom block: GMM. Within each panel the two rows compare the single-scale (top) and multiscale (bottom) parameterization, and the columns show (from left to right) the target distribution, the forward-mapped samples in latent space, the learned density in the data space, the generated samples obtained by sampling from the reference and applying the inverse map $\phi^{-1}$, and the particle transport (initial samples at time $t=0$ and the trajectories and samples at the final time $t=1$). The corresponding metrics are reported in \Cref{tab:8gaussians-genmod}.}
\label{fig:8gaussians-genmod}
\end{figure}

\begin{table}
\tabadjust
\setlength{\tabcolsep}{3pt}
\caption{Pinwheel test problem: complete control-point generative-modeling metrics ($m=256$). Each block reports one reference model (the UGM and the GMM) for the single-scale (S) and multiscale (M) parameterization. From left to right we report the reference model, the parameterization, the NLL on the training and test data, the final objective, the reduced-gradient (single-scale) or proximal-residual (multiscale) norm in the $\ell_2$ and $\ell_\infty$ norms, the test-set reconstruction error in the $\ell_2$ and Frobenius norms ($x_i \mapsto z_i = \phi(x_i) \mapsto \hat{x}_i = \phi^{-1}(z_i)$), and the time to solution (in seconds). All runs are iterated to the $200$-outer-step budget; the dimension-free $\ell_\infty$ norm is the reported convergence measure. The corresponding transport visualizations are shown in \Cref{fig:pinwheel-genmod}.}
\label{tab:pinwheel-genmod}
\begin{tabular}{llcccccccc}
\toprule
& & \multicolumn{2}{c}{\bf nll} & & \multicolumn{2}{c}{\bf gradient norm} & \multicolumn{2}{c}{\bf reconstruction error} & \\
\cmidrule(lr){3-4}\cmidrule(lr){6-7}\cmidrule(lr){8-9}
\bf reference & \bf scale & \bf train & \bf test & \bf objective & $\ell_2$ & $\ell_\infty$ & $e_2$ & $e_F$ & \bf time \\
\midrule
UGM & S & \fnum{2.347070} & \fnum{2.384624} & \fnum{495.998446} & \snum{2.804213e-01} & \snum{1.011390e-02} & \snum{8.162484e-03} & \snum{1.091919e-04} & \inum{646.893872} \\
    & M & \fnum{2.311491} & \fnum{2.347938} & \fnum{478.329493} & \snum{7.886548e-01} & \snum{1.321491e-02} & \snum{1.401059e-02} & \snum{1.842021e-04} & \inum{6159.024467} \\
\midrule
GMM & S & \fnum{2.323166} & \fnum{2.362617} & \fnum{469.612235} & \snum{4.334620e-01} & \snum{9.389634e-03} & \snum{1.291550e-03} & \snum{2.429312e-05} & \inum{783.292389} \\
    & M & \fnum{2.319794} & \fnum{2.354848} & \fnum{466.687110} & \snum{6.304373e-01} & \snum{1.039942e-02} & \snum{2.793385e-03} & \snum{4.397049e-05} & \inum{6003.647178} \\
\bottomrule
\end{tabular}
\end{table}

\begin{figure}
\centering
\includegraphics[width=\textwidth]{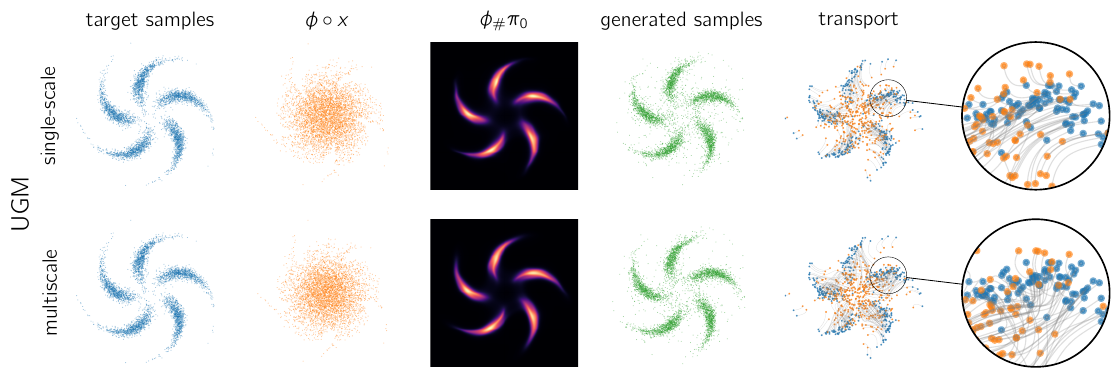}\\
\includegraphics[width=\textwidth]{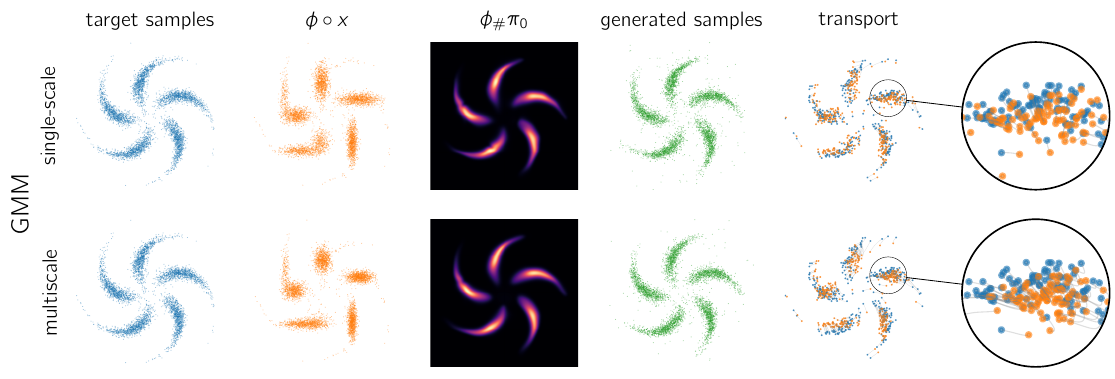}
\caption{Pinwheel test problem: transport visualizations for the two reference models. Top block: UGM; bottom block: GMM. Within each panel the two rows compare the single-scale (top) and multiscale (bottom) parameterization, and the columns show (from left to right) the target distribution, the forward-mapped samples in latent space, the learned density in the data space, the generated samples obtained by sampling from the reference and applying the inverse map $\phi^{-1}$, and the particle transport (initial samples at time $t=0$ and the trajectories and samples at the final time $t=1$). The corresponding metrics are reported in \Cref{tab:pinwheel-genmod}.}
\label{fig:pinwheel-genmod}
\end{figure}

\begin{table}
\tabadjust
\setlength{\tabcolsep}{3pt}
\caption{Checkerboard test problem: complete control-point generative-modeling metrics ($m=256$). Each block reports one reference model (the UGM and the GMM eliminated by VarPro) for the single-scale (S) and multiscale (M) parameterization. From left to right we report the reference model, the parameterization, the NLL on the training and test data, the final objective, the reduced-gradient (single-scale) or proximal-residual (multiscale) norm in the $\ell_2$ and $\ell_\infty$ norms, the test-set reconstruction error in the $\ell_2$ and Frobenius norms ($x_i \mapsto z_i = \phi(x_i) \mapsto \hat{x}_i = \phi^{-1}(z_i)$), and the time to solution (in seconds). All runs are iterated to the $200$-outer-step budget; the dimension-free $\ell_\infty$ norm is the reported convergence measure. The corresponding transport visualizations (UGM and VarPro GMM) are shown in \Cref{fig:checkerboard-varpro-gmm} in the main text.}
\label{tab:checkerboard-genmod}
\begin{tabular}{llcccccccc}
\toprule
& & \multicolumn{2}{c}{\bf nll} & & \multicolumn{2}{c}{\bf gradient norm} & \multicolumn{2}{c}{\bf reconstruction error} & \\
\cmidrule(lr){3-4}\cmidrule(lr){6-7}\cmidrule(lr){8-9}
\bf reference & \bf scale & \bf train & \bf test & \bf objective & $\ell_2$ & $\ell_\infty$ & $e_2$ & $e_F$ & \bf time \\
\midrule
UGM & S & \fnum{3.515578} & \fnum{3.539086} & \fnum{833.088328} & \snum{2.843342e-01} & \snum{9.441848e-03} & \snum{7.831668e-02} & \snum{1.282813e-03} & \inum{581.982760} \\
    & M & \fnum{3.511773} & \fnum{3.533090} & \fnum{744.755502} & \snum{2.781825e+00} & \snum{5.463982e-02} & \snum{6.529109e-02} & \snum{1.004729e-03} & \inum{6066.399942} \\
\midrule
GMM & S & \fnum{3.528293} & \fnum{3.555684} & \fnum{716.719958} & \snum{1.263378e-01} & \snum{4.869220e-03} & \snum{2.085118e-03} & \snum{2.414373e-05} & \inum{754.746375} \\
    & M & \fnum{3.514202} & \fnum{3.544360} & \fnum{711.237460} & \snum{1.515749e-01} & \snum{2.330397e-03} & \snum{5.240001e-03} & \snum{7.068113e-05} & \inum{6004.114405} \\
\bottomrule
\end{tabular}
\end{table}

\end{document}